\documentclass[english,10pt,leqno]{amsart}
 \usepackage{graphicx}    
 \usepackage{amsmath}
\usepackage{amssymb}
\usepackage[T1]{fontenc}
\usepackage[latin1]{inputenc}
\usepackage{babel}
\usepackage[arrow, matrix, curve]{xy}
\usepackage{amsthm}
\usepackage{amsmath,amscd}
\usepackage{epic,eepic}
\numberwithin{equation}{section}

\newcommand{\OK}{\mathcal{O}_K}
\newcommand{\F}{\mathbb{F}}
\newcommand{\VF}{V_\mathbb{F}}
\newcommand{\GR}{\mathcal{GR}_{\VF,0}}
\newcommand{\Fpot}{\mathbb{F}[\hspace{-0.5mm}[u]\hspace{-0.5mm}]}
\newcommand{\latM}{\mathfrak{M}}
\newcommand{\Phiquot}{\langle\Phi(\latM)\rangle/u^e\,\latM}
\newcommand{\GRloc}{\GR^{\mathbf{v},{\rm{loc}}}}

\newcommand{\Fbarpot}{\bar{\mathbb{F}}[\hspace{-0.5mm}[u]\hspace{-0.5mm}]}

\newcommand{\shOX}{\mathcal{O}_X}
\newcommand{\B}{\mathcal{B}}
\newcommand{\A}{\mathcal{A}}
\newcommand{\calL}{\mathcal{L}}
\newcommand{\T}{\mathcal{T}}
\newcommand{\bfv}{\mathbf{v}}
\newcommand{\G}{\mathcal{G}_{V_\F}}

\newtheorem{theo}{Theorem}[section]
\newtheorem{lem}[theo]{Lemma}
\newtheorem{prop}[theo]{Proposition}
\newtheorem{cor}[theo]{Corollary}
\theoremstyle{remark}
\newtheorem{rem}[theo]{Remark}
\theoremstyle{definition}
\newtheorem{defn}[theo]{Definition}

\def\plot#1{}

\begin{document}

\title[structure of moduli of finite flat group schemes]{On the structure of some moduli spaces of finite flat group schemes}
\author{Eugen Hellmann}

\maketitle

\section{Introduction and Notations}
Let $p$ be an odd prime and $k$ a finite field of characteristic $p$. Let $W=W(k)$ be the ring of Witt vectors with
coefficients in $k$ and $K_0=W[\tfrac{1}{p}]$ its fraction field.
We consider a finite, totally ramified extension $K/K_0$ and denote by $e=[K:K_0]$ the degree of the extension.
Let us fix a uniformizer $\pi\in\OK$ with minimal polynomial $E(u)\in W[u]$ over $K_0$.
Further we fix an algebraic closure $\bar K$ of $K$.

Let $\F$ be a finite field of characteristic $p$ and $\rho: G_K\rightarrow GL(V_\F)$ a continuous representation of
the absolute Galois group $G_K={\rm{Gal}}(\bar K/K)$ of $K$ in a finite dimensional $\F$-vector space $V_\F$ whose 
dimension will be denoted by $d$.\\
This datum is equivalent to a finite commutative group scheme $\widetilde{\mathcal{G}}\rightarrow {\rm{Spec}}\ K$ with an operation
of $\F$: The $\bar K$-valued points become an $\F$-vector space with a natural action of $G_K$ and we want
$\widetilde{\mathcal{G}}(\bar K)$ and $V_\F$ to be isomorphic as $\F[G_K]$-modules.\\
If $\F'$ is a finite extension of $\F$, the representation $\rho$ induces a representation $\rho'$ on
$V_{\F'}=V_\F\otimes_\F \F'$.

By the construction in Kisin's article \cite{Kisin}, there is a projective $\F$-scheme $\GR$ whose $\F'$-valued points
parametrize the isomorphism classes of finite flat models of $V_{\F'}$, i.e. finite flat group schemes 
$\mathcal{G}\rightarrow {\rm{Spec}}\ \OK$ with an operation of $\F'$ such that the generic fiber of $\mathcal{G}$
is the $G_K$-representation on $V_{\F'}$ in the above sense.\\
Our aim is to analyze the structure of (some stratification of) $\GR$ in the case $d=2$ and $k=\F_p$.

First we recall some constructions from \cite{Kisin}, see also \cite{phimod}. 
We assume $k=\F_p$ to simplify the situation.\\
For each $n$ let $\pi_n\in \bar K$ be a $p^n$-th root of the uniformizer $\pi$ such that $\pi_n^p=\pi_{n-1}$ for all $n$.
Define $K_\infty=\bigcup_{n\geq 1} K(\pi_n)$ and denote by $G_{K\infty}={\rm{Gal}}(\bar K/K_\infty)$ the absolute Galois group 
of $K_\infty$.\\
For each algebraic extension $\F'$ of $\F$ we denote by $\phi:\F'((u))\rightarrow \F'((u))$ the homomorphism which takes
$u$ to its $p$-th power and which is the identity on the coefficients:
\[\phi(\sum_i a_i u^i)=\sum a_i u^{pi}.\]   
Denote by ${\rm{Mod}}^{\ \phi}_{/\F'((u))}$ the category of finite dimensional $\F'((u))$-modules $M$
together with a $\phi$-linear map $\Phi: M\rightarrow M$ such that the linearization
${\rm{id}}\otimes\Phi: \phi^*M\rightarrow M$ is an isomorphism. The morphism are $\F'((u))$-linear maps commuting with $\Phi$. 
By (\cite{Kisin} 1.2.6, Lemma 1.2.7), there is an equivalence of abelian categories
\[
{\rm{Mod}}^{\ \phi}_{/\F'((u))}\longleftrightarrow \left\{ 
{\begin{array}{*{20}c}  \text{continuous}\ G_{K\infty}\text{-representations}\\
\text{on finite dimensional $\F'$-vector spaces} \end{array}}
\right\}
\]
which preserves the dimensions and is compatible with finite base change $\F''/\F'$.
This is a version with coefficients of the equivalence of categories of Fontaine (cf. \cite{Fontaine}, A3).

Denote by $(M_\F,\Phi)$ the $d$-dimensional $\F((u))$-vector space with semi-linear endomorphism $\Phi$,
associated to the restriction of the Tate-twist $V_\F(-1)$ to $G_{K_\infty}$ under the above equivalence.
By the descriptions in \cite{Kisin}, the finite flat models $\mathcal{G}\rightarrow {\rm{Spec}}\ \OK$ of $V_\F$
correspond to $\Fpot$-lattices $\latM\subset M_\F$ satisfying $u^e\latM\subset \langle\Phi(\latM)\rangle\subset \latM$.
Here $\langle\Phi(\latM)\rangle=(\rm{id}\otimes\Phi)\phi^*\latM$ is the $\Fpot$-lattice in $M_\F$ generated by $\Phi(\latM)$.\\
Under this description the multiplicative group schemes correspond to the lattices $\latM$ such that 
$\langle\Phi(\latM)\rangle=\latM$ and the étale group schemes correspond to the lattices with 
$u^e\latM=\langle\Phi(\latM)\rangle$. These lattices will be called multiplicative resp. étale.

This construction is compatible with base change in the following sense.
Suppose $\latM_\F\subset M_\F$ is a $\Fpot$-lattice corresponding to a finite flat model $\mathcal{G}$ of $V_\F$. 
If $\F'$ is a finite extension of $\F$ with $n=[\F':\F]$, then the $\F'[\hspace{-0.5mm}[u]\hspace{-0.5mm}]$ lattice 
\[\latM_{\F'}=\latM\widehat{\otimes}_\F \F'\subset M_{\F'}=M_{\F}\widehat{\otimes}_\F\F'\]
corresponds to the finite flat model $\mathcal{G}'=\mathcal{G}\boxtimes_\F\F'$ of $V_{\F'}$.
Here the exterior tensor product $\mathcal{G}\boxtimes_\F\F'$ is the following group scheme:
Choose a $\F$-basis $e_1\dots e_n$ of $\F'$. Then $\mathcal{G}\boxtimes_\F\F'=\prod_{i=1}^n\mathcal{G}$ 
and $z\in \F'$ operates via the matrix $A\in GL_n(\F)$ describing the multiplication by $z$ on $\F'$ in the fixed
$\F$-basis.

The scheme $\GR$ is constructed as a closed subscheme of the affine Grassmannian $\rm{Grass}\ M_\F$ for $GL(M_\F)$
and its closed points are given by
\begin{equation}\label{GRpts}
\GR(\F')=\{\F'[\hspace{-0.5mm}[u]\hspace{-0.5mm}]\text{-lattices}\ \latM\subset M_{\F'} 
\mid u^e\latM\subset \langle\Phi(\latM)\rangle\subset \latM \}
\end{equation} 
for every finite extension $\F'$ of $\F$.\\
In the following we will forget about the Galois representation and finite flat group schemes and will consider
lattices. We will drop the condition $p\neq 2$. All results hold for arbitrary $p$,
except those using the interpretation of the closed points as finite flat group schemes. 
We will always assume that there exists a finite flat model for $V_\F$ at least after extending
scalars.

For each $\mathbb{Q}_p$-algebra embedding $\psi:K\rightarrow \bar K_0$ we now fix an integer $v_\psi\in\{0,\dots,d\}$.
Denote by $\mathbf{v}=(v_\psi)_\psi$ the collection of the $v_\psi$ and by $\mathbf{r}=\check{\mathbf{v}}$ the dual
partition, i.e. $r_i=\sharp\{\psi\mid v_\psi\geq i\}$.\\
Kisin constructs closed, reduced subschemes 
\[\GRloc\subset\GR\]
whose $\F'$-valued points are given by
\begin{equation}\label{GRlocpts}
\GRloc(\F')=\{\latM\in\GR(\F')\mid J(u|_{\Phiquot})\leq \mathbf{r}\}
\end{equation}
for a finite extension $\F'$ of $\F$ (cf. \cite{Kisin}, Prop. 2.4.6).
Here $J(u|_{\Phiquot})$ denotes the Jordan type of the nilpotent endomorphism on $\Phiquot$ induced by the multiplication
with $u$. Recall that for $d=2$
\begin{equation}\label{order}
(a_1,b_1)\leq (a_2,b_2)\Leftrightarrow \left\{ \begin{array}{*{20}c}a_1\leq a_2,\\ a_1+b_1=a_2+b_2\end{array}\right.
\end{equation}
for pairs $(a_i,b_i)\in\mathbb{Z}^2$ with $a_i\geq b_i$.
The local structure of $\GRloc$ is linked to the structure of the local models studied in \cite{lokmod}.
These schemes are named \emph{" closed Kisin varieties"} in \cite{phimod}.

Kisin conjectures in (\cite{Kisin} 2.4.16) that, if $\rm{End}_{\F[G_K]}(V_\F)=\F$, 
the connected components of $\GRloc$ are given by the 
open and closed subschemes on which both the rank of the maximal multiplicative subobject and the rank of the 
maximal étale quotient are fixed.
In (\cite{Kisin}, 2.5) he proves this conjecture in the case $d=2$, $k=\F_p$ and $v_\psi=1$ for all $\psi$.
For $d=2$ and $v_\psi=1$ for all $\psi$ this result is generalized by Imai to the case of arbitrary $k$ (see \cite{Imai}).
In this paper we want to analyze the situation in the case $k=\F_p$, $d=2$ but arbitrary $\bfv$.
It turns out that the conjecture is not true in general. Our main results are as follows.

For $(a,b)\in\mathbb{Z}^2$ with $a\geq b$, we introduce a locally closed subscheme of the affine Grassmannian
\[\mathcal{G}_{V_\F}(a,b)\subset {\rm{Grass}}\, M_\F,\]
with closed points the lattices $\latM$ such that the elementary divisors of $\langle\Phi(\latM)\rangle$
with respect to $\latM$ are given by $(a,b)$.
\begin{theo}
Assume that $(M_{\F'},\Phi)=(M_{\F}\widehat{\otimes}_\F\F',\Phi)$ is simple for all finite extensions $\F'$ of $\F$.\\
\noindent (i) If\ \ $\mathcal{G}_{V_{\F}}(a,b)\neq \emptyset$, there exists a finite extension $\F'$ of $\F$ such that
\[\mathcal{G}_{V_{\F'}}(a,b)=\mathcal{G}_{V_{\F}}(a,b)\otimes_\F\F'\cong \mathbb{A}_{\F'}^n\]
for $n=\lfloor\tfrac{a-b}{p+1}\rfloor$.\\
\noindent (ii) The scheme $\GRloc$ is geometrically connected and irreducible. There exists a finite extension 
$\F'$ of $\F$ such that $\GRloc\otimes_\F\F'$ is isomorphic to a Schubert variety in the affine Grassmannian 
for $GL(M_{\F'})$.\\
The dimension of $\GRloc$ is either $\lfloor\tfrac{r_1-r_2}{p+1}\rfloor$ or $\lfloor\tfrac{r_1-r_2}{p+1}\rfloor-1$. 
Here $r_i=\sharp\{\psi\mid v_\psi\geq i\}$.
\end{theo} 
In the treatment of the reducible case we consider the set $\mathcal{S}(\mathbf{v})$ of isomorphism classes 
$[M']$ of one dimensional objects in ${\rm{Mod}}^{\ \phi}_{/\bar\F((u))}$ which admit an $\Fbarpot$-lattice $\latM_{[M']}\subset M'$
such that $\langle\Phi(\latM_{[M']})\rangle=u^{e-r_1}\latM_{[M']}$.
We will define subschemes
\[X^{\mathbf{v}}_{[M']}\subset \GRloc\otimes_\F\bar\F.\]
A lattice defines a closed point of $X^{\mathbf{v}}_{[M']}$ if it admits a $\Phi$-stable subobject isomorphic to 
$\latM_{[M']}$. A lattice $\latM$ is called $\mathbf{v}$-\emph{ordinary} iff it defines a closed point of 
$X^{\mathbf{v}}_{[M']}$ for some $[M']\in\mathcal{S}(\mathbf{v})$. The subscheme of non-$\mathbf{v}$-ordinary points
will be denoted by $X^{\mathbf{v}}_0$. We will prove the following Theorem.
\begin{theo} Assume that $(M_{\F'},\Phi)=(M_\F\widehat{\otimes}_\F\F',\Phi)$ is reducible for some finite extension $\F'$ of $\F$.\\
\noindent (i) The subschemes $X^{\mathbf{v}}_0$ and $X^{\mathbf{v}}_{[M']}$ 
are open and closed in $\GRloc\otimes_\F\bar\F$
for all isomorphism classes $[M']\in\mathcal{S}(\mathbf{v})$.\\
\noindent (ii) The scheme $X_0^{\mathbf{v}}$ is connected.\\
\noindent (iii) For each $[M']\in\mathcal{S}(\mathbf{v})$ the scheme $X^{\mathbf{v}}_{[M']}$ is connected.
If it is non empty, it is either a single point, or isomorphic to $\mathbb{P}^1_{\bar\F}$.\\
\noindent (iv) There are at most two isomorphism classes $[M']\in\mathcal{S}(\bfv)$ such that $X^\bfv_{[M']}\neq \emptyset$.
\end{theo}
The structure of the subscheme $X_0^\bfv$ of non-$\mathbf{v}$-ordinary lattices is much more complicated
than in the absolutely simple case.
In general $X_0^\bfv$ has many irreducible components of varying dimensions. 
The main result concerning the irreducible components of $X_0^\bfv$ is the following theorem.
\begin{theo}
If $(M_{\bar\F},\Phi)$ is not isomorphic to the direct sum of two isomorphic one-dimensional $\phi$-modules, then the irreducible components of $X_0^\bfv$ are Schubert varieties.
Especially they are normal.
\end{theo}
Theorem $1.2$ proves a modified version of Kisin's conjecture in the case $k=\F_p$ and $d=2$, as follows.\\
For an integer $s$ denote by 
\[\mathcal{GR}_{V_\F,0}^{\mathbf{v},{\rm{loc}},s}\subset \GRloc\]
the open and closed subscheme where the rank of the maximal $\Phi$-stable subobject $\latM_1$, satisfying
$\langle\Phi(\latM_1)\rangle=u^{e-r_1}\latM_1$, is equal to $s$.
\begin{cor}
Assume $p\neq 2$ and let $\rho:G_K\rightarrow V_\F$ be any two dimensional continuous representation of $G_K$.
Assume that ${\rm{End}}_{\F'[G_K]}(V_{\F'})$ is a simple algebra for all finite extensions $\F'$ of $\F$.
Then $\mathcal{GR}_{V_\F,0}^{\bfv,{\rm{loc}},s}$ is geometrically connected for all $s$.
Furthermore \\
\noindent (i) If $s=1$ and ${\rm{End}}_{\F'[G_K]}(V_{\F'})=\F'$ for all finite extensions $\F'$ of $\F$, 
then $\mathcal{GR}_{V_\F,0}^{\bfv,{\rm{loc}},s}$ is either empty or a single point.\\
If $s=1$ and ${\rm{End}}_{\F'[G_K]}(V_{\F'})=M_2(\F')$ for some finite extension $\F'$ of $\F$,
then $\mathcal{GR}_{V_\F,0}^{\bfv,{\rm{loc}},s}$ is either empty or becomes isomorphic to $\mathbb{P}^1_{\F'}$ after extending the
scalars to $\F'$.\\
\noindent (ii) If $s=2$, then $\mathcal{GR}_{V_\F,0}^{\bfv,{\rm{loc}},s}$ is either empty or a single point.
\end{cor}
{\bf{Acknowledgments:}} This paper is the author's diploma thesis written at the University of Bonn.
I want to thank M. Rapoport for introducing me into this subject
and for many helpful discussions. 
I also want to thank X. Caruso for lots of explanations on Kisin and Breuil modules and for his interest in my work.
\section{Some notations in the building}
The method of this paper is to determine all lattices in the building of $GL_2\left(\bar\F((u))\right)$ that correspond
to closed points of $\GRloc$. As we know that the scheme we study is a closed reduced subscheme of the affine 
Grassmannian, we can get information on the structure of $\GRloc$ by looking at its closed points.\\
For the rest of this paper, we fix the following notations:
Let $(M_{\F},\Phi)$ be the object in ${\rm{Mod}}^{\ \phi}_{/\F((u))}$ corresponding to the $2$-dimensional Galois representation
$\rho$ on $V_\F$.
Let $\bfv=(v_\psi)_\psi$ be a collection of integers $v_\psi\in\{0,1,2\}$ for every $\psi:K\rightarrow \bar K_0$.
Define 
\begin{equation}\label{dstrich}
d'=\sum_{\psi}v_\psi.
\end{equation}
Denote by $\mathbf{r}=\check{\bfv}$ the dual partition, i.e.
\begin{align*}
& r_1=\sharp\{\psi\mid v_\psi\geq 1\}\\
& r_2=\sharp\{\psi\mid v_\psi\geq 2\}.
\end{align*}
Denote by $\B$ the Bruhat-Tits building for $GL_2\left(\F((u))\right)$. For any finite extension $\F'$ of $\F$
the building for $GL_2\left(\F'((u))\right)$ will be denoted by $\B_{\F'}$. We write
\[\bar\B=\bigcup_{\F'/\F}\B_{\F'}\]
for the building for $GL_2\left(\bar\F((u))\right)$.\\
We choose an $\F((u))$-basis $e_1,e_2$ of $M_\F$. Denote by $\latM_0=\langle e_1,e_2\rangle$ the standard lattice in the 
standard apartment $\A_0$ determined by $e_1,e_2$. In this apartment we choose the following coordinates:\\
Let $(m,n)_0$ denote the lattice $\langle u^me_1,u^ne_2\rangle$. Further, we consider another set of coordinates given by 
$[x,y]_0=(\tfrac{x+y}{2},\tfrac{y-x}{2})_0$ for $x,y\in\mathbb{Z},\ x\equiv y\mod 2$; i.e. $(m,n)_0=[m-n,m+n]_0$.\\
Let $q\in\F((u))^\times$ and set $k=v_u(q)\in\mathbb{Z}$, where $v_u$ is the valuation on $\F((u))$ with $v_u(u)=1$.
The basis $e_1,qe_1+e_2$ of $M_\F$ defines another apartment $\A_q$ which is branching off from the standard apartment
at the line defined by $x=k$. Using the Iwasawa decomposition we find
\[\B=\bigcup_{q\in\F((u))}\A_q.\]
For arbitrary $q\in\F((u))$ we choose coordinates in the apartments $\A_q$, similar to the case of $\A_0$.
Define 
\[(m,n)_q=[m-n,m+n]_q:=\langle u^me_1,u^n(qe_1+e_2)\rangle\in\A_q.\] 
\begin{rem}\label{coord}
\noindent \emph{(i)} The systems of coordinates in the various apartments are compatible in the following sense:
For any $x,y\in\mathbb{Z}\ ,x\equiv y\mod 2$ and $q,q'\in\F((u))$ we have $[x,y]_q=[x,y]_{q'}$ iff $x\leq v_u(q-q')$ which implies
\[[x,y]_q=[x,y]_{q'}\Leftrightarrow [x,y]_q\in \A_q\cap\A_{q'}\Leftrightarrow [x,y]_{q'}\in\A_q\cap\A_{q'}.\]
\noindent \emph{(ii)} We will make use of these coordinates for arbitrary points in the building (not only points
corresponding to lattices). We see that $[x,y]_q$ defines a lattice if and only if $x,y\in\mathbb{Z}$ and 
$x\equiv y \mod 2$.\\
\noindent \emph{(iii)} We extend the above notations in the obvious way to the buildings $\bar\B$ and $\B_{\F'}$
for arbitrary finite extensions $\F'$ of $\F$.\\
\noindent \emph{(iv)} Two points $[x,y]_q,[x',y']_q\in\A_q$ define the same point in the building for 
$PGL_2\left(\F((u))\right)$ if and only if $x=x'$. Thus the projection from $\B$ onto the building for 
$PGL_2\left(\F((u))\right)$ is given by the projection onto the $x$-coordinate for every apartment $\A_q\subset\B$.
\end{rem}
\begin{defn}\label{defd1d2}
Let $\latM$ and $\latM'$ be lattices in $M_\F$. Let $a,b$ be the elementary divisors of $\latM'$ with respect to $\latM$,
i.e. there exists a basis $e_1',e_2'$ of $\latM$ such that $\latM'=\langle u^ae_1',u^be_2'\rangle$.
Define
\begin{align*}
&\ d_1(\latM,\latM')=|a-b|\\
&\ d_2(\latM,\latM')=a+b.
\end{align*}
\end{defn}
\begin{rem}\label{dbuilding}
These quantities have the following meaning in the building:\\
If $\latM=[x,y]_q$ and $\latM'=[x',y']_{q'}$, then $d_2(\latM,\latM')=y'-y$. Further $d_1(\latM,\latM')$
is the distance between $\latM$ and $\latM'$ in the building for $PGL_2(\F((u)))$. Here, the distance between two lattices
joined by an edge is equal to $1$. This can be seen as follows:
Assume that $\latM=\langle e_1,e_2\rangle$ is the standard lattice and $\latM'=A\latM=[x,y]_q$, with
\[A=(a_{ij})_{ij}=\begin{pmatrix}u^m & u^nq\\ 0& u^n\end{pmatrix}.\]
If $a\geq b$ are the elementary divisors of $\latM'$ with respect to $\latM$, then, by the theory of elementary divisors,
\[d_1(\latM,\latM')=a+b-2b=m+n-2\min_{i,j}v_u(a_{ij}).\]
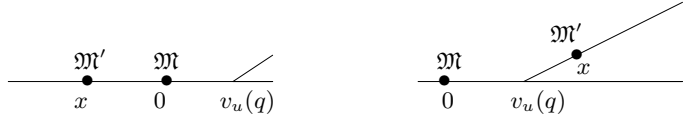
\begin{figure}[h]
\begin{picture}(350,25)
\drawline(55,0)(155,0)(140,0)(155,10) \jput(85,0){\makebox(0,0){$\bullet$}} \jput(115,0){\makebox(0,0){$\bullet$}}
\jput(80,5){{\small $\latM'$}} \jput(110,5){{\small $\latM$}}
\jput(80,-10){{\small $x$}} \jput(110,-10){{\small $0$}} \jput(135,-10){{\small $v_u(q)$}}

\drawline(210,0)(310,0)(250,0)(310,30) 
\jput(220,0){\makebox(0,0){$\bullet$}}  \jput(270,10){\makebox(0,0){$\bullet$}}
\jput(217,5){{\small $\latM$}} \jput(245,-10){{\small $v_u(q)$}} \jput(260,15){{\small $\latM'$}}
\jput(219,-10){{\small $0$}} \jput(270,3){{\small $x$}}
\end{picture}
\vspace{2mm}
\caption{The distance between two lattices in the building 
for $PGL_2(\F((u)))$ in the cases $x\leq v_u(q)$ and $x\geq v_u(q)\geq 0$.}
\end{figure}\\
If $x=m-n\leq v_u(q)$ or $v_u(q)\geq 0$, then $\min_{i,j}v_u(a_{ij})=\min\{m,n\}$ and hence 
\[d_1(\latM,\latM')=a-b=|m-n|=|x|.\]
If $x> v_u(q)$ and $v_u(q)<0$, then $\min_{i,j}v_u(a_{ij})=n+v_u(q)$. In this case we find 
\[d_1(\latM,\latM')=a-b=m-n-2v_u(q)=(x-v_u(q))+(0-v_u(q)).\] 
Compare also Fig. $1$ and Fig. $2$.\\
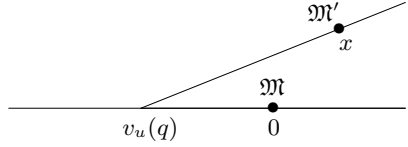
\begin{figure}[h]
\begin{picture}(350,25)
\drawline(90,0)(240,0)(140,0)(240,40)
\jput(190,0){\makebox(0,0){$\bullet$}} \jput(215,30){\makebox(0,0){$\bullet$}}
\jput(133,-10){{\small $v_u(q)$}} \jput(185,5){{\small $\latM$}} \jput(203,33){{\small $\latM'$}}
\jput(188,-10){{\small $0$}} \jput(215,22){{\small $x$}}
\end{picture}
\vspace{2mm}
\caption{The distance between two lattices in the building for $PGL_2(\F((u)))$
in the case $x>v_u(q)$ and $v_u(q)<0$.}
\end{figure}\\
We see that the distance $d_1(\latM,\latM')$ only depends on $x,x'$ (and on $v_u(q-q')$),
while $d_2(\latM,\latM')$ only depends on $y$ and $y'$.\\
Using this remark, we can extend the distances $d_1$ and $d_2$ in an obvious way to the whole building $\B$
(and to $\bar\B,\,\B_{\F'}$). For example 
\[d_1([x,y]_q,[0,0]_0)=
\begin{cases}\ \ x& \text{if}\ x\geq 0,\ v_u(q)\geq 0\\ \ -x& \text{if}\ x<0,\ x<v_u(q)\\
 \ x-2v_u(q)& \text{if}\ v_u(q)<x,\ v_u(q)<0. \end{cases}\]
\end{rem}
\begin{lem}\label{lemdist}
Define $d'$ as in $(\ref{dstrich})$.
The closed points $z\in\GRloc(\bar\F)$ correspond to the lattices $\latM\subset M_{\bar\F}$ which satisfy
\begin{align*}
&\ d_1(\latM,\langle\Phi(\latM)\rangle)\leq r_1-r_2\\
&\ d_2(\latM,\langle\Phi(\latM)\rangle)=2e-d'.
\end{align*}
\end{lem}
\begin{proof}
If $\latM\subset M_{\bar\F}$ is any lattice and if $a,b$ are the elementary divisors of $\langle\Phi(\latM)\rangle$
with respect to $\latM$, then the above conditions read
\begin{align*}
& a-b\leq r_1-r_2\\
& a+b=2e-d'=2e-(r_1+r_2).
\end{align*} 
This implies $u^e\latM\subset \langle\Phi(\latM)\rangle\subset\latM$.
The Jordan type of $u$ on $\Phiquot$ is given by 
\[J(u|_{\Phiquot})=(e-a,e-b).\]
Assuming $a\geq b$ we find:
\[J(u|_{\Phiquot})\leq\mathbf{r}\Leftrightarrow \begin{cases}
b\geq e-r_1 \\ a+b=2e-d'=2e-(r_1+r_2).
\end{cases}\]
The lemma follows easily from this.
\end{proof}
\begin{defn}\label{vadm}
A lattice $\latM$ is called $\mathbf{v}$-\emph{admissible} if it satisfies
\[d_1(\latM,\langle\Phi(\latM)\rangle)\leq r_1-r_2\hspace{5mm}\text{and}
\hspace{5mm}d_2(\latM,\langle\Phi(\latM)\rangle)=2e-d'.\]
\end{defn}
Let $\latM$ be a lattice in $M_\F$ and $A\in GL_2\left(\F((u))\right)$ be a matrix. We will use the notation $\latM\sim A$
if $\latM$ admits a $\Fpot$-basis $b_1,b_2$ satisfying $\Phi(b_i)=Ab_i$.
Similarly we will use the notation $M_\F\sim A$ (use a $\F((u))$-basis of $M_\F$).
\begin{lem}\label{ygleich}
For $i=1,2$, let $z_i\in\GRloc(\bar\F)$ be closed points corresponding to lattices $\latM_i=[x_i,y_i]_{q_i}\in\bar\B$.
Then $y_1=y_2$.
\end{lem}
\begin{proof}
Choose $A,B\in GL_2\left(\bar\F((u))\right)$ such that $\latM_2=A\latM_1$ and $\latM_1\sim B$. 
Then $\latM_2\sim \phi(A)BA^{-1}$. Using the theory of elementary divisors it follows that
\[v_u(\det\,B)=d_2(\latM_i,\langle\Phi(\latM_i)\rangle)=(p-1)v_u(\det\,A)+v_u(\det\,B)\]
and hence $v_u(\det\,A)=0$ which yields the claim.
\end{proof}
\begin{defn}\label{subtree}
For each $m\in\mathbb{Z}$ define the following subset of $\bar\B$:
\[\bar\B(m):=\bigcup_{q\in\bar\F((u))}\{[x,y]_q\in\A_q\mid y=m\}.\]
\end{defn}
Viewing $\GRloc(\bar\F)$ as a subset of $\bar\B$, Lemma $\ref{ygleich}$ implies: 
\[\GRloc(\bar\F)\subset \bar\B(m)\]
for some $m=m(\bfv)\in\mathbb{Z}$. The subset $\bar\B(m)$ is a tree which is (as a topological space) isomorphic to
the building for $PGL_2(\bar\F((u)))$. \\
The difference is that not every vertex represents a lattice:
A vertex $[x,m]_q\in\bar\B(m)$ represents a lattice $\latM\subset M_{\bar\F}$ iff $x\equiv m\mod 2$.
\begin{rem}\label{subgrass}
By construction we have
\[\GR\subset {\rm{Grass}}\,M_\F,\]
where $\rm{Grass}\,M_\F$ denotes the affine Grassmannian for $GL_2$. Since the determinant condition in $(\ref{GRlocpts})$ 
fixes the dimension 
\[\dim\,\Phiquot=\sum_\psi v_\psi=d',\]
the closed subscheme $\GRloc$ lies in a connected component of this Grassmannian:
If $\latM=A\latM_0$ defines a closed point (where $\latM_0$ is the standard lattice and $A$ is a matrix),
then the valuation of $\det\, A$ is determined by the dimension of $\Phiquot$.
\end{rem}
\begin{defn}
For a given collection $\bfv$ denote by $\GR^\bfv$ the closed subscheme of $\GR$ consisting of all lattices $\latM$
such that $\dim \Phiquot=\sum_\psi v_\psi$.
\end{defn}
\begin{prop} If any two of the $v_\psi$ differ at most by $1$, then 
\[\GR^\bfv=\GRloc.\]
\end{prop}
\begin{proof}
Let $\latM$ be a lattice and $a\geq b$ the elementary divisors of $\langle\Phi(\latM)\rangle$ with respect to $\latM$.
Then $\latM$ defines a closed point of $\GR^\bfv$ if and only if $0\leq b\leq a\leq e$ and $a+b=2e-d'$ with 
$d'=\sum_\psi v_\psi$. This is equivalent to 
\begin{align*}
& a+b=2e-d'\\
& a-b\leq \max\{d',2e-d'\}.
\end{align*}
Indeed, if $a+b=2e-d'$, then $a-b\leq d'$ is equivalent to $a\leq e$, while $a-b\leq 2e-d'$ is equivalent to $b\geq 0$.\\
Now if any of the $v_\psi$ differ at most by $1$ we must have $r_1=e$ (if all $v_\psi\geq 1$) or $r_2=0$ 
(if all $v_\psi\leq 1$). 
In both cases we find $r_1-r_2=\max\{d',2e-d'\}$.\\
(See also \cite{Kisin}, Prop. 2.4.6 (4) and \cite{lokmod}, Thm. B \emph{(iii)}).
\end{proof}
\section{The absolutely simple case}
In this section we will analyze the structure of $\GRloc$ in the case where $(M_\F,\Phi)$ is absolutely simple,
i.e. for every (finite) extension $\F'/\F$ there is no proper $\Phi$-stable subobject of $(M_{\F'},\Phi)$.
\begin{lem}\label{NFirred}
If $(M_\F,\Phi)$ is absolutely simple, there exists a finite extension $\F'$ of $\F$, a basis $e_1,e_2$ of $M_{\F'}$
and $a\in \F'^\times,\ s\in\mathbb{Z}$ satisfying 
\[0\leq s<p^2-1\ \text{and}\ s\not\equiv 0\mod (p+1)\] such that
\[M_{\F'}\sim\begin{pmatrix}0&au^s\\1&0\end{pmatrix}.\]
\end{lem}
\begin{proof}
This follows from (\cite{Caruso}, Cor. $8$), except that we need to check that $s\not\equiv 0\mod (p+1)$. If $p+1|s$, then there would be a proper $\Phi$-stable subspace
of $M_{\F''}$ for a quadratic extension $\F''$ of $\F'$, namely $\langle \sqrt{a}u^{s/p+1}e_1+e_2\rangle\subset M_{\F'}\widehat{\otimes}_{\F'}\F'[\sqrt{a}]$.\\
The constructions in \cite{Caruso} give a basis after extending scalars to the algebraic closure $\bar\F$
of $\F$, but of course this also gives a basis after finite field extension, as there are only finitely many equations
to solve. See also (\cite{Imai}, Lemma 1.2).
\end{proof}
For the rest of this section we fix the basis $e_1,e_2$ of Lemma $\ref{NFirred}$ as the standard basis of $M_{\bar\F}$
and use the coordinates introduced in section 2. Furthermore we fix the point
\begin{equation}\label{ptP}
P_{\rm{irred}}:=[\tfrac{s}{p+1},-\tfrac{s}{p-1}]_0\in\A_0\subset\bar\B.
\end{equation}
\begin{prop}\label{lemirred}
\noindent (i) The map $\Phi$ extends to a map $\bar\B\rightarrow \bar\B$ also denoted by $\Phi$.\\
\noindent (ii) Let $[x,y]_0\in\A_0$ be any point in the standard apartment. Then 
\[\Phi([x,y]_0)=[-px+s,py+s]_0.\]
\noindent (iii) For any $q\in\bar\F((u))^\times$ with $k=v_u(q)$ and $[x,y]\in\A_q\backslash\A_0$, i.e. $x>k$,
the map $\Phi$ is given by 
\[\Phi([x,y]_q)=[px-2pk+s,py+s]_{q'}\in\A_{q'}\]
for some $q'\in\bar\F((u))^\times$ with $v_u(q')=-pk+s\neq k$.  \\
\noindent (iv) The point $P_{\rm{irred}}$, as defined in $(\ref{ptP})$, satisfies $\Phi(P_{\rm{irred}})=P_{\rm{irred}}$.\\
\noindent (v) If $Q\in\B$ is an arbitrary point, then
\begin{align*}
&\ d_1(Q,\Phi(Q))=(p+1)d_1(Q,P_{\rm{irred}})\\
&\ d_2(Q,\Phi(Q))=(p-1)d_2(Q,P_{\rm{irred}}).
\end{align*}
\end{prop}
\begin{proof}
\noindent \emph{(i)} We can use the expressions in \emph{(ii)} and \emph{(iii)} to extend $\Phi$.\\
\noindent \emph{(ii)} We have 
\begin{align*}
&\ \Phi(u^me_1)=u^{pm}\Phi(e_1)=u^{pm}e_2\\
&\ \Phi(u^ne_2)=u^{pn}\Phi(e_2)=au^{pn+s}e_1
\end{align*}
and hence $\Phi((m,n)_0)=(pn+s,pm)_0$. The statement follows from this.\\
\emph{(iii)} We put $v_u(q)=k$ and $\phi(q)=\alpha u^{pk}$ for some $\alpha\in\Fbarpot^{\times}$. If $\latM=(m,n)_q$, then
\[\langle\Phi(\latM)\rangle=\langle u^{pm}e_2\ ,\ u^{pn}\phi(q)e_2+au^{pn+s}e_1\rangle.\]
As $\latM=[m-n,m+n]_q\notin\A_0$ we have $m> n+k$. Hence
\begin{align*}
\langle\Phi(\latM)\rangle&=\langle u^{pm}e_2-\alpha^{-1}u^{p(m-n-k)}(u^{pn}\phi(q)e_2+au^{pn+s}e_1),
u^{pn}\phi(q)e_2+au^{pn+s}e_1\rangle\\
&=\langle u^{p(m-k)+s}e_1, u^{p(n+k)}(q'e_1+e_2)\rangle
\end{align*}
with $q'=\alpha^{-1}au^{-pk+s}$. And thus $\Phi((m,n)_q)=(p(m-k)+s,p(n+k))_{q'}$ with $v_u(q')=-pk+s\neq k$
, as $k\not\equiv 0\mod (p+1)$.\\
\noindent \emph{(iv)} Obvious.\\
\noindent \emph{(v)}
If $\latM=[x,y]_q$, then the statement on $d_2$ follows immediately from \emph{(ii)} and \emph{(iii)}:
\[d_2(\latM,\langle\Phi(\latM)\rangle)=(p-1)y+s=(p-1)d_2([x,y]_q,[\tfrac{s}{p+1},-\tfrac{s}{p-1}]_0).\]
\begin{figure}[h]
\begin{picture}(350,25)
\drawline(70,0)(220,0)(100,0)(130,40)
\drawline(190,0)(220,40)
\jput(98,-10){{\small $k$}} \jput(170,-10){{\small $-p\,k+s$}}
\jput(120,26){\makebox(0,0){$\bullet$}} \jput(200,13){\makebox(0,0){$\bullet$}}
\jput(123,21){{\small $x$}} \jput(105,27){{\small $\latM$}}
\jput(130,-10){{\small $s/p+1$}} \jput(135,0){$|$}
\jput(205,10){{\small $p\,x-2p\,k+s$}} \jput(165,17){{\small $\langle\Phi(\latM)\rangle$}}
\end{picture}
\vspace{3mm}
\caption{The images of $\latM$ and $\Phi(\latM)$ in the building for $PGL_2(\bar\F((u)))$
in the case $\latM\notin\A_0$.}
\end{figure}
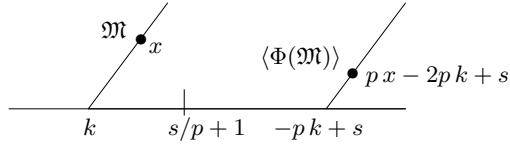\\
For the statement on $d_1$ first assume that $\latM=[x,y]_0\in\A_0$. Then \emph{(ii)} implies
\[d_1(\latM,\langle\Phi(\latM)\rangle)=|(p+1)x-s|=(p+1)d_1([x,y]_0,[\tfrac{s}{p+1},-\tfrac{s}{p-1}]_0).\] 
If $\latM=[x,y]_q\in\A_q\backslash \A_0$, then $x>k$ and $px-2pk+s>-pk+s$ which implies $\langle\Phi(\latM)\rangle\notin\A_0$.
Now \emph{(iii)} implies
\begin{align*}
d_1(\latM,\langle\Phi(\latM)\rangle)&=(x-k)+|-pk+s-k|+(px-2pk+s-(-pk+s))\\
&= (p+1)(x-k)+(p+1)\left|k-\tfrac{s}{p+1}\right|\\ &=
\begin{cases} (p+1)(x-\tfrac{s}{p+1}) & \text{if}\ k>\tfrac{s}{p+1}\\
(p+1)(x-k+\tfrac{s}{p+1}-k) & \text{if}\ k<\tfrac{s}{p+1},
\end{cases}
\end{align*}
using $k\neq -pk+s$. In both cases the claim follows. (See also Fig. 3)
\end{proof}
\begin{rem}
The Proposition shows that the absolutely simple case is exactly the case discussed in (\cite{phimod}, 6.d, A1).
The fixed point in the building is the point $P_{\rm{irred}}$ and its projection onto the building for $PGL_2(\bar\F((u)))$
lies on the edge between two vertices. The set of lattices $\latM$ with $d_1(\latM,\langle\Phi(\latM)\rangle)\leq r_1-r_2$
is identified with a ball around this fixed point.  
\end{rem}
Let $\bfv$ be a collection of integers as in the introduction. By Lemma $\ref{lemdist}$ and Proposition 
$\ref{lemirred}$\emph{(v)} we find
$\GRloc(\bar\F)\subset \bar\B(m(\bfv))$, where 
\begin{equation}\label{mvirred}
m(\bfv)=(2e-d'-s)/(p-1).
\end{equation}
\begin{cor}\label{GRempty}
The scheme $\GRloc$ is empty if $2e-d'\not\equiv s\mod (p-1)$.
\end{cor}
\begin{proof}
This follows from Lemma $\ref{lemdist}$ and Proposition $\ref{lemirred}$.
\end{proof}
\begin{figure}[h]
\begin{picture}(400,100)
\drawline(100,0)(270,0)
\jput(160,0){\makebox(0,0){$\bullet$}}
\jput(210,0){\makebox(0,0){$\bullet$}}
\jput(210,25){\makebox(0,0){$\bullet$}}
\jput(210,-25){\makebox(0,0){$\bullet$}}
\jput(260,0){\makebox(0,0){$\bullet$}}
\jput(110,0){\makebox(0,0){$\bullet$}}
\drawline(260,0)(270,10)
\drawline(260,0)(270,-10)
\drawline(235,0)(260,25) \drawline(260,25)(270,30) \drawline(260,25)(270,35) \drawline(260,25)(270,40)
			\jput(260,25){\makebox(0,0){$\bullet$}}

\drawline(210,0)(260,50)
\drawline(235,25)(260,40) \drawline(235,25)(260,60) \jput(260,40){\makebox(0,0){$\bullet$}} 
\jput(260,50){\makebox(0,0){$\bullet$}} \jput(260,60){\makebox(0,0){$\bullet$}}
\drawline(260,40)(270,42) \drawline(260,40)(270,46) \drawline(260,40)(270,50)
\drawline(260,50)(270,53) \drawline(260,50)(270,56) \drawline(260,50)(270,60)
\drawline(260,60)(270,64) \drawline(260,60)(270,68) \drawline(260,60)(270,72)

\drawline(185,0)(210,25) \drawline(210,25)(235,50) \drawline(210,25)(235,60) \drawline(210,25)(235,70) 
\drawline(235,50)(260,70) \drawline(235,50)(260,75) \drawline(235,50)(260,80)
\drawline(235,60)(260,85) \drawline(235,60)(260,90) \drawline(235,60)(260,95)
\drawline(235,70)(260,100) \drawline(235,70)(260,105) \drawline(235,70)(260,110)
\jput(260,70){\makebox(0,0){$\bullet$}}\jput(260,75){\makebox(0,0){$\bullet$}} \jput(260,80){\makebox(0,0){$\bullet$}}
\jput(260,85){\makebox(0,0){$\bullet$}} \jput(260,90){\makebox(0,0){$\bullet$}} \jput(260,95){\makebox(0,0){$\bullet$}}
\jput(260,100){\makebox(0,0){$\bullet$}} \jput(260,105){\makebox(0,0){$\bullet$}} \jput(260,110){\makebox(0,0){$\bullet$}}
\drawline(260,80)(270,81) \drawline(260,80)(270,82) \drawline(260,80)(270,83)
\drawline(260,85)(270,86) \drawline(260,85)(270,87) \drawline(260,85)(270,88)
\drawline(260,90)(270,91) \drawline(260,90)(270,92) \drawline(260,90)(270,93)
\drawline(260,70)(270,71) \drawline(260,70)(270,72) \drawline(260,70)(270,73)
\drawline(260,75)(270,76) \drawline(260,75)(270,77) \drawline(260,75)(270,78)
\drawline(260,95)(270,96) \drawline(260,95)(270,97) \drawline(260,95)(270,98)
\drawline(260,100)(270,101) \drawline(260,100)(270,102) \drawline(260,100)(270,103)
\drawline(260,105)(270,106) \drawline(260,105)(270,107) \drawline(260,105)(270,108)
\drawline(260,110)(270,111) \drawline(260,110)(270,112) \drawline(260,110)(270,113)

\drawline(160,0)(110,50)
\drawline(135,25)(110,40) \drawline(135,25)(110,60)
\drawline(135,0)(110,25)
\jput(110,40){\makebox(0,0){$\bullet$}} \jput(110,50){\makebox(0,0){$\bullet$}} \jput(110,60){\makebox(0,0){$\bullet$}}
\jput(110,25){\makebox(0,0){$\bullet$}}
\drawline(110,40)(100,42) \drawline(110,40)(100,46) \drawline(110,40)(100,50)
\drawline(110,50)(100,53) \drawline(110,50)(100,53) \drawline(110,50)(100,60)
\drawline(110,60)(100,64) \drawline(110,60)(100,68) \drawline(110,60)(100,72)
\drawline(110,25)(100,30) \drawline(110,25)(100,35) \drawline(110,25)(100,40)
\drawline(110,0)(100,10)

\drawline(235,0)(260,-25) \drawline(260,-25)(270,-30) \drawline(260,-25)(270,-35) \drawline(260,-25)(270,-40)
			\jput(260,-25){\makebox(0,0){$\bullet$}}

\drawline(210,0)(260,-50)
\drawline(235,-25)(260,-40) \drawline(235,-25)(260,-60) \jput(260,-40){\makebox(0,0){$\bullet$}} 
\jput(260,-50){\makebox(0,0){$\bullet$}} \jput(260,-60){\makebox(0,0){$\bullet$}}
\drawline(260,-40)(270,-42) \drawline(260,-40)(270,-46) \drawline(260,-40)(270,-50)
\drawline(260,-50)(270,-53) \drawline(260,-50)(270,-56) \drawline(260,-50)(270,-60)
\drawline(260,-60)(270,-64) \drawline(260,-60)(270,-68) \drawline(260,-60)(270,-72)

\drawline(185,0)(210,-25) \drawline(210,-25)(235,-50) \drawline(210,-25)(235,-60) \drawline(210,-25)(235,-70) 
\drawline(235,-50)(260,-70) \drawline(235,-50)(260,-75) \drawline(235,-50)(260,-80)
\drawline(235,-60)(260,-85) \drawline(235,-60)(260,-90) \drawline(235,-60)(260,-95)
\drawline(235,-70)(260,-100) \drawline(235,-70)(260,-105) \drawline(235,-70)(260,-110)
\jput(260,-70){\makebox(0,0){$\bullet$}}\jput(260,-75){\makebox(0,0){$\bullet$}} \jput(260,-80){\makebox(0,0){$\bullet$}}
\jput(260,-85){\makebox(0,0){$\bullet$}} \jput(260,-90){\makebox(0,0){$\bullet$}} \jput(260,-95){\makebox(0,0){$\bullet$}}
\jput(260,-100){\makebox(0,0){$\bullet$}} \jput(260,-105){\makebox(0,0){$\bullet$}} \jput(260,-110){\makebox(0,0){$\bullet$}}
\drawline(260,-80)(270,-81) \drawline(260,-80)(270,-82) \drawline(260,-80)(270,-83)
\drawline(260,-85)(270,-86) \drawline(260,-85)(270,-87) \drawline(260,-85)(270,-88)
\drawline(260,-90)(270,-91) \drawline(260,-90)(270,-92) \drawline(260,-90)(270,-93)
\drawline(260,-70)(270,-71) \drawline(260,-70)(270,-72) \drawline(260,-70)(270,-73)
\drawline(260,-75)(270,-76) \drawline(260,-75)(270,-77) \drawline(260,-75)(270,-78)
\drawline(260,-95)(270,-96) \drawline(260,-95)(270,-97) \drawline(260,-95)(270,-98)
\drawline(260,-100)(270,-101) \drawline(260,-100)(270,-102) \drawline(260,-100)(270,-103)
\drawline(260,-105)(270,-106) \drawline(260,-105)(270,-107) \drawline(260,-105)(270,-108)
\drawline(260,-110)(270,-111) \drawline(260,-110)(270,-112) \drawline(260,-110)(270,-113)

\drawline(160,0)(110,-50)
\drawline(135,-25)(110,-40) \drawline(135,-25)(110,-60)
\drawline(135,0)(110,-25)
\jput(110,-40){\makebox(0,0){$\bullet$}} \jput(110,-50){\makebox(0,0){$\bullet$}} \jput(110,-60){\makebox(0,0){$\bullet$}}
\jput(110,-25){\makebox(0,0){$\bullet$}}
\drawline(110,-40)(100,-42) \drawline(110,-40)(100,-46) \drawline(110,-40)(100,-50)
\drawline(110,-50)(100,-53) \drawline(110,-50)(100,-53) \drawline(110,-50)(100,-60)
\drawline(110,-60)(100,-64) \drawline(110,-60)(100,-68) \drawline(110,-60)(100,-72)
\drawline(110,-25)(100,-30) \drawline(110,-25)(100,-35) \drawline(110,-25)(100,-40)
\drawline(110,0)(100,-10)
\jput(170,2){\makebox(0,0){$|$}}
\jput(154,-12){{\tiny $[\tfrac{s}{p+1},m(\bfv)]_0$}}
\end{picture}
\vspace{40mm}
\caption{This picture illustrates the subset of $\bfv$-admissible lattices in the case $p=3$ and $\F=\F_3$.
This subset is given by all lattices $\latM\in \bar\B(m(\bfv))$ satisfying 
$d_1(\latM,P_{\rm{irred}})\leq (r_1-r_2)/(p+1)$. The fat points correspond to 
$\bfv$-admissible lattices.}
\end{figure}
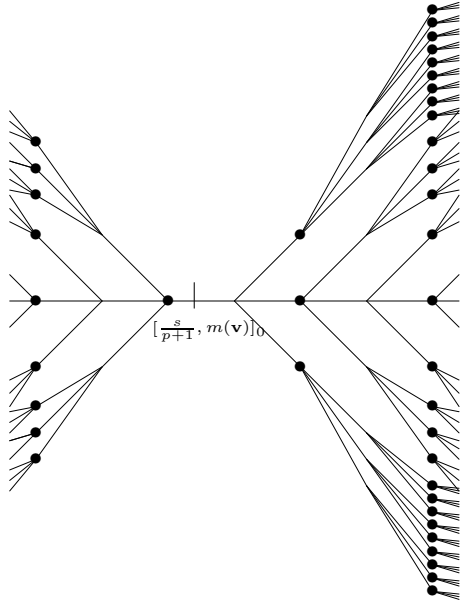
Now we want to define locally closed subschemes of ${\rm{Grass}}\,M_{\F}$ on which the elementary divisors of 
$\langle\Phi(\latM)\rangle$ with respect to $\latM$ are fixed.
Define a function
\[E:{\rm{Grass}}\,M_\F\rightarrow \mathbb{Z}^2.\]
For an extension field $L$ of $\F$ and an $L$-valued point $z\in ({\rm{Grass}}\,M_\F)(L)$ consider the 
$\Fpot\widehat{\otimes}_\F L$-lattice $\latM_z$ in $M_\F\widehat{\otimes}_\F L$ corresponding to $z$.
Then $E(z)=(j_1,j_2)$, where $j_1\geq j_2$ are the elementary divisors of $\langle\Phi(\latM_z)\rangle$
with respect to $\latM_z$. 
Recall that there is a partial order on the pairs $(a,b)\in\mathbb{Z}^2$ given by $(\ref{order})$.
\begin{lem}\label{lowersemcont}
The function $E$ is lower semi-continuous with respect to the Zariski topology on ${\rm{Grass}}\,M_{\F}$.
\end{lem}
\begin{proof}
Let $\eta\rightsquigarrow z$ be a specialization and let $\latM_\eta$ and $\latM_z$ be the lattices corresponding to the points $\eta$ and $z$. Denote by $E(\eta)=(a(\eta),b(\eta))$ and $E(z)=(a(z),b(z))$ the elementary divisors of $\langle\Phi(\latM_\eta)\rangle$ with respect to $\latM_\eta$ (resp. the elementary divisors of $\langle\Phi(\latM_z)\rangle$ with respect to $\latM_z$). 
We mark the specialization by a morphism 
$f:{\rm{Spec}}\,R\rightarrow\GRloc$, where $R$ is a discrete valuation ring with uniformizer $t$. The morphism $f$
defines a $R[\hspace{-0.5mm}[u]\hspace{-0.5mm}]$-lattice $\latM_R$ in $M_\F\widehat{\otimes}_\F R$.
After choosing a basis we find a matrix $C=(c_{ij})_{ij}\in GL_2(R((u)))\cap M_2(R[\hspace{-0.5mm}[u]\hspace{-0.5mm}])$  
such that $\latM_R\sim C$. Denote by $\bar c_{ij}$ the reduction $\mod t$ of the matrix coefficients. Using the theory of elementary divisors we find
\[b(\eta)=\min_{i,j}v_u(c_{ij})\leq \min_{i,j}v_u(\bar c_{ij})=b(z)\]
and hence $E(\eta)\geq E(z)$ which yields the claim. 
\end{proof}
\begin{defn}\label{Gab}
Let $(a,b)\in\mathbb{Z}^2$ such that $a\geq b$. The \emph{"Kisin variety"} associated to $(a,b)$ is 
\[\G(a,b)=E^{-1}(a,b)\subset {\rm{Grass}}\,M_{\F}.\]
By Lemma $\ref{lowersemcont}$, this is a locally closed subset 
and it will be considered as a subscheme with the reduced scheme structure (See also \cite{phimod}).
\end{defn}
Now we want to analyze the structure of $\G(a,b)$ and $\GRloc$. We will make use of the following fact.
\begin{lem}\label{P1nachgrass}
Let $b_1,b_2$ be any basis of $M_{\bar\F}$. There exists a morphism
\[\chi:\mathbb{A}^1_{\bar\F}\rightarrow {\rm{Grass}}\,M_{\bar\F}\]
such that $\chi(z)=\langle b_1,zu^{-1}b_1+b_2\rangle$ for every closed point $z\in\mathbb{A}^1_{\bar\F}$.
The morphism $\chi$ extends in a unique way to a morphism
\[\bar\chi:\mathbb{P}^1_{\bar\F}\rightarrow {\rm{Grass}}\,M_{\bar\F}.\]
The image of the point at infinity is given by $\bar\chi(\infty)=\langle u^{-1}b_1,ub_2\rangle$.
\end{lem}
\begin{proof}
Consider the family  
\[\langle b_1,Tu^{-1}b_1+b_2\rangle_{\bar\F[T][\hspace{-0.5mm}[u]\hspace{-0.5mm}]}
\subset M_{\bar\F}\widehat{\otimes}_{\bar\F}\bar\F[T]\]
of lattices on $\mathbb{A}^1={\rm{Spec}}\,\bar\F[T]$. This family defines the morphism $\chi$.\\ 
Let $X$ be the closed subscheme of ${\rm{Grass}}\,M_{\bar\F}$ consisting of all lattices $\latM$ that satisfy
$u\langle b_1,b_2\rangle\subset \latM\subset u^{-1}\langle b_1,b_2\rangle$ and that lie in the same connected component
of ${\rm{Grass}}\,M_{\bar\F}$ as $\langle b_1,b_2\rangle$. The scheme $X$ is identified with a closed subscheme of the
(ordinary) Grassmann variety $\rm{Grass}_{\bar\F}(4,2)$ of $2$-dimensional subspaces in $\bar\F^4$.
The morphism $\chi$ factors as follows:
\[\begin{xy}
\xymatrix{
\mathbb{A}^1_{\bar\F} \ar[rr]^{\chi}\ar@{.>}[rd]^{\chi'}& & **[r]{\rm{Grass}}\,M_{\bar\F} \\
& X \ar@^{(->}[ru] \ar[r] &  \rm{Grass}_{\bar\F}(4,2) 
\ar[r]^/3mm/{\iota}  & \mathbb{P}^5_{\bar\F}
}
\end{xy}\]
where $\iota$ is the Pl\"ucker embedding. As $X$ is projective, the valuative criterion shows that $\chi$ extends 
in a unique way to $\mathbb{P}^1$.\\
We view $\rm{Grass}_{\bar\F}(4,2)$ as the quotient $GL_{2,\bar\F}\backslash V$, where $V$ is the scheme of $2\times 4$
matrices of rank $2$ and $GL_{2,\bar\F}$ acts on $V$ by left multiplication.
Now, the computations using Plücker coordinates gives
\[\chi'(z)=\begin{pmatrix} 0&1&0&0\\ z& 0&0&1\end{pmatrix}\ ,\ \iota(\chi'(z))=(-z:0:0:0:1:0)\]
for all closed points $z\in\mathbb{A}^1(\bar\F)$. Hence the extension to $\mathbb{P}^1$ is 
\[(z_1:z_2)\mapsto (-z_1:0:0:0:z_2:0).\]
The image of the point at infinity is
\[\bar\chi((1:0))=\begin{pmatrix}0&1&0&0\\ 1&0&0&0\end{pmatrix}\mapsto (-1:0:0:0:0:0).\]
This is the lattice $\langle u^{-1}b_1,ub_2\rangle$.
\end{proof}
\begin{rem}\label{P1inbuilding}
In the building, the $\F$-valued points of the image of the morphism $\bar\chi$ can be illustrated in the following way
(if the morphism is defined over $\F$):
\end{rem}
\begin{figure}[h]
\vspace{5mm}
\begin{picture}(300,20)
\drawline(100,0)(200,0)
\jput(125,0){\makebox(0,0){$\bullet$}} \jput(175,0){\makebox(0,0){$\bullet$}}
\drawline(125,0)(175,0)
\drawline(150,0)(180,24) \drawline(150,0)(180,36)
\drawline(150,0)(180,-24) \drawline(150,0)(180,-36)
\jput(175,20){\makebox(0,0){$\bullet$}} \jput(175,30){\makebox(0,0){$\bullet$}}
\jput(175,-20){\makebox(0,0){$\bullet$}} \jput(175,-30){\makebox(0,0){$\bullet$}}
\end{picture}
\vspace{12mm}\\
\caption{The morphism $\bar\chi$ in the building for $p=5$ and $\F=\F_5$.}
\end{figure}
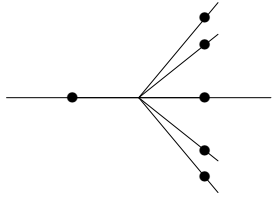
Similarly, we can define morphisms $\bar\chi_1,\bar\chi_2:\mathbb{P}^1_{\bar\F}\rightarrow {\rm{Grass}}\,M_{\bar\F}$ such that
\begin{align*}
&\ {\rm{im}}(\bar\chi_1)=\{\langle u^{n-1}b_1,u^{-(n-1)}(zu^{-1}b_1+b_2)\rangle\mid z\in\bar\F\}
\cup\{\langle u^{-n}b_1,u^nb_2\rangle\}\\
&\ {\rm{im}}(\bar\chi_2)=\{\langle u^nb_1,u^{-n}(zb_1+b_2)\rangle\mid z\in\bar\F\}
\cup\{\langle u^{-n}b_1,u^nb_2\rangle\}
\end{align*}
\begin{theo}\label{theoirred}
Assume that $(M_\F,\Phi)$ is absolutely simple. Fix a finite extension $\F'$ of $\F$ such that the normal form for $\Phi$
of Lemma $\ref{NFirred}$ is defined over $\F'$.\\ 
\noindent (a) For any $(a,b)\in\mathbb{Z}^2$ with $a\geq b$:
\[\G(a,b)\neq \emptyset\Leftrightarrow a+b\equiv s\mod(p-1)\ ,\ 
\begin{cases} & pa+b\equiv s\ \mod  (p^2-1)\\ \ \text{or} & pa+b\equiv ps\mod (p^2-1). \end{cases}\]
This condition being satisfied, there exists an isomorphism 
\[\G(a,b)\otimes_\F\F'\cong \mathbb{A}^n_{\F'}\]
with $n=\lfloor\tfrac{a-b}{p+1}\rfloor$. Further
\[\overline{\G(a,b)}=\bigcup_{(a',b')\leq (a,b)}\G(a',b').\]
\noindent (b) The scheme $\GRloc$ is geometrically connected and irreducible. After extending the scalars to $\F'$
it becomes isomorphic to a Schubert variety in the affine Grassmannian for $M_{\F'}=M_\F\widehat{\otimes}_\F\F'$
with dimension given by
\[\dim\,\GRloc=
\left\lfloor\tfrac{r_1-r_2}{p+1}-(-1)^\epsilon\tfrac{s}{p+1}\right\rfloor
+\left\lfloor(-1)^\epsilon\tfrac{s}{p+1}\right\rfloor,\]
with $\epsilon=\lfloor\tfrac{r_1-r_2}{p+1}\rfloor+\lfloor\tfrac{s}{p+1}\rfloor+\tfrac{2e-d'-s}{p-1}$.\\
\rm{[Here as in the rest of the paper, $\lfloor x \rfloor$ denotes the integral part of a real number $x$.]}
\end{theo}
\begin{proof}
\noindent {\emph{(a)}}
Assume $\latM\in\G(a,b)\neq \emptyset$. Without loss of generality, we may assume $\latM=[x,y]_0\in\A_0$:
if $\latM$ is an arbitrary lattice, then there exists a lattice $\latM'\in\A_0$ such that 
$d_i(\latM,P_{\rm{irred}})=d_i(\latM',P_{\rm{irred}})$ for $i=1,2$ (compare Fig. 4, for example).\\ 
By Lemma $\ref{lemirred}$\emph{(v)} and Definition $\ref{defd1d2}$, the condition for $\latM=[x,y]_0\in\G(a,b)$ is
\[(p+1)d_1(\latM,P_{\rm{irred}})=a-b\hspace{5mm}\hspace{5mm}(p-1)d_2(\latM,P_{\rm{irred}})=a+b.\]
By an explicit computation of these distances, this is equivalent to
\[\left|x-\tfrac{s}{p+1}\right|=\tfrac{a-b}{p+1}\hspace{5mm}\hspace{5mm}y+\tfrac{s}{p-1}=\tfrac{a+b}{p-1}.\]
The second equation gives $s\equiv a+b\mod (p-1)$ and the sum of both equations gives $s\equiv pa+b\mod (p^2-1)$
if $(p+1)x>s$ and $ps\equiv pa+b\mod (p^2-1)$ if $(p+1)x<s$ (using the fact that $x+y$ and $x-y$ are even).\\
Conversely, suppose $s\equiv a+b\mod (p-1)$ and $s\equiv pa+b\mod (p^2-1)$ and define 
\[x=\tfrac{a-b+s}{p+1}\hspace{10mm}y=\tfrac{a+b-s}{p-1}.\]
Then we have $y\in\mathbb{Z}$ and $x+y\in 2\mathbb{Z}$. Thus $[x,y]_0$ defines a lattice $\latM\in\G(a,b)$.
If $ps\equiv pa+b\mod (p^2-1)$ we use
\[x=\tfrac{s-(a-b)}{p+1}\hspace{10mm}y=\tfrac{a+b-s}{p-1}.\]
Now fix the sum $a+b$ and denote by $y$ the integer solving the equation 
\[(p-1)y+s=a+b.\]
Let us assume that $x_0:=\lfloor\tfrac{s}{p+1}\rfloor\equiv y\mod 2$
(the case $x_0\not\equiv y\mod 2$ admits a similar treatment).
In this case $[x_0,y]_0$ defines a lattice $\latM_0$ 
and we denote by $X$ the connected component of ${\rm{Grass}}\,M_{\F'}$ containing $\latM_0$, i.e. $X(\bar\F)=\{\latM\in\bar\B(y)\}$.\\
For each $m\geq 0$, there is a morphism $f_m:\mathbb{A}^{2m+1}_{\F'}\rightarrow X$ given by the family of lattices
\begin{align*}
\left\langle u^{(x_0+y)/2}u^{m+1}e_1\ ,\ 
u^{(y-x_0)/2}u^{-(m+1)}\left((\sum_{i=1}^{2m+1} T_{i}u^{i+x_0})e_1+e_2\right)\right\rangle\\
\subset M_{\F'}\widehat{\otimes}_{\F'}(\F'[T_1,\dots,T_{2m+1}])
\end{align*}
on $\mathbb{A}^{2m+1}_{\F'}={\rm{Spec}}\, \F'[T_1,\dots,T_{2m+1}]$. Let $V_m\cong \mathbb{A}^{2m+1}_{\F'}$ be its image. We have
\begin{equation}\label{affodd}
\begin{aligned}
&\bar\B(y)\supset V_m(\bar\F)=\\
&\{\langle u^{(x_0+y)/2}u^{m+1}e_1\,,\,u^{(y-x_0)/2}u^{-(m+1)}(qe_1+e_2)\rangle\mid 
q=\sum_{i=1}^{2m+1}a_iu^{i+x_0}\},
\end{aligned}
\end{equation}
with $a_1\dots a_{2m+1}\in\bar\F$.
Similarly, define for $m\geq 0$ a morphism $g_m:\mathbb{A}^{2m}_{\F'}\rightarrow X$ given by the family of lattices 
\begin{align*}
\left\langle u^{(x_0+y)/2}u^{-m}\left(e_1+(\sum_{i=0}^{2m-1} T_iu^{i-x_0}) e_2\right)\ ,\ u^{(y-x_0)/2}u^me_2\right\rangle \\
\subset M_{\F'}\widehat{\otimes}_{\F'}(\F'[T_0,\dots,T_{2m-1}])
\end{align*}
and let $U_m\cong \mathbb{A}^{2m}_{\F'}$ be its image. We have
\begin{equation}\label{affeven}
\begin{aligned}
&\bar\B(y)\supset U_m(\bar\F)=\\
&\{\langle u^{(x_0+y)/2}u^{-m}(e_1+qe_2)\,,\,u^{(y-x_0)/2}u^me_2\rangle\mid q=\sum_{i=0}^{2m-1}a_iu^{i-x_0}\},
\end{aligned}
\end{equation}
with $a_0\dots a_{2m-1}\in\bar\F$. 
It is easy to see that every lattice $\latM\in\bar\B(y)$ is either of the form $(\ref{affodd})$ or of the form $(\ref{affeven})$ 
for some $m\geq 0$. Thus
\[X=(\bigcup_{m\geq 0}V_m)\cup (\bigcup_{m\geq 0}U_m).\]
We claim
\begin{align}\label{dodd}
&\ \latM\in V_m(\bar\F)\Longrightarrow d_1(\latM,P_{\rm{irred}})=2m+2-\xi \\
&\ \latM\in U_m(\bar\F)\Longrightarrow d_1(\latM,P_{\rm{irred}})=2m+\xi,\label{deven}
\end{align}
where $\xi=\tfrac{s}{p+1}-x_0$ denotes the fractional part of $\tfrac{s}{p+1}$.\\
Indeed, if $\latM\in V_m(\bar\F)$, then $\latM=[x_0+2m+2,y]_q$ for some $q\in\bar\F((u))$ with $v_u(q)>x_0$
and hence 
\[d_1(\latM,P_{\rm{irred}})=x_0+2m+2-\tfrac{s}{p+1}=2m+2-\xi.\]
The statement on $U_m$ follows by a more complicated computation or by a symmetry argument:
The choice of apartments $\A_q$ and coordinates $[-,-]_q$ depends on the order of $e_1$ and $e_2$.
Interchanging $e_1$ and $e_2$ yields expressions for the lattices $\latM\in U_m$ similar to the above expressions for $V_m$
(if $\latM\in U_m$ is a lattice, then $\latM=[-x_0+2m,y]_q$ for some $q$) while it maps the point $P_{\rm{irred}}$ 
to $[-\tfrac{s}{p+1},-\tfrac{s}{p-1}]_0$ and hence the claim follows by the same computation. \\
Now equation $(\ref{dodd})$ and $(\ref{deven})$  together with Proposition $\ref{lemirred}$ \emph{(v)} imply 
\begin{equation}\label{subodd}
\begin{aligned}
&\ V_m(\bar\F)\subset \G(a_{\rm{odd}}(m),b_{\rm{odd}}(m))(\bar\F)\\
&\ U_m(\bar\F)\subset \G(a_{\rm{even}}(m),b_{\rm{even}}(m))(\bar\F)
\end{aligned}
\end{equation}  
for some $(a_{\rm{odd}}(m),b_{\rm{odd}}(m))\,,\,(a_{\rm{even}}(m),b_{\rm{even}}(m))\in\mathbb{Z}^2$ with 
\begin{align}
&\ a_{\rm{odd}}(m)+b_{\rm{odd}}(m) && =a_{\rm{even}}(m')+b_{\rm{even}}(m')=(p-1)y+s \notag\\
&\ a_{\rm{odd}}(m)-b_{\rm{odd}}(m) && =(p+1)(2m+2-\xi) \label{diffodd}\\
&\ a_{\rm{even}}(m)-b_{\rm{even}}(m) && =(p+1)(2m+\xi)\notag
\end{align}
and $0<\xi<1$ implies that all these pairs are pairwise distinct when $m$ runs over all positive integers.\\
As $U_m$ and $V_m$ cover $X$, the inclusions in $(\ref{subodd})$ are actually equalities. Furthermore 
$V_m=\G(a_{\rm{odd}}(m),b_{\rm{odd}}(m))$ as schemes, as both are reduced locally closed subschemes of ${\rm{Grass}}\,M_{\F'}$
with the same underlying point set. Finally $(\ref{diffodd})$ yields
\[\dim\,V_m=2m+1=\lfloor2m+2-\xi\rfloor=\lfloor\tfrac{a_{\rm{odd}}(m)-b_{\rm{odd}}(m)}{p+1}\rfloor.\]
The conclusion for $U_m$ is similar.\\
To finish the proof of \noindent {\emph{(a)}}, it remains to show that $U_m\subset \overline{V_m}$ and $V_{m-1}\subset \overline{U_m}$.
We will prove the first assertion: the second is proved in the same way.\\
Let $z_1\in U_m$ be an arbitrary point corresponding to a lattice 
\[\latM_1=\langle u^{(x_0+y)/2}u^{-m}(e_1+qe_2)\,,\,u^{(y-x_0)/2}u^me_2\rangle\]
with $q=\sum_{i=0}^{2m-1}a_iu^{i-x_0}$ and let $z_2\in V_m$ be the point corresponding to 
\[\latM_2=\langle u^{(x_0+y)/2}u^{m+1}e_1\,,\,u^{(y-x_0)/2}u^{-(m+1)}e_2\rangle.\]
There exists a basis $b_1$ and $b_2$ of $M_{\bar\F}$ such that 
\begin{align*}
&\ \langle b_1,b_2\rangle=\latM_0=[x_0,y]_0,\\
&\ \langle u^{-m}b_1,u^mb_2\rangle=\latM_1,\\
&\ \langle u^{m+1}b_1,u^{-(m+1)}b_2\rangle=\latM_2.
\end{align*}
Explicitly, we may choose 
\[b_1=u^{(x_0+y)/2}(e_1+qe_2)\ ,\ b_2=u^{(y-x_0)/2}e_2.\]
Applying Lemma $\ref{P1nachgrass}$ (resp. Remark $\ref{P1inbuilding}$) with the basis $ub_1,u^{-1}b_2$, we obtain a morphism 
$\chi:\mathbb{A}^1_{\bar\F}\rightarrow {\rm{Grass}}\,M_{\bar\F}$ that is given by 
$\chi(z)=\langle{u^{m+1}b_1,u^{-(m+1)}(zub_1+b_2)}\rangle$ on closed points
and we easily find ${\rm{im}}\,\chi\subset V_m\otimes_{\F'}\bar\F$. As $\overline{V_m}\otimes_{\F'}{\bar\F}$ is projective,
the morphism $\chi$ extends to a morphism from $\mathbb{P}^1$ to $\overline{V_m}\otimes_{\F'}{\bar\F}$ 
and the point at infinity is mapped to $z_1$
(Fig. 6 illustrates the image of the morphism $\bar\chi$ in the building. 
The fat points are the lattices in the image of $\bar\chi$).
Hence $z_1\in\overline{V_m}(\bar\F)$ and the claim follows.\\
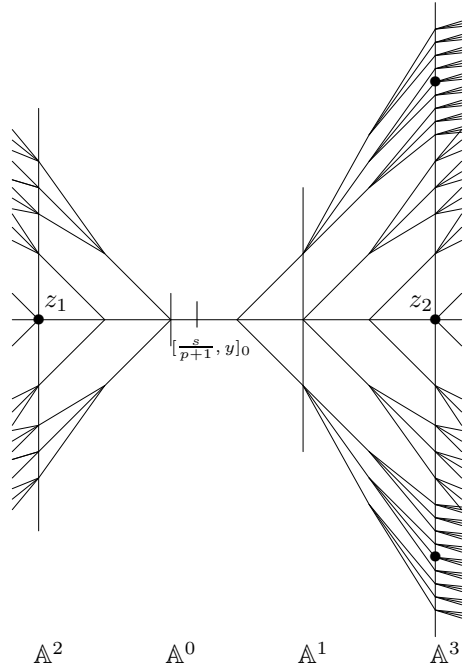
\begin{figure}[h]
\begin{picture}(350,100)
\drawline(100,0)(270,0)
\jput(260,0){\makebox(0,0){$\bullet$}} \jput(250,5){$z_2$}
\jput(110,0){\makebox(0,0){$\bullet$}} \jput(112,5){$z_1$}
\drawline(260,0)(270,10)
\drawline(260,0)(270,-10)
\drawline(235,0)(260,25) \drawline(260,25)(270,30) \drawline(260,25)(270,35) \drawline(260,25)(270,40)

\drawline(210,0)(260,50)
\drawline(235,25)(260,40) \drawline(235,25)(260,60) 
\drawline(260,40)(270,42) \drawline(260,40)(270,46) \drawline(260,40)(270,50)
\drawline(260,50)(270,53) \drawline(260,50)(270,56) \drawline(260,50)(270,60)
\drawline(260,60)(270,64) \drawline(260,60)(270,68) \drawline(260,60)(270,72)

\drawline(185,0)(210,25) \drawline(210,25)(235,50) \drawline(210,25)(235,60) \drawline(210,25)(235,70) 
\drawline(235,50)(260,70) \drawline(235,50)(260,75) \drawline(235,50)(260,80)
\drawline(235,60)(260,85) \drawline(235,60)(260,90) \drawline(235,60)(260,95)
\drawline(235,70)(260,100) \drawline(235,70)(260,105) \drawline(235,70)(260,110)
\jput(260,90){\makebox(0,0){$\bullet$}}
\drawline(260,80)(270,81) \drawline(260,80)(270,82) \drawline(260,80)(270,83)
\drawline(260,85)(270,86) \drawline(260,85)(270,87) \drawline(260,85)(270,88)
\drawline(260,90)(270,91) \drawline(260,90)(270,92) \drawline(260,90)(270,93)
\drawline(260,70)(270,71) \drawline(260,70)(270,72) \drawline(260,70)(270,73)
\drawline(260,75)(270,76) \drawline(260,75)(270,77) \drawline(260,75)(270,78)
\drawline(260,95)(270,96) \drawline(260,95)(270,97) \drawline(260,95)(270,98)
\drawline(260,100)(270,101) \drawline(260,100)(270,102) \drawline(260,100)(270,103)
\drawline(260,105)(270,106) \drawline(260,105)(270,107) \drawline(260,105)(270,108)
\drawline(260,110)(270,111) \drawline(260,110)(270,112) \drawline(260,110)(270,113)

\drawline(160,0)(110,50)
\drawline(135,25)(110,40) \drawline(135,25)(110,60)
\drawline(135,0)(110,25)
\drawline(110,40)(100,42) \drawline(110,40)(100,46) \drawline(110,40)(100,50)
\drawline(110,50)(100,53) \drawline(110,50)(100,53) \drawline(110,50)(100,60)
\drawline(110,60)(100,64) \drawline(110,60)(100,68) \drawline(110,60)(100,72)
\drawline(110,25)(100,30) \drawline(110,25)(100,35) \drawline(110,25)(100,40)
\drawline(110,0)(100,10)

\drawline(235,0)(260,-25) \drawline(260,-25)(270,-30) \drawline(260,-25)(270,-35) \drawline(260,-25)(270,-40)

\drawline(210,0)(260,-50)
\drawline(235,-25)(260,-40) \drawline(235,-25)(260,-60) 
\drawline(260,-40)(270,-42) \drawline(260,-40)(270,-46) \drawline(260,-40)(270,-50)
\drawline(260,-50)(270,-53) \drawline(260,-50)(270,-56) \drawline(260,-50)(270,-60)
\drawline(260,-60)(270,-64) \drawline(260,-60)(270,-68) \drawline(260,-60)(270,-72)

\drawline(185,0)(210,-25) \drawline(210,-25)(235,-50)  \drawline(210,-25)(235,-60) \drawline(210,-25)(235,-70)
\drawline(235,-50)(260,-70) \drawline(235,-50)(260,-75) \drawline(235,-50)(260,-80)
\drawline(235,-60)(260,-85) \drawline(235,-60)(260,-90) \drawline(235,-60)(260,-95)
\drawline(235,-70)(260,-100) \drawline(235,-70)(260,-105) \drawline(235,-70)(260,-110)
\jput(260,-90){\makebox(0,0){$\bullet$}}
\drawline(260,-80)(270,-81) \drawline(260,-80)(270,-82) \drawline(260,-80)(270,-83)
\drawline(260,-85)(270,-86) \drawline(260,-85)(270,-87) \drawline(260,-85)(270,-88)
\drawline(260,-90)(270,-91) \drawline(260,-90)(270,-92) \drawline(260,-90)(270,-93)
\drawline(260,-70)(270,-71) \drawline(260,-70)(270,-72) \drawline(260,-70)(270,-73)
\drawline(260,-75)(270,-76) \drawline(260,-75)(270,-77) \drawline(260,-75)(270,-78)
\drawline(260,-95)(270,-96) \drawline(260,-95)(270,-97) \drawline(260,-95)(270,-98)
\drawline(260,-100)(270,-101) \drawline(260,-100)(270,-102) \drawline(260,-100)(270,-103)
\drawline(260,-105)(270,-106) \drawline(260,-105)(270,-107) \drawline(260,-105)(270,-108)
\drawline(260,-110)(270,-111) \drawline(260,-110)(270,-112) \drawline(260,-110)(270,-113)

\drawline(160,0)(110,-50)
\drawline(135,-25)(110,-40) \drawline(135,-25)(110,-60)
\drawline(135,0)(110,-25)
\drawline(110,-40)(100,-42) \drawline(110,-40)(100,-46) \drawline(110,-40)(100,-50)
\drawline(110,-50)(100,-53) \drawline(110,-50)(100,-53) \drawline(110,-50)(100,-60)
\drawline(110,-60)(100,-64) \drawline(110,-60)(100,-68) \drawline(110,-60)(100,-72)
\drawline(110,-25)(100,-30) \drawline(110,-25)(100,-35) \drawline(110,-25)(100,-40)
\drawline(110,0)(100,-10)
\jput(170,2){\makebox(0,0){$|$}}
\jput(160,-12){{\tiny $[\tfrac{s}{p+1},y]_0$}}

\drawline(110,80)(110,-80)
\drawline(260,120)(260,-120)
\drawline(210,50)(210,-50)
\drawline(160,10)(160,-10)
\jput(108,-130){$\mathbb{A}^2$}
\jput(158,-130){$\mathbb{A}^0$}
\jput(208,-130){$\mathbb{A}^1$}
\jput(258,-130){$\mathbb{A}^3$}

\end{picture}
\vspace{50mm}
\caption{The stratification with affine spaces in the building. Fat points mark the image of 
an exemplary morphism $\bar\chi$.}
\end{figure}\\
\noindent {\emph{(b)}} For a given collection $\bfv$ we have
\begin{equation}\label{GRdecomp}
\GRloc=\bigcup_{\substack{a+b=2e-d'\\e-r_1\leq b\leq a\leq e-r_2}}\G(a,b),
\end{equation}
where $d'$ is the integer defined in $(\ref{dstrich})$.
Hence the scheme is geometrically irreducible, because the restriction of the order "$\leq$" on the pairs
\[\{(a,b)\in\mathbb{Z}^2\mid a+b=m(\bfv)\},\] where $m(\bfv)$ is given by $(\ref{mvirred})$, is a total order. 
Of course this also implies connectedness.\\
The dimension of $\GRloc$ is given by the dimension of the maximal affine space in $(\ref{GRdecomp})$.
We assume that $\epsilon$ is even, i.e. $\lfloor\tfrac{r_1-r_2}{p+1}\rfloor+x_0\equiv m(\bfv)\mod 2$. 
The computations in the other case are similar. \\
In this case the affine subspace of maximal dimension consists of all lattices \\$\latM\in\bar\B(m(\bfv))$ with 
\[d_1(\latM,P_{\rm{irred}})=d_1([x_0-\lfloor\tfrac{r_1-r_2}{p+1}\rfloor,m(\bfv)]_0,P_{\rm{irred}}).\]
(if the latter distance is $\leq \tfrac{r_1-r_2}{p+1}$) 
or of the lattices with 
\[d_1(\latM,P_{\rm{irred}})=d_1([x_0+\lfloor\tfrac{r_1-r_2}{p+1}\rfloor,m(\bfv)]_0,P_{\rm{irred}})\]
(if $d_1([x_0-\lfloor\tfrac{r_1-r_2}{p+1}\rfloor,m(\bfv)]_0,P_{\rm{irred}})>\tfrac{r_1-r_2}{p+1}$).
Hence its dimension is either $n:=\lfloor\tfrac{r_1-r_2}{p+1}\rfloor$ (in the first case) or $n-1$ (in the second case).
This yields the claim on the dimension:
\begin{align*}
\dim\,\GRloc &=\lfloor\tfrac{r_1-r_2}{p+1}-\tfrac{s}{p+1}\rfloor+\lfloor\tfrac{s}{p+1}\rfloor\\
&=\begin{cases}
n &\text{if}\ \tfrac{s}{p+1}-\lfloor\tfrac{s}{p+1}\rfloor\leq \tfrac{r_1-r_2}{p+1}-\lfloor\tfrac{r_1-r_2}{p+1}\rfloor\\
n-1 &\text{if}\ \tfrac{s}{p+1}-\lfloor\tfrac{s}{p+1}\rfloor> \tfrac{r_1-r_2}{p+1}-\lfloor\tfrac{r_1-r_2}{p+1}\rfloor.
\end{cases}
\end{align*}
We further see that the set of $\bfv$-admissible lattices is 
exactly the set of lattices in $\bar\B(m(\bfv))$ with 
\begin{equation}\label{fall1}
\begin{aligned}
&\ d_1(\latM,[x_0,m(\bfv)]_0)\leq n && \text{if}\ x_0+\lfloor\tfrac{r_1-r_2}{p+1}\rfloor-\tfrac{s}{p+1}\leq\tfrac{r_1-r_2}{p+1}\\
&\ d_1(\latM,[x_0+1,m(\bfv)]_0)\leq n-1 && \text{otherwise}
\end{aligned}
\end{equation}
and hence this is the set of lattices whose elementary divisors $(a,b)$ with respect to a lattice $\mathfrak{N}$ satisfy 
$(a,b)\leq (a_{\max},b_{\max})$ for some given integers $a_{\max},b_{\max}$.
For $\mathfrak{N}$ we choose one of the lattices 
\[[x_0,m(\bfv)]_0\ ,\ [x_0,m(\bfv)-1]_0\hspace{3mm} \text{or}\hspace{3mm}[x_0+1,m(\bfv)]_0\ ,\ [x_0+1,m(\bfv)-1]_0\]
depending on the cases as listed in $(\ref{fall1})$ and on $x_0-m(\bfv)\mod 2$.
Since we know that $\GRloc$ is reduced, we find that it is isomorphic to a Schubert variety in the affine Grassmannian
after extending the scalars to $\F'$. \\
If $\epsilon$ is odd, then the maximal affine subspace consists of all lattices $\latM\in\bar\B(m(\bfv))$ with
\[d_1(\latM,P_{\rm{irred}})=d_1([x_0+1+\lfloor\tfrac{r_1-r_2}{p+1}\rfloor,m(\bfv)]_0,P_{\rm{irred}})\]
(if the latter distance is $\leq \tfrac{r_1-r_2}{p+1}$) or of the lattices with
\[d_1(\latM,P_{\rm{irred}})=d_1([x_0+1-\lfloor\tfrac{r_1-r_2}{p+1}\rfloor,m(\bfv)]_0,P_{\rm{irred}})\]
(if $d_1([x_0+1+\lfloor\tfrac{r_1-r_2}{p+1}\rfloor,m(\bfv)]_0,P_{\rm{irred}})>\tfrac{r_1-r_2}{p+1}$).
We find
\begin{align*}
\dim\,\GRloc &= \lfloor\tfrac{r_1-r_2}{p+1}+\tfrac{s}{p+1}\rfloor+\lfloor -\tfrac{s}{p+1}\rfloor\\
&=\begin{cases}
n &\text{if}\ 1-(\tfrac{s}{p+1}-\lfloor\tfrac{s}{p+1}\rfloor)\leq \tfrac{r_1-r_2}{p+1}-\lfloor\tfrac{r_1-r_2}{p+1}\rfloor  \\
n-1 &\text{if}\ 1-(\tfrac{s}{p+1}-\lfloor\tfrac{s}{p+1}\rfloor)> \tfrac{r_1-r_2}{p+1}-\lfloor\tfrac{r_1-r_2}{p+1}\rfloor
\end{cases}
\end{align*}
and the conclusion for the isomorphism with a Schubert variety is similar.
\end{proof}
As a consequence of the theorem, we may determine the cases when $\GRloc$ is a single point.
\begin{cor}\label{singptirred}
Denote by $\xi=\tfrac{s}{p+1}-\lfloor\tfrac{s}{p+1}\rfloor$ the fractional part of $\tfrac{s}{p+1}$.
\[\GRloc =\{\ast\}\Leftrightarrow 
\begin{cases}
\ 0+\xi\leq \tfrac{r_1-r_2}{p+1}<2-\xi & \text{if}\ \lfloor\tfrac{s}{p+1}\rfloor\equiv \tfrac{2e-d'-s}{p-1}\mod 2\\
\ 1-\xi\leq \tfrac{r_1-r_2}{p+1}<1+\xi & \text{if}\ \lfloor\tfrac{s}{p+1}\rfloor\not\equiv \tfrac{2e-d'-s}{p-1}\mod 2.
\end{cases}\]
\end{cor}
\begin{proof}
This is just the case where the dimension of $\GRloc$ is zero. More explicitly:\\
If $x_0=\lfloor\tfrac{s}{p+1}\rfloor\equiv m(\bfv)\mod 2$, then
$[x_0,m(\bfv)]$ is the unique lattice with minimal distance $d_1$ from $P_{\rm{irred}}$. We have 
\[d_1([x_0,m(\bfv)]_0,P_{\rm{irred}})=\xi.\]
Thus this lattice is $\bfv$-admissible if and only if $\tfrac{r_1-r_2}{p+1}\geq \xi$.
There is no other $\bfv$-admissible lattice iff the lattices $\latM$ with 
$d_1(\latM,P_{\rm{irred}})=2-\xi$ are not $\bfv$-admissible. This yields the claim.\\
The case $x_0\not\equiv m(\bfv)\mod 2$ is similar. 
Instead of $[x_0,m(\bfv)]_0$ we have to consider the lattice $[x_0+1,m(\bfv)]_0$.
\end{proof}
\numberwithin{equation}{subsection}
\section{The reducible case}
In this section we want to analyze the case, where $(M_\F,\Phi)$ admits a proper $\Phi$-stable subobject,
at least after extending the scalars to some finite extension of $\F$. Before we start to determine the set of
$\bfv$-admissible lattices in the building, we want to formulate the precise statement on 
the connected components of $\GRloc$. We first define some open and closed subschemes of $\GRloc$.
\begin{defn}\label{1dimM}
For $a\in\bar\F^\times$ and $j\in\mathbb{Z}_{\geq 0}$ define $(\latM^j(a),\Phi_a^j)$ by
\[\latM^j(a)=\Fbarpot\ ,\hspace{5mm} \Phi_a^j(1)=au^j.\]
\end{defn}
\begin{defn}\label{defvord}
A $\bfv$-admissible lattice $\latM\subset M_{\bar\F}$ is called $\bfv$-{\emph{ordinary}} if there exists a short exact sequence
\begin{equation}\label{extension}
0\rightarrow (\latM^{e-r_1}(a),\Phi_a^{e-r_1})\rightarrow (\latM,\Phi)
\rightarrow (\latM^{e-r_2}(b),\Phi_b^{e-r_2})\rightarrow 0
\end{equation}
for some $a,b\in\bar\F^\times$.
\end{defn}
\begin{rem}\label{remvord}
The determinant condition in $(\ref{GRlocpts})$ implies that
\begin{equation}\label{det}
u^{e-r_1}\latM\subset \langle\Phi(\latM)\rangle\subset u^{e-r_2}\latM
\end{equation}
for all $\bfv$-admissible lattices $\latM$. Hence the $\bfv$-ordinary lattices are the lattices which admit a $\Phi$-stable
subobject with the minimal possible elementary divisors.
If a $\bfv$-admissible lattice $(\latM,\Phi)$ admits a subobject isomorphic to $(\latM^{e-r_1}(a),\Phi_a^{e-r_1})$
for some $a\in\bar\F^\times$, then the quotient has no $u$-torsion by $(\ref{det})$ and is 
isomorphic to $(\latM^{e-r_2}(b),\Phi_b^{e-r_2})$ for some $b\in\bar\F^\times$,
because the sum of the elementary divisors is fixed by $(\ref{GRlocpts})$. Hence $(\latM,\Phi)$ is $\bfv$-ordinary in 
this case.
\end{rem}
Denote by $\mathcal{S}(\bfv)$ the set of isomorphism classes of one dimensional $\bar\F((u))$-modules $M'$ with 
$\phi$-linear map $\Phi'\neq 0$ such that $M'$ admits a (unique) lattice $\latM_{[M']}\subset M'$ with 
$\langle\Phi(\latM_{[M']})\rangle=u^{e-r_1}\latM_{[M']}$. The elements of $\mathcal{S}(\bfv)$ are in bijection with the 
elements of $\bar\F^\times$: For each $a\in\bar\F^\times$ there is a unique isomorphism class represented by 
\begin{equation}\label{MaPhia}
(M_a,\Phi_a)=(\latM^{e-r_1}(a)[\tfrac{1}{u}],\Phi_a^{e-r_1}).
\end{equation}
Set $X=\GRloc\otimes_\F\bar\F$. On $X$ there is a universal sheaf of 
$\Fbarpot\widehat{\otimes}_{\bar\F}\shOX=\shOX[\hspace{-0.5mm}[u]\hspace{-0.5mm}]$-lattices 
$\mathcal{M}\subset M_{\bar\F}\widehat{\otimes}_{\bar\F}\shOX$ satisfying
\[u^e\mathcal{M}\subset (\rm{id}\otimes\Phi)\phi^*\mathcal{M}\subset\mathcal{M}.\]
For each $[M']\in\mathcal{S}(\bfv)$ define a sheaf of $\shOX$-modules 
\begin{equation}\label{defFM}
\mathcal{F}_{[M']}=\mathcal{H}om_{\shOX[\hspace{-0.5mm}[u]\hspace{-0.5mm}],\Phi}
(\latM_{[M']}\widehat{\otimes}_{\bar\F}\shOX,\mathcal{M})
\end{equation}
where the subscript $\Phi$ indicates that the homomorphism have to commute with the semi-linear maps that are part of the data.
\begin{prop}\label{propcoh}
\noindent (i) For each $[M']\in\mathcal{S}(\bfv)$ the sheaf $\mathcal{F}_{[M']}$ is a coherent $\shOX$-module.\\
\noindent (ii) A closed point $x\in\GRloc$ corresponds to a non-$\bfv$-ordinary lattice if and only if 
$\mathcal{F}_{[M']}\otimes \kappa(x)=0$ for all $[M']\in\mathcal{S}(\bfv)$.
\end{prop}
\begin{proof}
\noindent \emph{(i)} For the isomorphism class $[M']$ we choose a representative of the form $M_a$ defined in $(\ref{MaPhia})$.
Let $U= {\rm{Spec}}\,A\subset X$ an affine open. We claim\\
\emph{\noindent (a) ${\rm{Hom}}_{A[\hspace{-0.5mm}[u]\hspace{-0.5mm}],\Phi}(\latM_{[M']}\widehat{\otimes}_{\bar\F}A,\mathcal{M}(U))$
is a finitely generated $A$-module.\\
\noindent (b) If $V={\rm{Spec}}\,B\subset U$ is an affine open we have}
\begin{equation}\label{localization}
{\rm{Hom}}_{B[\hspace{-0.5mm}[u]\hspace{-0.5mm}],\Phi}(\latM_{[M']}\widehat{\otimes}_{\bar\F}B,\mathcal{M}(V))
\cong {\rm{Hom}}_{A[\hspace{-0.5mm}[u]\hspace{-0.5mm}],\Phi}(\latM_{[M']}\widehat{\otimes}_{\bar\F}A,\mathcal{M}(U))\otimes_AB.
\end{equation}
This implies the first part of the Proposition.\\
\emph{Proof of (a):}
Because $\latM_{[M']}\widehat{\otimes}_{\bar\F}A$ is a free $A[\hspace{-0.5mm}[u]\hspace{-0.5mm}]$-module of rank one,
a morphism is given by the image of $1$ and hence
\[{\rm{Hom}}_{A[\hspace{-0.5mm}[u]\hspace{-0.5mm}],\Phi}(\latM_{[M']}\widehat{\otimes}_{\bar\F}A,\mathcal{M}(U))
\cong N_A\subset \mathcal{M}(U),\]
where $N_A$ is the $A$-submodule of all $v\in\mathcal{M}(U)$ satisfying $\Phi(v)=au^{e-r_1}v$.
We claim that the reduction modulo $u^{e+1}$ induces an injective homomorphism
\[N_A\hookrightarrow\mathcal{M}(U)/u^{e+1}\mathcal{M}(U)\]
and hence $N_A$ is finitely generated as an $A$-module, because the scheme $X$ is noetherian.
Now, if $0\neq v=u^nw\in N_A$ with $n\geq 0$ and $w\in\mathcal{M}(U)\backslash u\mathcal{M}(U)$, then
\[u^{pn}\Phi(w)=\Phi(u^nw)=au^{e-r_1+n}w\]
and hence $0\leq e-r_1-(p-1)n\leq e$ which implies $n\leq e$.\\
{\emph{Proof of (b):}} We have the following commutative diagram
\[\begin{xy}
\xymatrix{
{\rm{Hom}}_{A[\hspace{-0.5mm}[u]\hspace{-0.5mm}],\Phi}(\latM_{[M']}\widehat{\otimes}_{\bar\F}A,\mathcal{M}(U)) 
\ar[d]\ar[r]^/1cm/{\cong}
& N_A\ar[d] \ar@^{(->}[r] & \mathcal{M}(U)\ar[d]\\
{\rm{Hom}}_{B[\hspace{-0.5mm}[u]\hspace{-0.5mm}],\Phi}(\latM_{[M']}\widehat{\otimes}_{\bar\F}B,\mathcal{M}(V))
\ar[r]^/1cm/{\cong}
& N_B \ar@^{(->}[r] & **[r]\mathcal{M}(V)\cong\mathcal{M}(U)\widehat{\otimes}_AB}
\end{xy}
\]
As $N_A$ is a finitely generated $A$-module, we do not need to complete the tensor product
to obtain $N_B$ from $N_A$ (there are only finitely many denominators). Hence $(\ref{localization})$ is an isomorphism.\\
\noindent \emph{(ii)} Let $[M']\in\mathcal{S}(\bfv)$ be an isomorphism class and suppose that $x\in\GRloc$ is a closed
point corresponding to a lattice $\latM$ such that $\mathcal{F}_{[M']}\otimes\kappa(x)\neq 0$, i.e. there exists 
a non trivial morphism
\[f:\latM_{[M']}\rightarrow \latM.\]
As both sides are free $\Fbarpot$-modules and the morphism is non trivial, it is injective.
We have to convince ourselves that ${\rm{coker}}\,f$ has no $u$-torsion: in this case $\latM$ is the extension of 
free $\Fbarpot$-modules of rank $1$ (an extension of ${\rm{coker}}\,f$ by ${\rm{im}}\,f$), and hence $\bfv$-ordinary.\\
We write $f(1)=u^nv$ for some $n\in\mathbb{Z}$ and $v\in\latM\backslash u\latM$ and claim $n=0$.
Because of $\Phi(f(1))=f(\Phi(1))$ we find $\Phi(v)=au^{e-r_1-(p-1)n}\in\Phi(\latM)\subset u^{e-r_1}\latM$
for some $a\in\bar\F^\times$ and hence $n=0$.\\
Conversely, if $\latM$ is $\bfv$-ordinary, then the inclusion of the $\Phi$-stable subobject defines a nontrivial
morphism $\latM_{[M']}\rightarrow \latM$ for some $[M']\in\mathcal{S}(\bfv)$.  
\end{proof}
\begin{defn}\label{defX0}
For each isomorphism class $[M']\in\mathcal{S}(\bfv)$ define 
\[X_{[M']}^{\bfv}=\{x\in\GRloc\otimes_\F\bar\F\mid \mathcal{F}_{[M']}\otimes\kappa(x)\neq 0\}.\]
Further define 
\[X_0^{\bfv}=\GRloc\otimes_\F\bar\F\ \backslash\bigcup_{[M']\in\mathcal{S}(\bfv)}X_{[M']}^{\bfv}.\]
By the Proposition below these subsets are open and closed and hence they come along with a canonical scheme structure.
\end{defn}
\begin{prop}\label{openclosed}
\noindent (i) The subset $X_{[M']}^\bfv$ is open and closed for each $[M']\in\mathcal{S}(\bfv)$.\\
\noindent (ii) The subset $X_0^\bfv$ is open and closed.
\end{prop}
\begin{proof}
\noindent \emph{(i)} It is clear that $X^\bfv_{[M']}$ is closed, as $\mathcal{F}_{[M']}$ is coherent.
We show that it is closed under cospecialization.\\
Let $\eta\rightsquigarrow x$ be a specialization with $x\in X^\bfv_{[M']}$ and assume that $x$ is a closed point.
We mark this specialization by ${\rm{Spec}}\,R\rightarrow X$, where $R$ is a discrete valuation ring with uniformizer $t$
and residue field $\bar\F$. Denote by $\latM_R$ the $R[\hspace{-0.5mm}[u]\hspace{-0.5mm}]$-lattice in 
$M_{\bar\F}\widehat{\otimes}_{\bar\F}R$ defined by this morphism.
Because of $\mathcal{F}_{[M']}\otimes\kappa(x)\neq 0$, there is a non trivial morphism $\latM_{[M']}\rightarrow \latM_x$
and hence there is a basis vector $b_1\in M_{\bar\F}$ such that $\Phi(b_1)=au^{e-r_1}b_1$ for some $a\in\bar\F^\times$.
As $\latM_R$ is a free $R[\hspace{-0.5mm}[u]\hspace{-0.5mm}]$-module, there is a basis of $\latM_R$ such that
\[\latM_R\sim\begin{pmatrix}\alpha & \gamma\\ 0 & \delta \end{pmatrix},\]
for $\alpha,\gamma,\delta\in R[\hspace{-0.5mm}[u]\hspace{-0.5mm}]$,
with $\alpha\equiv au^{e-r_1}\mod t$. But the determinant condition in $(\ref{GRlocpts})$ implies $v_u(\alpha)\geq e-r_1$.
Hence $v_u(\alpha)=e-r_1$ and $\eta\in X^\bfv_{[M'']}$ for some $[M'']\in\mathcal{S}(\bfv)$.
If $[M']=[M'']$ we are done.\\
Assume $[M']\neq [M'']$. As $X^\bfv_{[M'']}$ is closed, we have $x\in X^\bfv_{[M']}\cap X^\bfv_{[M'']}$.
In this case $\latM_x$ admits two linear independent subspaces:
\[\latM_x\sim\begin{pmatrix} au^{e-r_1} & 0\\ 0 & bu^{e-r_1}\end{pmatrix}\]
and hence $e-r_1=e-r_2$. Now we easily deduce $\GRloc=\{\latM_x\}$ and the claim follows.\\
\noindent \emph{(ii)} This follows from the first part of the Proposition together with the fact that the one-dimensional
$\Phi$-invariant subspaces of $M_{\bar\F}$ which admit an integral model $\latM$ with 
$\langle\Phi(\latM)\rangle=u^{e-r_1}\latM$ run over a finite set of isomorphism classes of one-dimensional objects:\\
Assume that there are two different one-dimensional $\Phi$-stable subspaces $\langle b_1\rangle$ and $\langle b_2\rangle$
of $M_{\bar\F}$ such that $\Phi(b_i)=a_iu^{e-r_1}b_i$, for $i=1,2$. Then $b_1$ and $b_2$ are linear independent.\\
If $a_1\neq a_2$, then $\langle b_1+qb_2\rangle$ is not $\Phi$-stable for all $q\in\bar\F((u))^\times$
and hence there are only two isomorphism classes.\\
If $a_1=a_2$, then there is a unique such isomorphism class given by $[M_a]$.
\end{proof}
We will see below that the open and closed subschemes $X_{[M']}^\bfv$ and $X_0^\bfv$ of $\GRloc\otimes_\F\bar\F$
are connected and hence turn out to be the connected components of $\GRloc\otimes_{\F}\bar\F$.

Now we want to determine the subset of $\bfv$-admissible lattices in the building. As we are assuming that $(M_\F,\Phi)$
is reducible, at least after extending scalars, there exists a finite extension $\F'$ of $\F$ and a basis $e_1,e_2$
of $M_{\F'}=M_\F\widehat{\otimes}_\F\F'$ such that
\[M_{\F'}\sim \begin{pmatrix}au^s & \gamma\\0&bu^t\end{pmatrix}\]
for some $a,b\in\F'^\times,\,\gamma\in\F'((u))$ and $s,t\in\mathbb{Z}$ with $0\leq s,t<p-1$.
We choose this basis to be the standard basis.
\begin{lem}\label{lemred}
\noindent (i) The map $\Phi$ extends to a map $\bar\B\rightarrow \bar\B$ also denoted by $\Phi$.\\
\noindent (ii) For $q\in\bar\F((u))$ and $[x,y]_q\in\A_q$ the map $\Phi$ is given by 
\[\Phi([x,y]_q)=[px+s-t,py+s+t]_{q'}\]
with $q'=b^{-1}u^{-t}(au^s\phi(q)+\gamma)$.
\end{lem}
\begin{proof}
\noindent \emph{(i)} We can use the expressions in \emph{(ii)} to extend $\Phi$.\\
\noindent \emph{(ii)} We have $\Phi(u^me_1)=au^{pm+s}e_1$ and
\begin{align*}
\Phi(u^n(qe_1+e_2))&=u^{pn}(au^s\phi(q)e_1+\gamma e_1+bu^te_2)\\
&=bu^{pn+t}(b^{-1}u^{-t}(\gamma +au^s\phi(q))e_1+e_2).
\end{align*}
The Lemma follows from this.
\end{proof}
\begin{cor}\label{emptyred}
The scheme $\GRloc$ is empty if $2e-d'\not\equiv s+t\mod (p-1)$.
\end{cor}
\begin{proof}
This follows from Lemma $\ref{lemdist}$ and Lemma $\ref{lemred}$: We have
\[d_2([x,y]_q,\Phi([x,y]_q))=(p-1)y+s+t\]
and this distance must be equal to $2e-d'$ if $[x,y]_q$ is $\bfv$-admissible.
\end{proof}
We assume that the scheme is non empty and define 
\begin{align}
&\ P_{\rm{red}}=[\tfrac{t-s}{p-1},-\tfrac{t+s}{p-1}]\in\A_0\subset\bar\B \label{ptPred}\\
&\ m(\bfv)=\tfrac{2e-d'-(s+t)}{p-1}\in\mathbb{Z}. \label{mvred}
\end{align}
These definitions imply $\GRloc(\bar\F)\subset\bar\B(m(\bfv))$.\\
There are three different cases which we have to study in order to determine the set of $\bfv$-admissible lattices.
It makes a difference whether $(M_{\bar\F},\Phi)$ is a split or a non-split extension of two one-dimensional objects. 
In the split case there are two possibilities: Either the direct summands are isomorphic or non-isomorphic.
\subsection{The case $(M_{\bar\F},\Phi)\cong (M_1,\Phi_1)\oplus (M_1,\Phi_1)$.}
In this section we want to analyze the case where $(M_\F,\Phi)$ becomes isomorphic to a direct sum of two isomorphic
one-dimensional objects after possibly extending the scalars to some finite extension of $\F$, i.e. we want to assume
that there exists $\F'/\F$ and an $\F'((u))$-basis $e_1,e_2$ of $M_{\F'}$ such that
\begin{equation}\label{phisplit1}
M_{\F'}\sim\begin{pmatrix} au^s&0\\0&au^s\end{pmatrix}
\end{equation}
with $a\in\F'^\times$ and $0\leq s<p-1$. We immediately find $\Phi(P_{\rm{red}})=P_{\rm{red}}$.\\
For each $z\in\mathbb{P}^1(\bar\F)$ we define a (half)-line $\mathcal{L}_z\subset \bar\B(m(\bfv))$ by
\begin{align*}
&\ \calL_z\ =\{[x,m(\bfv)]_z\mid x\geq 0\}\subset \bar\B(m(\bfv)) && \text{if}\ z\in\bar\F=\mathbb{A}^1(\bar\F)\\
&\ \calL_{\infty}=\{[x,m(\bfv)]_0\mid x\leq 0\}\subset \bar\B(m(\bfv)). 
\end{align*} 
These lines are defined in such a way that 
\begin{equation}\label{excepttree}
\T:=\bigcup_{z\in\mathbb{P}^1(\bar\F)}\calL_z=\bigcup_{z\in\bar\F}\A_z\cap\bar\B(m(\bfv)).
\end{equation}
The apartments on the right hand side are given by the basis $e_1, ze_1+e_2$ and in this basis the semi-linear endomorphism
$\Phi$ is of the form $(\ref{phisplit1})$.
\begin{lem}\label{distred1}
Let $Q\in\bar\B(m(\bfv))$ be an arbitrary point. Let $Q'\in\T$ be the unique point satisfying
$d_1(Q,Q')=d_1(Q,\T)$. Then
\begin{align*}
&\ d_1(Q,\Phi(Q))= (p+1)d_1(Q,P_{\rm{red}})-2d_1(Q',P_{\rm{red}})\\
&\ d_2(Q,\Phi(Q))= (p-1)d_2(Q,P_{\rm{red}}).
\end{align*}
\end{lem}
\begin{proof}
The statement on $d_2$ follows immediately from Lemma $\ref{lemred}$. For the statement on $d_1$ we assume $Q'\in\calL_0$.
The cases $Q'\in\calL_z$ for $z\in\bar\F$ are analogous and the case $Q'\in\calL_\infty$ is obtained by interchanging
$e_1$ and $e_2$.\\
First assume $Q=Q'$, i.e. $Q=[x,m(\bfv)]_0\in\calL_0$. Then Lemma $\ref{lemred}$ implies \\ $\Phi(Q)=[px,pm(\bfv)+2s]_0$
and hence $d_1(Q,\Phi(Q))=(p-1)x=(p-1)d_1(Q,P_{\rm{red}})$.\\
Now assume $Q\neq Q'$. We write
\[Q=[x,m(\bfv)]_q\ ,\hspace{5mm} Q'=[x',m(\bfv)]_0\]
with $x>x'=v_u(q)\in\mathbb{Z}_{>0}$. Then $\Phi(Q)=[px,pm(\bfv)+2s]_{\phi(q)}$ by Lemma $\ref{lemred}$. Using 
$v_u(\phi(q))=px'$, we find
\begin{align*}
d_1(Q,\Phi(Q))&=(x-x')+(px'-x')+(px-px')\\ &=(p+1)x-2x'=(p+1)d_1(Q,P_{\rm{red}})-2d_1(Q',P_{\rm{red}}).
\end{align*}
\end{proof}
\begin{rem}
This Lemma shows that the case of a direct sum of two isomorphic objects corresponds to the 
case B 2 in \cite{phimod} 6.d:\\
The unique point fixed by $\Phi$ is the point $P_{\rm{red}}$ and the projection of this point to the building for $PGL_2(\bar\F((u)))$ is a vertex. The link of this vertex is the projection of $\T$
and all the half-lines $\calL_z$ of $\T$ (for $z\in\mathbb{P}^1(\bar\F)$) are fixed by $\Phi$.
\end{rem}
\begin{prop}\label{propred1}
With the notations of Definition $\ref{1dimM}$ and $(\ref{MaPhia})$, $(\ref{mvred})$ assume that 
\[(M_{\bar\F},\Phi)\cong (\latM^s(a)[\tfrac{1}{u}],\Phi^s_a)\oplus(\latM^s(a)[\tfrac{1}{u}],\Phi^s_a)\]
for some $a\in\bar\F^\times$ and $0\leq s<p-1$.\\
\noindent (i) The schemes $X^\bfv_{[M']}$ are empty for all $[M']\in\mathcal{S}(\bfv)\backslash\{[M_a]\}$.\\
\noindent (ii) The scheme $X^\bfv_{[M_a]}$ is given by
\[X^\bfv_{[M_a]}\cong\begin{cases}
\ \emptyset & \text{if}\ m(\bfv)+\tfrac{r_1-r_2}{p-1}\notin 2\mathbb{Z}\\
\{\ast\} & \text{if}\ 0=\tfrac{r_1-r_2}{p-1}\in\mathbb{Z}\ \text{and} \ \tfrac{r_1-r_2}{p-1}\equiv m(\bfv) \mod 2\\
\ \mathbb{P}^1_{\bar\F} & \text{if}\ 0\neq \tfrac{r_1-r_2}{p-1}\in\mathbb{Z}\ \text{and}\ \tfrac{r_1-r_2}{p-1}\equiv m(\bfv)\mod 2.
\end{cases}\]
\noindent (iii) If non empty, the scheme $X^\bfv_0$ is connected. 
\end{prop}
\begin{proof}
We first claim that every $\bfv$-admissible lattice $\latM$ can be linked to a $\bfv$-admissible lattice $\latM'\in\T$
by a chain of $\mathbb{P}^1$.\\
Assume $\latM=[x,m(\bfv)]_q\notin\T$ and let $Q'\in\T$ be the unique point satisfying $d_1(\latM,Q')=d_1(\latM,\T)$.
Without loss of generality, we may again assume that $Q'=[x',m(\bfv)]_0\in\calL_0$.
By construction we have $\latM,Q'\in\A_q$ and we choose the following basis $b_1,b_2$ of $\latM$: 
\[b_1=u^{(x+m(\bfv))/2}e_1\ ,\hspace{5mm} b_2=u^{(m(\bfv)-x)/2}(qe_1+e_2).\]
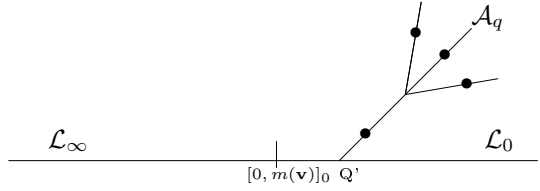
\begin{figure}[h]
\begin{picture}(200,50)
\drawline(0,0)(200,0)
\drawline(125,0)(175,50)
\drawline(150,25)(156,60)(150,25)(185,31)
\jput(135,10){\makebox(0,0){$\bullet$}}
\jput(165,40){\makebox(0,0){$\bullet$}}
\jput(154,48.3){\makebox(0,0){$\bullet$}} \jput(173.3,29){\makebox(0,0){$\bullet$}}
\jput(100,0){|}
\jput(90,-7){{\tiny $[0,m(\bfv)]_0$}}
\jput(125,-7){{\tiny Q'}}
\jput(180,5){$\calL_0$}\jput(15,5){$\calL_\infty$}\jput(175,52){$\A_q$}
\end{picture}
\caption{The fat points mark the image of the morphism $\bar\chi$ in the building in the case $p=3$ and $\F=\F_3$.}
\end{figure}\\
Applying Lemma $\ref{P1nachgrass}$ with this basis yields a morphism 
$\bar\chi:\mathbb{P}^1_{\bar\F}\rightarrow {\rm{Grass}}\, M_{\bar\F}$ with $\bar\chi(z)=[x,m(\bfv)]_{q+zu^{x-1}}$ for 
$z\in\bar\F=\mathbb{A}^1(\bar\F)$ and $\bar\chi(\infty)=[x-2,m(\bfv)]_q$. We have 
\[d_1(\bar\chi(\infty),P_{\rm{red}})<d_1(\bar\chi(z),P_{\rm{red}})=d_1(\latM,P_{\rm{red}})\]
for all $z\in\bar\F$, while $d_1(\bar\chi(z),\T)\leq d_1(\latM,\T)$ for all $z\in\mathbb{P}^1(\bar\F)$
and, by construction, $d_1(\bar\chi(\infty),\T)<d_1(\latM,\T)$.\\
By Lemma $\ref{distred1}$ and Lemma $\ref{lemdist}$, the morphism $\bar\chi$ factors through 
$\GRloc$ and the claim follows by induction on the distance $d_1(\latM,\T)$.\\
Now we assume that $\latM\in\T$ is a $\bfv$-admissible lattice and we are looking for a $\bfv$-admissible lattice
$\latM'$ that can be linked with $\latM$ by a $\mathbb{P}^1$ and that has strictly smaller distance $d_1$
from $P_{\rm{red}}=[0,\tfrac{-2s}{p-1}]_0$ than $\latM$, i.e. $d_1(\latM',P_{\rm{red}})<d_1(\latM,P_{\rm{red}})$.\\
We may assume $\latM=[x,m(\bfv)]_0\in\calL_0$. Assuming $x>1$, our candidate for $\latM'$ is $[x-2,m(\bfv)]_0$.
Fixing a basis 
\[b_1=u^{(x+m(\bfv))/2}e_1\ ,\ b_2=u^{(m(\bfv)-x)/2}e_2\]
of $\latM$ so that $\latM'=\langle u^{-1}b_1,ub_2 \rangle$, yields a morphism
$\bar\chi:\mathbb{P}^1_{\bar\F}\rightarrow {\rm{Grass}}\,M_{\bar\F}$ with $\bar\chi(0)=\latM$ and $\bar\chi(\infty)=\latM'$.
This morphism factors through $\GRloc$ iff the lattices $\bar\chi(z)=[x,m(\bfv)]_{zu^{x-1}}$ are $\bfv$-admissible
for all $z\in\bar\F\backslash\{0\}$.This is the case iff
\begin{align*}
d_1(\bar\chi(z),\Phi(\bar\chi(z)))& =(p+1)d_1(\bar\chi(z),P_{\rm{red}})-2d_1([x-1,m(\bfv)]_0,P_{\rm{red}})\\
& =(p+1)x-2(x-1)=(p-1)x+2\leq r_1-r_2.
\end{align*}
Consider the following subset of $\bfv$-admissible lattices
\[\mathcal{N}=\{\latM\in\GRloc(\bar\F)\mid 
\latM\notin\T\ \text{or}\ (\latM\in\T\ \text{and}\ d_1(\latM,P_{\rm{red}})\leq\tfrac{r_1-r_2-2}{p-1})\}.\]
So far, we have shown that all $\bfv$-admissible lattices $\latM\in\mathcal{N}$ can either be linked to the lattice $[0,m(\bfv)]_0$ or 
to one of the lattices
\begin{equation}\label{P1center}
\{[1,m(\bfv)]_z\mid z\in\bar\F\}\cup\{[-1,m(\bfv)]_0\}=\{\latM\in\bar\B(m(\bfv))\mid d_1(\latM,P_{\rm{red}})=1\}.
\end{equation}
by a chain of $\mathbb{P}^1$. Here, the two different cases depend on $m(\bfv)\mod 2$.
Hence the subset of $\GRloc\otimes_\F\bar\F$ given by the lattices in $\mathcal{N}$ is connected:
using Remark $\ref{P1inbuilding}$ again, the set in $(\ref{P1center})$ forms a $\mathbb{P}^1$. The Proposition now 
follows from the following two facts:\\
\emph{\noindent (a) If $\latM\in\mathcal{N}$, then $\latM$ is not $\bfv$-ordinary,
i.e. $\mathcal{N}\subset X^\bfv_0(\bar\F)$.\\
\noindent (b) If $\latM\notin\mathcal{N}$ is $\bfv$-admissible, then $\latM\in X^\bfv_{[M_a]}(\bar\F)$ and}
\begin{equation}\label{exceptP1}
\latM\in\{[\tfrac{r_1-r_2}{p-1},m(\bfv)]_z\mid z\in\bar\F\}\cup\{[-\tfrac{r_1-r_2}{p-1},m(\bfv)]_0\}.
\end{equation}
By Remark $\ref{P1inbuilding}$, this set forms a $\mathbb{P}^1$ if $r_1\neq r_2$. Otherwise it is a single point.\\
{\emph{Proof of (a):}} If $\latM\in\T$, then $d_1(\latM,\langle\Phi(\latM)\rangle)\leq r_1-r_2-2<r_1-r_2$ 
and hence the elementary divisors of $\langle\Phi(\latM)\rangle$ with respect to $\latM$ are not given by $(e-r_2,e-r_1)$.\\
If $\latM\notin\T$, say $\latM=[x,m(\bfv)]_q$ with $x>v_u(q)>0$ for example, then $\latM=\langle b_1,b_2\rangle$ with
\[b_1=u^{(x+m(\bfv))/2}e_1\ ,\hspace{5mm} b_2=u^{(m(\bfv)-x)/2}(qe_1+e_2)\]
and one finds
\[\latM\sim (a_{ij})_{ij}=\begin{pmatrix} au^{\tfrac{p-1}{2}(x+m(\bfv))+s} & a\phi(q)u^{\tfrac{p-1}{2}m(\bfv)-\tfrac{p+1}{2}x+s}\\
0 & au^{\tfrac{p-1}{2}(m(\bfv)-x)+s}  \end{pmatrix}\]
with $v_u(a_{12})<v_u(a_{11})$, because $v_u(q)<x$, and hence the minimal elementary divisor of $\langle\Phi(\latM)\rangle$ 
with respect to $\latM$ is not given by a $\Phi$-stable subspace.\\
{\emph{Proof of (b):}} Let $\latM\notin\mathcal{N}$ be $\bfv$-admissible. Then $\latM\in\T$ and 
\[\tfrac{r_1-r_2-2}{p-1}<d_1(\latM,P_{\rm{red}})\leq \tfrac{r_1-r_2}{p-1}.\]
We show that $d_1(\latM,P_{\rm{red}})=\tfrac{r_1-r_2}{p-1}$ which implies $(\ref{exceptP1})$.\\
Suppose that $\latM=[x,m(\bfv)]_z$ with $z\in\bar\F$ and 
\[x=\pm\tfrac{r_1-r_2-1}{p-1}\in\mathbb{Z}\ ,\hspace{5mm} m(\bfv)=\tfrac{2e-d'-2s}{p-1}=\tfrac{2e-r_1-r_2-2s}{p-1}.\]
In this case we find 
\[x+m(\bfv)=\tfrac{2e-2s-(r_1+r_2)\pm(r_1-r_2)\mp 1}{p-1}\notin 2\mathbb{Z},\]
contradiction. We are left to show that $\latM\in X^\bfv_{[M_a]}(\bar\F)$, i.e. that there exists a vector $e_\latM\in\latM$
and a $\Phi$-stable subspace $\Fbarpot e_\latM\subset \latM$ with $\Phi(e_\latM)=au^{e-r_1}e_\latM$.
An easy computation shows that we may choose
\begin{align*}
& e_\latM=u^{\tfrac{e-r_1-s}{p-1}}(ze_1+e_2) && \text{if}\ \latM=[\tfrac{r_1-r_2}{p-1},m(\bfv)]_z ,\ z\in\bar\F\\
& e_\latM=u^{\tfrac{e-r_1-s}{p-1}}e_1 && \text{if}\ \latM=[-\tfrac{r_1-r_2}{p-1},m(\bfv)]_0.
\end{align*}
\end{proof}
We conclude the discussion by determining the cases where $\GRloc$ is reduced to a single point.
\begin{cor}\label{oneptred1}
\noindent (i) If $m(\bfv)\equiv 0\mod 2$, then $\GRloc=\{\ast\}$ iff $\tfrac{r_1-r_2}{p-1}<2$.\\
\noindent (ii) If $m(\bfv)\equiv 1\mod 2$, then $\GRloc$ can not be a single point.
\begin{align*}
& \GRloc=\emptyset &&\Leftrightarrow 0\leq \tfrac{r_1-r_2}{p-1}<1\\
& \GRloc\otimes_\F\bar\F\cong\mathbb{P}^1_{\bar\F}&&\Leftrightarrow 1\leq \tfrac{r_1-r_2}{p-1}<3. 
\end{align*}
\end{cor}
\begin{proof}
\noindent \emph{(i)} As $m(\bfv)\equiv 0\mod 2$, the lattice $[0,m(\bfv)]_0$ is always $\bfv$-admissible.
It is the unique point of $\GRloc$ if the lattices $\latM$ with $d_1(\latM,P_{\rm{red}})=2$ are not $\bfv$-admissible.
By Lemma $\ref{distred1}$ this is the case iff $\tfrac{r_1-r_2}{p-1}<2$.\\
\noindent \emph{(ii)} The scheme is empty if the lattices $\latM$ with $d_1(\latM,P_{\rm{red}})=1$ are not $\bfv$-admissible.
By Lemma $\ref{distred1}$, this is the case iff $\tfrac{r_1-r_2}{p-1}<1$.\\
If $\tfrac{r_1-r_2}{p-1}\geq 1$, then the lattices
\[\{[1,0]_z\mid z\in\bar\F\}\cup \{[-1,0]_0\}\]
are $\bfv$-admissible and form a $\mathbb{P}^1_{\bar\F}$.
Again by Lemma $\ref{distred1}$ there are no other $\bfv$-admissible lattices iff $\tfrac{r_1-r_2}{p-1}<3$.
\end{proof}
\subsection{The case $(M_{\bar\F},\Phi)\cong (M_1,\Phi_1)\oplus(M_2,\Phi_2)$}
In this section we treat the case where $(M_\F,\Phi)$ becomes isomorphic to the direct sum of two non-isomorphic one-dimensional
objects after extending the scalars to some finite extension. The situation is the following:
There exists a finite extension $\F'$ of $\F$ and a basis $e_1,e_2$ of $M_{\F'}$ such that
\[M_{\F'}\sim\begin{pmatrix}au^s&0\\0&bu^t\end{pmatrix}\]
with $a,b\in\F'^\times$ and $0\leq s,t<p-1$. As we are assuming that the direct summands are not isomorphic, we further have 
$s\neq t$ or $a\neq b$. Again we find $\Phi(P_{\rm{red}})=P_{\rm{red}}$.
\begin{lem}\label{lemred2}
Let $Q\in\bar\B(m(\bfv))$ be an arbitrary point. Let $Q'\in\A_0\cap\bar\B(m(\bfv))$ be the unique point satisfying
$d_1(Q,Q')=d_1(Q,\A_0)$. Then
\begin{align*}
&\ d_1(Q,\Phi(Q))=(p+1)d_1(Q,P_{\rm{red}})-2d_1(Q',P_{\rm{red}})\\
&\ d_2(Q,\Phi(Q))=(p-1)d_2(Q,P_{\rm{red}}).
\end{align*}
\end{lem}
\begin{proof}
This is similar to Lemma $\ref{distred1}$. Again the statement on $d_2$ is an immediate consequence of 
Lemma $\ref{lemred}$. Let $Q$ be any point. 
We may assume that the unique point $Q'\in\A_0\cap\bar\B(m(\bfv))$ satisfying $d_1(Q,Q')=d_1(Q,\A_0)$
is given by $[x,m(\bfv)]_0$ with $x\geq \tfrac{t-s}{p-1}$. The case $x\leq \tfrac{t-s}{p-1}$
is obtained by interchanging $e_1$ and $e_2$.\\
If $Q=Q'$, then $Q\in\A_0$ and the statement is a consequence of Lemma $\ref{lemred}$.\\
Assume $Q\neq Q'$ and put
\[Q=[x,m(\bfv)]_q\ ,\hspace{5mm} Q'=[x',m(\bfv)]_0\]
with $x>x'=v_u(q)\geq\tfrac{t-s}{p-1}$. Now $\Phi(Q)=[px+s-t,pm(\bfv)+s+t]_{q'}$ with $q'=ab^{-1}u^{-(t-s)}\phi(q)$.\\
If $s\neq t$, then $x'=v_u(q)>\tfrac{t-s}{p-1}$ or equivalently $x'=v_u(q)< v_u(q')=pv_u(q)-(t-s)$, and we find
\begin{align*}
d_1(Q,\Phi(Q)) &= (x-x')+(px'-(t-s)-x')+(px+s-t-(px'+s-t))\\
 &= (p+1)(x-\tfrac{t-s}{p-1})-2(x'-\tfrac{t-s}{p-1})\\
 &=(p+1)d_1(Q,P_{\rm{red}})-2d_1(Q',P_{\rm{red}}).
\end{align*}
If $s=t$, then $a\neq b$. We find $v_u(q)\neq v_u(q')$ if $v_u(q)\neq 0$ and in this case the computation is the same as above.\\
If $q=a_0+a_1u+\dots$ with $a_0\neq 0$, then $q'=ab^{-1}a_0+\dots$ and hence the absolute coefficient of $q$
is different from the absolute coefficient of $q'$. The geodesic between $Q$ and the projection of $\Phi(Q)$ 
to $\bar\B(m(\bfv))$ contains the point $Q'=[0,m(\bfv)]_0=[\tfrac{t-s}{p-1},m(\bfv)]_0$. Hence
\begin{align*}
d_1(Q,\Phi(Q))&=x+(px+s-t)=(p+1)x\\
&= (p+1)d_1(Q,P_{\rm{red}})-2d_1(Q',P_{\rm{red}}).
\end{align*}
\end{proof}
\begin{rem}
Again, this Lemma shows the connection to \cite{phimod} 6.d. The point fixed by $\Phi$ is again the point $P_{\rm{red}}$.\\
If $s=t$, then we are in the case B 2 of loc. cit.:
The projection of the fixed point to the building for $PGL_2(\bar\F((u)))$ is a vertex.
Exactly two of the half-lines of the link of this vertex are fixed by $\Phi$. \\
If $s\neq t$ we are in the case A 2 of loc. cit.:
The projection of the fixed point $P_{\rm{red}}$ is not a vertex but it lies on an edge and the projections of the two half-lines 
$\{[x,m(\bfv)]_0\mid x\leq \tfrac{t-s}{p-1}\}$ and $\{[x,m(\bfv)]_0\mid x\geq \tfrac{t-s}{p-1}\}$
to the building for $PGL_2(\bar\F((u)))$ are fixed by $\Phi$.
\end{rem}
\begin{prop}\label{propred2}
With the notations of Definition $\ref{1dimM}$ and $(\ref{MaPhia})$, $(\ref{mvred})$ assume that
\[(M_{\bar\F},\Phi)\cong (\latM^s(a)[\tfrac{1}{u}],\Phi^s_a)\oplus(\latM^t(b)[\tfrac{1}{u}],\Phi^t_b)\]
with $a,b\in\bar\F^\times$ and $0\leq s,t<p-1$. Further assume $a\neq b$ or $s\neq t$.\\
\noindent (i) The schemes $X^\bfv_{[M']}$ are empty for all $[M']\in\mathcal{S}(\bfv)\backslash\{[M_a],[M_b]\}$.\\
\noindent (ii) If $s=t$, then
\[X^\bfv_{[M_a]}\cong X^\bfv_{[M_b]}=
\begin{cases}\ \emptyset & \text{if}\ m(\bfv)+\tfrac{r_1-r_2}{p-1}\notin 2\mathbb{Z}\\
\{\ast\} & \text{if}\ m(\bfv)+\tfrac{r_1-r_2}{p-1}\in 2\mathbb{Z},
\end{cases}\]
further $X^\bfv_{[M_a]}=X^\bfv_{[M_b]}$ iff $r_1=r_2$.\\
\noindent (iii) If $s\neq t$, then 
\begin{align*}
& X^\bfv_{[M_a]}=\begin{cases} \ \emptyset & \text{if}\ \tfrac{t-s}{p-1}-\tfrac{r_1-r_2}{p-1}+m(\bfv)\notin 2\mathbb{Z}\\
 \{\ast\} & \text{if}\ \tfrac{t-s}{p-1}-\tfrac{r_1-r_2}{p-1}+m(\bfv)\in 2\mathbb{Z} \end{cases} \\
& X^\bfv_{[M_b]}=\begin{cases} \ \emptyset & \text{if}\ \tfrac{t-s}{p-1}+\tfrac{r_1-r_2}{p-1}+m(\bfv)\notin 2\mathbb{Z}\\
 \{\ast\} & \text{if}\ \tfrac{t-s}{p-1}+\tfrac{r_1-r_2}{p-1}+m(\bfv)\in 2\mathbb{Z}. \end{cases}
\end{align*}
\noindent (iv) If non empty the scheme $X^\bfv_0$ is connected.
\end{prop}
\begin{proof}
Again, this is similar to the proof of Proposition $\ref{propred1}$.\\
First, we link any $\bfv$-admissible lattice to a $\bfv$-admissible lattice in $\A_0$ by a chain of $\mathbb{P}^1$.\\
Let $\latM$ be any $\bfv$-admissible lattice and let $Q'\in\A_0\cap\bar\B(m(\bfv))$
be the unique point with $d_1(Q,Q')=d_1(Q,\A_0)$. Again, we may assume without loss of generality
$Q'=[x',m(\bfv)]_0$ with $x'\geq \tfrac{t-s}{p-1}$. Completely analogous to the proof of Proposition $\ref{propred1}$, 
we find a morphism 
\[\bar\chi:\mathbb{P}^1_{\bar\F}\rightarrow \GRloc\]
such that $\bar\chi(0)=\latM$ and $d_1(\bar\chi(\infty),\A_0)<d_1(\latM,\A_0)$.
By induction on the distance $d_1(\latM,\A_0)$, we find that we can link any $\bfv$-admissible lattice
to a $\bfv$-admissible lattice in $\A_0$.\\
Now assume $\latM=[x,m(\bfv)]_0\in\A_0$.  \\
If $x>\tfrac{t-s}{p-1}$, we find a map $\mathbb{P}^1\rightarrow \GRloc$,
as in the proof of Proposition $\ref{propred1}$, whose image contains $[x,m(\bfv)]_0$ and $[x-2,m(\bfv)]_0$, if 
$(p-1)(x-\tfrac{t-s}{p-1})\leq r_1-r_2-2.$\\
If $x<\tfrac{t-s}{p-1}$, we find a map $\mathbb{P}^1\rightarrow \GRloc$
whose image contains $[x,m(\bfv)]_0$ and $[x+2,m(\bfv)]_0$, if 
$(p-1)(\tfrac{t-s}{p-1}-x)\leq r_1-r_2-2.$\\
Similarly as in Proposition $\ref{propred1}$, one can proof the following two facts:\\
\noindent \emph{(a)} The set
\[\mathcal{N}=\{\latM\in\GRloc(\bar\F)\mid \latM\notin\A_0\ \text{or}\ 
(\latM\in\A_0\ \text{and}\ d_1(\latM,P_{\rm{red}})\leq \tfrac{r_1-r_2-2}{p-1})\}\]
consists of non-$\bfv$-ordinary lattices.\\
\noindent \emph{(b)} If $\latM\in\A_0$ is $\bfv$-admissible, then
\[d_1(\latM,P_{\rm{red}})>\tfrac{r_1-r_2-2}{p-1}\Rightarrow 
d_1(\latM,P_{\rm{red}})=\tfrac{r_1-r_2}{p-1}.\]
Now $\mathcal{N}$ defines a connected subset of $X^\bfv_0$: \\
If $s\neq t$, then there is a unique lattice with minimal distance from $[\tfrac{t-s}{p-1},m(\bfv)]_0$
and every $\bfv$-admissible lattice in $\mathcal{N}$ can be linked to this lattice by a chain of $\mathbb{P}^1$.\\
If $s=t$, then either $[\tfrac{t-s}{p-1},m(\bfv)]_0$ is a lattice itself and any $\bfv$-admissible lattice can
be linked to this lattice by a chain of $\mathbb{P}^1$,
or there are two $\bfv$-admissible lattices $[\pm 1,m(\bfv)]_0$ in $\mathcal{N}$ with distance $1$ from $[0,m(\bfv)]_0$
and by the above there is a morphism 
\[\mathbb{P}^1\rightarrow \GRloc\] 
containing both lattices in its image. Thus $\mathcal{N}$ defines a connected subset.\\
Consider the following points:
\begin{align*}
& Q_+=[\tfrac{t-s}{p-1}+\tfrac{r_1-r_2}{p-1},m(\bfv)]_0\\
& Q_-=[\tfrac{t-s}{p-1}-\tfrac{r_1-r_2}{p-1},m(\bfv)]_0.
\end{align*}
We are left to show that, if one of these points defines a lattice $\latM_+=Q_+$ (resp. $\latM_-=Q_-$),
then this point lies in $X^\bfv_{[M_b]}$ (resp $X^\bfv_{[M_a]}$).\\
If $\tfrac{t-s}{p-1}-\tfrac{r_1-r_2}{p-1}+m(\bfv)\in 2\mathbb{Z}$, then $Q_-=\latM_-$ is a lattice
and 
\[e_{\latM_-}=u^{(e-r_1-s)/(p-1)}e_1\] 
defines a $\Phi$-stable subspace satisfying $\Phi(e_{\latM_-})=au^{e-r_1}e_{\latM_-}$,
i.e. $\latM_-\in X^\bfv_{[M_a]}$.\\
If $\tfrac{t-s}{p-1}+\tfrac{r_1-r_2}{p-1}+m(\bfv)\in 2\mathbb{Z}$, then $Q_+=\latM_+$ is a lattice
and 
\[e_{\latM_+}=u^{(e-r_1-t)/(p-1)}e_2\] 
defines a $\Phi$-stable subspace satisfying $\Phi(e_{\latM_+})=bu^{e-r_1}e_{\latM_+}$,
i.e. $\latM_+\in X^\bfv_{[M_b]}$.\\
We have two different cases: \\
If $s=t$ and $\tfrac{r_1-r_2}{p-1}+m(\bfv)\in 2\mathbb{Z}$, then the lattices 
$\latM_+$ and $\latM_-$ define points $\latM_-\in X^\bfv_{[M_a]}$ and
$\latM_+\in X^\bfv_{[M_b]}$ which coincide iff $r_1=r_2$.\\
If $s\neq t$, then   
\begin{align*}
& Q_-\ \text{defines an isolated point in}\ X^\bfv_{[M_a]}
\Leftrightarrow \tfrac{t-s-(r_1-r_2)}{p-1}+m(\bfv)\in 2\mathbb{Z}\\
& Q_+\ \text{defines an isolated point in}\ X^\bfv_{[M_b]}
\Leftrightarrow \tfrac{t-s+(r_1-r_2)}{p-1}+m(\bfv)\in 2\mathbb{Z}.
\end{align*}
This cannot happen at the same time, as $\tfrac{t-s}{p-1}\notin\mathbb{Z}$. This finishes the proof of the Proposition.
\end{proof}
\begin{cor}\label{oneptred2}
\noindent (i) Assume $s=t$. \\
\noindent (a) If $m(\bfv)\equiv 0\mod 2$, then $\GRloc=\{\ast\}$ iff $\tfrac{r_1-r_2}{p-1}<2$.\\
\noindent (b) If $m(\bfv)\equiv 1\mod 2$, then $\GRloc$ cannot be a single point.
\begin{align*}
& \GRloc=\emptyset && \Leftrightarrow 0\leq\tfrac{r_1-r_2}{p-1}<1\\
& \GRloc=\{\ast\}\cup\{\ast\} && \Leftrightarrow 1\leq\tfrac{r_1-r_2}{p-1}<3.
\end{align*}
\noindent (ii) Assume $s\neq t$. \\ 
Define $x_0=\lfloor\tfrac{t-s}{p-1}\rfloor$ and write $\xi=\tfrac{t-s}{p-1}-x_0$
for the fractional part of $\tfrac{t-s}{p-1}$.
\[\GRloc=\{\ast\}\Leftrightarrow\begin{cases}
0+\xi\leq \tfrac{r_1-r_2}{p-1}<2-\xi&if\ m(\bfv)\equiv x_0\mod 2\\
1-\xi\leq \tfrac{r_1-r_2}{p-1}<1+\xi&if\ m(\bfv)\not\equiv x_0\mod 2.
\end{cases}\]
\end{cor}
\begin{proof}
\noindent \emph{(i)} This is nearly identical to Corollary $\ref{oneptred1}$.\\
\emph{(ii)} Assume $m(\bfv)\equiv x_0\mod 2$. Then $[x_0,m(\bfv)]_0$ is the unique lattice with minimal distance $d_1$
from $P_{\rm{red}}$. By Lemma $\ref{lemred2}$ it is $\bfv$-admissible iff $\tfrac{r_1-r_2}{p-1}\geq \xi$.\\
Again by Lemma $\ref{lemred2}$ it is the only $\bfv$-admissible lattice iff $[x_0+2,m(\bfv)]_0$
is not $\bfv$-admissible. This is the case iff $\tfrac{r_1-r_2}{p-1}<2-\xi$.\\
The case $m(\bfv)\not\equiv x_0\mod 2$ is similar.
\end{proof}
\subsection{The case of a non split extension}
Finally, we analyze the case where $(M_{\bar\F},\Phi)$ is a non split extension of two one dimensional objects. 
There is a basis $e_1,e_2$ such that
\[M_{\bar\F}\sim\begin{pmatrix}au^s&\gamma\\0&bu^t\end{pmatrix}\]
with $0\leq s,t<p-1$ and $a,b\in \bar\F^\times$, $\gamma\in\bar\F((u))$. 
In any basis of the form $e_1,qe_1+e_2$ defining the apartment $\A_q$, the 
endomorphism $\Phi$ is upper triangular with diagonal entries $au^s$ and $bu^t$, 
and we fix the basis such that the valuation of the upper right entry $k:=v_u(\gamma)$ is maximal. 
\begin{lem}\label{lemred3}
\noindent (i) The integer $k=v_u(\gamma)$ satisfies
\[k\leq \tfrac{pt-s}{p-1}.\]
\noindent (ii) If $\latM=[x,y]_q$ with $\min\{x,v_u(q)\}\geq\tfrac{k-s}{p}$, then
\begin{align*}
& d_1(\latM,\langle\Phi(\latM)\rangle)=(p+1)x+s+t-2k\\
& d_2(\latM,\langle\Phi(\latM)\rangle)=(p-1)d_2(\latM,P_{\rm{red}}).
\end{align*}
\noindent (iii) If $\latM=[x,y]_q$ with $x<\tfrac{k-s}{p}$ or $v_u(q)<\tfrac{k-s}{p}$, let $Q'\in\A_0\cap\bar\B(y)$ be the 
unique point such that $d_1(\latM,Q')=d_1(\latM,\A_0)$. Then
\begin{align*}
& d_1(\latM,\langle\Phi(\latM)\rangle)=(p+1)d_1(\latM,P_{\rm{red}})-2d_1(Q',P_{\rm{red}})\\
& d_2(\latM,\langle\Phi(\latM)\rangle)=(p-1)d_2(\latM,P_{\rm{red}}).
\end{align*}
\end{lem}
\begin{proof}
\noindent \emph{(i)} This follows from the maximality of $k=v_u(\gamma)$:
We have 
\begin{equation}\label{gammabw}\Phi(qe_1+e_2)=(\gamma+au^s\phi(q)-bu^tq)e_1+bu^t(qe_1+e_2).\end{equation}
And 
\[v_u(au^s\phi(q)-bu^tq)=\begin{cases}v_u(q)+t & \text{if}\ v_u(q)>\tfrac{t-s}{p-1}\\
pv_u(q)+s & \text{if}\ v_u(q)<\tfrac{t-s}{p-1}.\end{cases}\]
If we had $k=v_u(\gamma)=v_u(q)+t$ for any $q$ with $v_u(q)>\tfrac{t-s}{p-1}$, we
could delete the leading coefficient of $\gamma$ in $(\ref{gammabw})$ which contradicts the maximality of $v_u(\gamma)$.
Hence we have $k<v_u(q)+t$ for all $q$ with $v_u(q)>\tfrac{t-s}{p-1}$ which yields the first claim.\\
\noindent \emph{(ii)} The first part of the lemma implies $\tfrac{k-s}{p}\geq k-t$ and hence our assumptions on $v_u(q)$ imply 
$k\leq \min\{v_u(q)+t,pv_u(q)+s\}$. We find $v_u(\gamma+au^s\phi(q)-bu^tq)=v_u(\gamma)$ and we may assume $q=0$, i.e. $\latM\in\A_0$,
as the situation is the same as in the standard apartment.
Now we have $\langle\Phi(\latM)\rangle=[px+s-t,py+s+t]_{b^{-1}u^{-t}\gamma}$, and $x\geq \tfrac{k-s}{p}$ implies
\[px+s-t\geq k-t\ ,\hspace{5mm} x\geq k-t.\]
Thus $d_1(\latM,\langle\Phi(\latM)\rangle)=(px+s-t-(k-t))+(x-(k-t))=(p+1)x+s+t-2k$. The statement on $d_2$ is easy.\\
\noindent \emph{(iii)} If $\latM\notin\A_0$, then $v_u(q)<\tfrac{k-s}{p}\leq \tfrac{t-s}{p-1}$ and hence
\[v_u(\gamma+au^s\phi(q)-bu^tq)=v_u(au^s\phi(q)-bu^tq)\]
and the situation is the same as in the split case, i.e. the case $\gamma=0$.
If $\latM\in\A_0$, then $\langle\Phi(\latM)\rangle\in\A_0$ and the statement is easy.
\end{proof}
\begin{rem}
In the case of a non split extension we are in the case B 2 or A 3 of \cite{phimod} 6.d.
More precisely, if $\tfrac{k-s}{p}\notin\mathbb{Z}$, then the unique fixed point of (\cite{phimod}, Prop. 6.1)
is not in the building $\bar\B$. It is only visible after extending $\bar\F((u))$ to some separable wildly ramified extension
(the apartment containing the fixed point will branch off from $\A_0$ at the line $x=\tfrac{k-s}{p}$,
because we can successively delete the leading coefficient of $\gamma$ in $(\ref{gammabw})$ if there is 
some $q$ with $v_u(q)=\tfrac{k-s}{p}$).
The image of the half line $\{[x,m(\bfv)]_0\mid x\leq \tfrac{k-s}{p}\}$ in the building for $PGL_2(\bar\F((u)))$ is stable under $\Phi$
and the geodesic between $[\lfloor\tfrac{k-s}{p}\rfloor+1,m(\bfv)]_0$ and its image under $\Phi$
contains the (projection of the) point $[x_0,m(\bfv)]_0$
in the building for $PGL_2(\bar\F((u)))$.
This is the case A 3 of \cite{phimod} 6.d.\\
If $\tfrac{k-s}{p}\in\mathbb{Z}$, then we are in the case B 2 of \cite{phimod} 6.d.:
In this case the maximality of $k=v_u(q)$ implies $k-t=\tfrac{k-s}{p}=\tfrac{t-s}{p-1}=0$
(otherwise we could delete the leading coefficient of $\gamma$) and we find that
$P_{\rm{red}}$ is the fixed point in the building.
In this case there is a unique half-line in the apartment for $PGL_2(\bar\F((u)))$ that is fixed by $\Phi$,
namely the image of the half-line $\{[x,m(\bfv)]_0\mid x\leq 0\}$ under the projection.
\end{rem}
\begin{prop}\label{propred3}
With the notations of Definition $\ref{1dimM}$ and $(\ref{MaPhia})$, $(\ref{mvred})$, assume that $(M_{\bar\F},\Phi)$ is a non split 
extension 
\[0\rightarrow (\latM^{s}(a)[\tfrac{1}{u}],\Phi_a^s) \rightarrow (M_{\bar\F},\Phi) \rightarrow
(\latM^{t}(b)[\tfrac{1}{u}],\Phi_b^t) \rightarrow 0\]
for some $a,b\in\bar\F^\times$ and $0\leq s,t<p-1$.\\
\noindent (i) The schemes $X^\bfv_{[M']}$ are empty for all $[M']\in\mathcal{S}(\bfv)\backslash\{[M_a]\}$.\\
\noindent (ii) For $X^\bfv_{[M_a]}$ the following holds:
\[X^\bfv_{[M_a]}=\begin{cases}
\ \emptyset & \text{if}\ \tfrac{t-s}{p-1}-\tfrac{r_1-r_2}{p-1}+m(\bfv)\notin 2\mathbb{Z} \\
\{\ast\}    & \text{if}\ \tfrac{t-s}{p-1}-\tfrac{r_1-r_2}{p-1}+m(\bfv)\in 2\mathbb{Z}.
\end{cases}\]
\noindent (iii) If non empty, the scheme $X^\bfv_0$ is connected.
\end{prop}
\begin{proof}
Lemma $\ref{lemred3}$ \emph{(i)} implies $\tfrac{k-s}{p}\leq \tfrac{t-s}{p-1}$ and an easy computation using the
same inequality shows that
\[\tfrac{t-s}{p-1}-\tfrac{r_1-r_2}{p-1}\leq \tfrac{k-s}{p}\Leftrightarrow 
\tfrac{k-s}{p}\leq \tfrac{1}{p+1}(r_1-r_2-s-t+2k),\]
and hence
\[\tfrac{t-s}{p-1}-\tfrac{r_1-r_2}{p-1}\leq \tfrac{k-s}{p}\leq \tfrac{1}{p+1}(r_1-r_2-s-t+2k),\]
if $\GRloc\neq \emptyset$. Further denote by $\widetilde{\mathcal{N}}$ the set of $\bfv$-admissible lattices
$\latM=[x,m(\bfv)]_q$ with $\tfrac{k-s}{p}\leq \min\{x,v_u(q)\}$. \\
As we have seen above, the situation for the $\bfv$-admissible lattices $\latM\notin\widetilde{\mathcal{N}}$
is the same as in the split case.
Hence we can link all $\bfv$-admissible lattices to $\bfv$-admissible lattices in $\A_0$ by a chain of $\mathbb{P}^1$.
If $\latM=[x,m(\bfv)]_0$ is a $\bfv$-admissible lattice in $\A_0$ with $x+2<\tfrac{k-s}{p}$,
then there is a $\mathbb{P}^1$ in $\GRloc$ containing $\latM=[x,m(\bfv)]_0$ and $[x+2,m(\bfv)]_0$ 
, except if $\latM=\latM_-=[\tfrac{t-s-(r_1-r_2)}{p-1},m(\bfv)]_0$ which defines an isolated point in $X^\bfv_{[M_a]}$
if $\tfrac{t-s}{p-1}-\tfrac{r_1-r_2}{p-1}+m(\bfv)\in 2\mathbb{Z}$ (compare Proposition $\ref{propred2}$).\\
Let $\latM'=[x_0,m(\bfv)]_0$ be the lattice where $x_0$ is the maximal integer smaller than $\tfrac{k-s}{p}$
that is congruent to $m(\bfv)\mod 2$. We claim:\\
\emph{\noindent (a) If $\latM'$ is $\bfv$-admissible  and $x_0\neq \tfrac{t-s}{p-1}-\tfrac{r_1-r_2}{p-1}$,
then any lattice in $\widetilde{\mathcal{N}}$ can be linked to $\latM'$ by a chain of $\mathbb{P}^1$. \\
\noindent (b) The lattices $\latM\in\widetilde{\mathcal{N}}$ are non-$\bfv$-ordinary.}\\
This finishes the proof of the proposition.\\ 
{\emph{Proof of (a)}}: Let $\latM=[x,m(\bfv)]_q\in\widetilde{\mathcal{N}}$ be a lattice. Without loss of generality, we may assume
$\latM\in\A_0$, as the situation is the same in all apartments $\A_q$ with $v_u(q)\geq \tfrac{k-s}{p}$.
By Lemma $\ref{lemred3}$, we have
\[\tfrac{k-s}{p}\leq x\leq \tfrac{1}{p+1}(r_1-r_2-s-t+2k).\]
We consider the basis
\[b_1=u^{(x+m(\bfv))/2}e_1\ ,\ b_2=u^{(m(\bfv)-x)/2}e_2\]
of $\latM$ and by Lemma $\ref{P1nachgrass}$, there is a morphism 
\[\bar\chi:\mathbb{P}^1_{\bar\F}\rightarrow {\rm{Grass}}\,M_{\bar\F}\]
with $\bar\chi(z)=[x,m(\bfv)]_{z^{x-1}}$ for $z\in\bar\F$ and $\bar\chi(\infty)=[x-2,m(\bfv)]_0$.\\
If $x-1\geq \tfrac{k-s}{p}$, then the morphism factors over $\GRloc$.
Consider the following two cases:\\
If $\tfrac{k-s}{p}\leq x_0+1$, then this argument shows that we can link all $\latM\in\widetilde{\mathcal{N}}$
to the lattice $[x_0,m(\bfv)]_0$ by a chain of $\mathbb{P}^1$.\\
If $\tfrac{k-s}{p}>x_0+1$, then this argument shows that we can link all $\latM\in\widetilde{\mathcal{N}}$
to the lattice $\latM''=[x_0+2,m(\bfv)]_0$ by a chain of $\mathbb{P}^1$.
We can link the lattice $\latM''$ to the lattice $\latM'=[x_0,m(\bfv)]_0$ if the lattices $\latM_z=[x_0,m(\bfv)]_{zu^{x-1}}$
are $\bfv$-admissible for all $z\in\bar\F$. For $z\neq 0$ we have
\[d_1(\latM_z,\langle\Phi(\latM_z)\rangle)=(p+1)d_1(\latM_z,P_{\rm{red}})-2d_1(Q',P_{\rm{red}}),\]
where $Q'=[x_0+1,m(\bfv)]_0$ is the unique point in $\A_0$ with minimal distance from $\latM_z$. Hence
$d_1(\latM_z,\langle\Phi(\latM_z)\rangle)=d_1(\latM',\langle\Phi(\latM')\rangle)+2$
and the morphism factors through $\GRloc$ as $x_0\neq \tfrac{t-s}{p-1}-\tfrac{r_1-r_2}{p-1}$.
(Otherwise $\latM'$ is the unique isolated point in $X^\bfv_{[M_a]}$).\\
{\emph{Proof of (b)}}: Let $\latM\in\widetilde{\mathcal{N}}$ be a lattice. 
Similarly to the proof of Proposition $\ref{propred1}$, we find
\[\latM \sim\begin{pmatrix} a_{11} & a_{12}\\ 0 & a_{22}\end{pmatrix}\]
with $v_u(a_{12})<v_u(a_{11})$ and hence the minimal elementary divisor of $\langle\Phi(\latM)\rangle$
with respect to $\latM$ is not defined by a $\Phi$-stable subspace.
\end{proof}
Summarizing the results on the connected components we find the following Theorem.
\begin{theo}\label{theored}
Assume that $(M_{\F},\Phi)$ becomes reducible after extending the scalars to some finite extension $\F'$ of $\F$.\\
\noindent (i) The subschemes $X^\bfv_0$ and $X^\bfv_{[M']}$ are open and closed in $\GRloc\otimes_\F\bar\F$ for all 
isomorphism classes $[M']\in\mathcal{S}(\bfv)$.\\
\noindent (ii) If non empty, the scheme $X^\bfv_0$ is connected.\\
\noindent (iii) For each $[M']\in\mathcal{S}(\bfv)$ the scheme $X^\bfv_{[M']}$ is either empty, a single point or
isomorphic to $\mathbb{P}^1_{\bar\F}$.\\
\noindent (iv) There are at most two isomorphism classes $[M']\in\mathcal{S}(\bfv)$ such that $X^\bfv_{[M']}\neq \emptyset$.
\end{theo} 
\begin{proof}
This is a summary of the Propositions $\ref{propred1}$, $\ref{propred2}$ and $\ref{propred3}$.
\end{proof}
This theorem implies a modified version of the conjecture of Kisin stated in (\cite{Kisin}, 2.4. 16).
\begin{defn}\label{defconcomp}
For an integer $s$ denote by $\mathcal{GR}_{V_\F,0}^{\bfv,{\rm{loc}},s}$ the open and closed subscheme of $\GRloc$
consisting of all $\bfv$-admissible lattices $\latM$, where the rank of the maximal $\Phi$-stable subobject
$\latM_1$ satisfying $\langle\Phi(\latM_1)\rangle=u^{e-r_1}\latM_1$ is equal to $s$.
\end{defn}
\begin{cor}\label{concomp}
Assume $p\neq 2$ and let $\rho:G_K\rightarrow V_\F$ be any two-dimensional continuous representation of $G_K$
that admits a finite flat model after possibly extending the scalars to some finite extension of $\F$.\\
Assume that ${\rm{End}}_{\F'[G_K]}(V_{\F'})$ is a simple algebra for all finite extensions $\F'$ of $\F$.
Then $\mathcal{GR}_{V_\F,0}^{\bfv,{\rm{loc}},s}$ is geometrically connected for all $s$.
Furthermore \\
\noindent (i) If $s=1$ and ${\rm{End}}_{\F'[G_K]}(V_{\F'})=\F'$ for all finite extensions $\F'$ of $\F$, 
then $\mathcal{GR}_{V_\F,0}^{\bfv,{\rm{loc}},s}$ is either empty or a single point.\\
If $s=1$ and ${\rm{End}}_{\F'[G_K]}(V_{\F'})=M_2(\F')$ for some finite extension $\F'$ of $\F$,
then $\mathcal{GR}_{V_\F,0}^{\bfv,{\rm{loc}},s}$ is either empty or becomes isomorphic to $\mathbb{P}^1_{\F'}$ after extending the
scalars to $\F'$.\\
\noindent (ii) If $s=2$, then $\mathcal{GR}_{V_\F,0}^{\bfv,{\rm{loc}},s}$ is either empty or a single point.
\end{cor}
\begin{proof}
Our definitions imply
\[\mathcal{GR}_{V_\F,0}^{\bfv,{\rm{loc}},0}\otimes_\F\bar\F=X^\bfv_0.\]
Further 
\[\bigcup_{[M']\in\mathcal{S}(\bfv)}X^\bfv_{[M']}=\begin{cases}
\mathcal{GR}_{V_\F,0}^{\bfv,{\rm{loc}},1} & \text{if}\ r_1>r_2\\
\mathcal{GR}_{V_\F,0}^{\bfv,{\rm{loc}},2} & \text{if}\ r_1=r_2.
\end{cases}\]
By (\cite{Breuil}, Thm. 3.4.3) we have ${\rm{End}}_{\F'[G_K]}(V_{\F'})= {\rm{End}}_{\F'((u)),\Phi}(M_{\F'})$.
The same Theorem implies that the image of the category of finite flat $G_K$-representations
on finite length $\mathbb{Z}_p$-algebras under the restriction to $G_{K_\infty}$ is closed under subobjects and quotients. 
Hence $V_{\F'}$ is irreducible (resp. reducible, resp. split reducible) if and only if $(M_{\F'},\Phi)$ is.
An easy computation yields:\\
${\rm{End}}_{\F'[G_K]}(V_{\F'})=\F'$ if $V_{\F'}$ is irreducible or non-split reducible.\\
${\rm{End}}_{\F'[G_K]}(V_{\F'})=\F'\times \F'$ if $V_{\F'}$ is the direct sum of two non-isomorphic 
one-dimensional representations.\\
${\rm{End}}_{\F'[G_K]}(V_{\F'})=M_2(\F')$ if $V_{\F'}$ is the direct sum of two isomorphic 
one-dimensional representations.\\
The Corollary now follows from Theorem $\ref{theored}$ and Propositions $\ref{propred1}$, $\ref{propred2}$ and $\ref{propred3}$.
\end{proof}
\section{The structure of $X^\bfv_0$}
In this section we want to analyze the structure of the connected component $X^\bfv_0$ of non-$\bfv$-ordinary lattices.
In the absolutely simple case we have 
\[X_0^\bfv=\GRloc\otimes_\F\bar\F\]
and this is isomorphic to a Schubert variety.
In the reducible case it turns out that this component has a quite complicated structure. 
It is in general not irreducible and its irreducible
components have varying dimensions.
\subsection{The case $(M_{\bar\F},\Phi)\cong (M_1,\Phi_1)\oplus (M_1,\Phi_1)$.}
We assume that $(M_{\bar\F},\Phi)$ is a direct summand of two isomorphic one-dimensional objects and we will use the notations
of section $4.1$.
First we define some subsets of the affine Grassmannian.\\
Denote by $n$ the maximal integer congruent to $m(\bfv)\mod 2$, such that 
\begin{equation}\label{n1}
n\leq \tfrac{r_1-r_2+2}{p+1}.
\end{equation}
Denote by $l$ the minimal integer such that 
\begin{equation}\label{l1}
n+2\leq \tfrac{r_1-r_2+2l}{p+1}.
\end{equation}
For $z\in\mathbb{P}^1(\bar\F)$ and $j\geq 0$ we define the following points:
\begin{align*}
& Q_j^z\ =[l+(p+1)j,m(\bfv)]_z \hspace{5mm} \text{if}\ z\in\bar\F\\
& Q_j^\infty=[-l-(p+1)j,m(\bfv)]_0.
\end{align*}
We define the following subschemes $Z,Z_j\subset X^\bfv_0$ for $j\geq 0$ by specifying its closed points:
\begin{equation}\label{Zred1}
\begin{aligned}
& Z(\bar\F)\ =\{\latM\in\bar\B(m(\bfv))\mid d_1(\latM,P_{\rm{red}})\leq n\}\\
& Z_j(\bar\F)=\bigcup_{z\in\mathbb{P}^1(\bar\F)}\{\latM\in\bar\B(m(\bfv))\mid d_1(\latM,Q_j^z)\leq n+2-l-(p-1)j\}.
\end{aligned}
\end{equation}
We want to consider these subsets as subschemes with the reduced scheme structure.
\begin{lem}\label{lemX01}
With the notations of $(\ref{n1})$-$(\ref{Zred1})$: 
\[X^\bfv_0(\bar\F)=(\bigcup_{j\geq 0} Z_j(\bar\F))\cup Z(\bar\F).\]
\end{lem}
\begin{proof}
Let $\latM=[x,m(\bfv)]_q$ be a non-$\bfv$-ordinary lattice and denote by $Q'=[x',m(\bfv)]_z$ the unique lattice with 
$d_1(\latM,Q')=d_1(\latM,\T)$. Without loss of generality, we may assume that
$Q'\in\calL_0$, i.e. $z=0$ and $v_u(q)=x'>0$.\\
If $1\leq x'=v_u(q)<l$, then by Lemma $\ref{distred1}$ and the definition of $n$ and $l$ we find that
$\latM$ is $\bfv$-admissible (and non-$\bfv$-ordinary) iff $d_1(\latM,P_{\rm{red}})\leq n$.\\
If $v_u(q)\geq l$, then there is a unique $j$ such that $l+(p+1)j\leq x'<l+(p+1)(j+1)$. By Lemma $\ref{distred1}$, we find
that $\latM$ is $\bfv$-admissible iff $d_1(\latM,P_{\rm{red}})\leq \tfrac{r_1-r_2+2x'}{p+1}$. Now 
\[d_1(\latM,P_{\rm{red}})=d_1(\latM,Q_j^0)+(l+(p+1)j)\]
and hence $\latM$ is $\bfv$-admissible iff
\[d_1(\latM,Q_j^0)\leq \tfrac{r_1-r_2+2l}{p+1}+\tfrac{2(x'-l-(p+1)j)}{p+1}-l-(p-1)j.\]
By the definition of $n$ and $l$ and the fact $x'-l-(p+1)j<(p+1)$ we find that $\latM$ is $\bfv$-admissible iff
\[d_1(\latM,Q_j^0)\leq n+2-l-(p-1)j.\]
This yields the claim.
\end{proof}
\begin{prop}\label{propX01}
With the notation of $(\ref{n1})$-$(\ref{Zred1})$:
\begin{equation}\label{decompX01}
X^\bfv_0=(\bigcup_{j\geq 0} Z_j)\cup Z.
\end{equation}
\noindent (i) The scheme $Z$ is isomorphic to an $n$-dimensional Schubert variety.\\
\noindent (ii) For $j\geq 0$ there is a projective, surjective and birational morphism
\[f_j: \mathbb{P}^1_{\bar\F}\times Y_j\rightarrow Z_j\]
where $Y_j$ is an $n+2-l-(p-1)j$ dimensional Schubert variety. Especially $Z_j$ is closed and irreducible.\\
\noindent (iii) If $l\neq 2$, then $(\ref{decompX01})$ is the decomposition of $X^\bfv_0$ into its irreducible components.\\ 
\noindent (iv) If $l=2$, then the decomposition of $X^\bfv_0$ into its irreducible components is given by
\[X^\bfv_0=\bigcup_{j\geq 0}Z_j.\]
\noindent (v) The dimension of $X_0^\bfv$ is given by
\[\dim\, X^\bfv_0=\begin{cases}n+1&\text{if}\ l=2\\ n& \text{if}\ l\neq 2.\end{cases}\]
\end{prop}
\begin{proof}
\noindent \emph{(i)} The closed points of the scheme $Z$ are the lattices with distance smaller than $n$ from the 
point $[0,m(\bfv)]_0$.
By the same argument as in the proof of Theorem $\ref{theoirred}$ \emph{(b)}, this is an $n$-dimensional Schubert variety.\\
\noindent \emph{(ii)} The scheme $Z_j$ is the union of the Schubert varieties consisting of the lattices $\latM$ with distance $d_1(\latM,Q_j^z)\leq n+2-l-(p-1)j=:n_j$ for $z\in\mathbb{P}^1(\bar\F)$.\\
Let us first assume that $m(\bfv)\equiv x_j:=l+(p+1)j\mod 2$, i.e. $Q_j^z$ is a lattice for all $z\in\mathbb{P}^1(\bar\F)$.
For any linearly independent vectors $b_1$ and $b_2$ denote by
\[\psi(b_1,b_2):Y_j\hookrightarrow {\rm{Grass}}\,M_{\bar\F}\]
the inclusion of the Schubert variety of lattices $\latM$ with
\begin{align*}
& d_1(\latM,\langle b_1,b_2\rangle)\leq n_j \\
& d_2(\latM,\langle b_1,b_2\rangle)=0 .
\end{align*}
First we construct a morphism
\[\widetilde{f}_j:\mathbb{P}^1_{\bar\F}\times Y_j\rightarrow {\rm{Grass}}\,M_{\bar\F}.\]
The inclusion $\psi(e_1,e_2)$ defines a sheaf $\mathcal{M}_{Y_j}$ of $\mathcal{O}_{Y_j}[\hspace{-0.5mm}[u]\hspace{-0.5mm}]$-lattices in $M_{\bar\F}\widehat{\otimes}_{\bar\F}\mathcal{O}_{Y_j}$.
If $U={\rm{Spec}}\,A\subset Y_j$ is an affine open we write $\latM_A=\Gamma(U,\mathcal{M}_{Y_j})$
for the $A[\hspace{-0.5mm}[u]\hspace{-0.5mm}]$-lattice in $M_{\bar\F}\widehat{\otimes}_{\bar\F}A$ defined by $\mathcal{M}_{Y_j}$. 
To define the morphism $\widetilde{f}_j$ we define a sheaf $\widetilde{\mathcal{M}}$ of $\mathcal{O}_{\mathbb{P}^1_{\bar\F}\times Y_j}[\hspace{-0.5mm}[u]\hspace{-0.5mm}]$-lattices in $M_{\bar\F}\widehat{\otimes}_{\bar\F}\mathcal{O}_{\mathbb{P}^1_{\bar\F}\times Y_j}$.
Let $\mathbb{P}^1_{\bar\F}=V_0\cup V_\infty$ with $V_0={\rm{Spec}}\,\bar\F[T]$ and $V_\infty={\rm{Spec}}\,\bar\F[T^{-1}]$. We define $\widetilde{\mathcal{M}}$ by specifying its sections over the open subsets $V\times U$ of $\mathbb{P}^1_{\bar\F}\times Y_j$ where $V\subset \mathbb{P}^1_{\bar\F}$ and $U={\rm{Spec}}\,A\subset Y_j$ are affine open subschemes.
If $V'={\rm{Spec}}\,\bar\F[T]_g\subset V_0$ for some $g\in\bar\F[T]$, then $\Gamma(V'\times U,\widetilde{\mathcal{M}})$ is the pushout of $\latM_A\widehat{\otimes}_A A[T]_g$ via the endomorphism of $M_{\bar\F}\widehat{\otimes}_{\bar\F}A[T]_g$ defined by the matrix
\[C_{A}^0=\begin{pmatrix}u^{(m(\bfv)+x_j)/2}& T u^{(m(\bfv)-x_j)/2}\\ Tu^{(m(\bfv)+x_j)/2} & u^{(m(\bfv)-x_j)/2}\end{pmatrix}.\]
If $V''={\rm{Spec}}\,\bar\F[T^{-1}]_h\subset V_\infty$ for some $h\in\bar\F[T^{-1}]$, then $\Gamma(V''\times U,\widetilde{\mathcal{M}})$ is the pushout of $\latM_A\widehat{\otimes}_A A[T^{-1}]_h$ via the endomorphism of $M_{\bar\F}\widehat{\otimes}_{\bar\F}A[T^{-1}]_h$ defined by the matrix
\[C_{A}^\infty=\begin{pmatrix}T^{-1}u^{(m(\bfv)+x_j)/2}&  u^{(m(\bfv)-x_j)/2}\\ u^{(m(\bfv)+x_j)/2} & T^{-1}u^{(m(\bfv)-x_j)/2}\end{pmatrix}.\]
These definitions are compatible: if $V'\subset V_0\cap V_\infty={\rm{Spec}}\,\bar\F[T,T^{-1}]$,
then the matrices $C_A^0$ and $C_A^\infty$ differ by a unit (namely $T$ resp. $T^{-1}$).
Further this definitions are compatible with localization in the following sense.\\
If $U'={\rm{Spec}}\,B\subset U={\rm{Spec}}\,A$ is an affine open, then
\[\Gamma(V'\times U',\widetilde{\mathcal{M}})=\Gamma(V'\times U,\widetilde{\mathcal{M}})\widehat{\otimes}_A B,\]
as $\latM_B=\latM_A\widehat{\otimes}_AB$.
And similarly for $V''$ and for localization on $\mathbb{P}^1_{\bar\F}$.
As the sets $\{V'\times U,V''\times U\mid V'\subset V_0 \ ,\ V''\subset V_\infty\ ,\  U\subset Y_j\ \text{affine open}\}$ form a basis of the topology this indeed defines a sheaf of $\mathcal{O}_{\mathbb{P}^1_{\bar\F}\times Y_j}[\hspace{-0.5mm}[u]\hspace{-0.5mm}]$-lattices on $\mathbb{P}^1_{\bar\F}\times Y_j$.\\

\begin{figure}[h]
\vspace{3cm}
\begin{picture}(225,-100)
\drawline(0,0)(200,0)
\drawline(100,0)(175,75)
\drawline(100,0)(175,-75)
\drawline(150,0)(170,20)\drawline(150,0)(170,-20)
\drawline(160,0)(170,7)\drawline(160,0)(170,-7)
\drawline(160,10)(170,25)\drawline(160,10)(170,15)
\drawline(160,-10)(170,-25)\drawline(160,-10)(170,-15)
\drawline(140,0)(130,7)\drawline(140,0)(130,-7)

\drawline(50,0)(30,20)\drawline(50,0)(30,-20)
\drawline(40,0)(30,7)\drawline(40,0)(30,-7)
\drawline(60,0)(70,7)\drawline(60,0)(70,-7)
\drawline(40,10)(30,25)\drawline(40,10)(30,15)
\drawline(40,-10)(30,-25)\drawline(40,-10)(30,-15)

\drawline(150,50)(145,75)\drawline(150,50)(175,40)
\drawline(148,62)(140,75) \drawline(148,62)(150,75)
\drawline(162,45)(175,33)\drawline(162,45)(175,47)
\drawline(160,60)(160,70)\drawline(160,60)(170,60)
\drawline(140,40)(140,30)\drawline(140,40)(130,40)

\drawline(150,-50)(145,-75)\drawline(150,-50)(175,-40)
\drawline(148,-62)(140,-75) \drawline(148,-62)(150,-75)
\drawline(162,-45)(175,-33)\drawline(162,-45)(175,-47)
\drawline(160,-60)(160,-70)\drawline(160,-60)(170,-60)
\drawline(140,-40)(140,-30)\drawline(140,-40)(130,-40)

\jput(150,0){\makebox(0,0){$\bullet$}}
\jput(50,0){\makebox(0,0){$\bullet$}}
\jput(150,50){\makebox(0,0){$\bullet$}}
\jput(150,-50){\makebox(0,0){$\bullet$}}
\end{picture}
\vspace{3cm}
\caption{The closed points of $Z_j$ in the building in the case $p=3$, $\F=\F_3$.
The fat points mark the points $Q_j^z$ for $z\in\mathbb{P}^1(\F).$}
\end{figure}
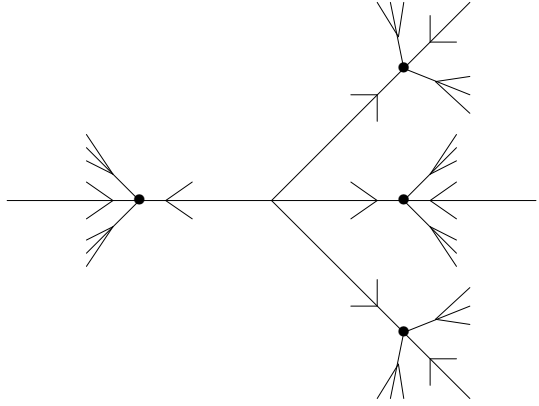
By construction the values of $\widetilde{f}_j$ on closed points are given by
\[\widetilde{f}_j\left((z_1:z_2),x\right)= \psi\left(u^{(m(\bfv)+x_j)/2}(z_1e_2+z_2e_1),u^{(m(\bfv)-x_j)/2}(z_1e_1+z_2e_2)\right)(x).\]
If we set $T=z\in\bar\F$ (resp. $T^{-1}=0$), then we pushout the Schubert variety $Y_j$ along the automorphism 
\begin{align*}
& e_1\mapsto u^{(m(\bfv)+x_j)/2}e_1\\
& e_2\mapsto u^{(m(\bfv)-x_j)/2}(ze_1+e_2).
\end{align*}
This is the Schubert variety consisting of the lattices $\latM$ with  $d_1(\latM,Q_j^z)\leq n_j$, where 
\[Q_j^{z}=\langle u^{(m(\bfv)+x_j)/2}e_1,u^{(m(\bfv)-x_j)/2}(ze_1+e_2)\rangle.\]
The conclusion for the point at infinity in $\mathbb{P}^1(\bar\F)$ is similar.
This also shows that the image of $\widetilde{f}_j$ is $Z_j$.
As $\mathbb{P}^1_{\bar\F}\times Y_j$ is reduced, the morphism $\widetilde{f}_j$ factors through $Z_j$ and we obtain a surjective morphism $f_j:\mathbb{P}^1_{\bar\F}\times Y_j\rightarrow Z_j$. As the source of this morphism is projective, it follows that $Z_j$ is a closed irreducible subset of the affine Grassmannian and that the morphism $f_j$ is projective.\\
We have to show that it is birational. Denote by $\widetilde{U}\subset Y_j$ the subset of all lattices 
\[\{\latM=\langle u^{n_j/2}e_1,u^{-n_j/2}(qe_1+e_2)\mid q=\sum_{i=0}^{n_j-1}a_i u^i\rangle\},\]
(our assumptions guarantee that $n_j$ is even in this case).
This subscheme is isomorphic to the affine space $A^{n_j}_{\bar\F}$ and is a maximal dimensional affine subspace of $Y_j$. Now $\widetilde{V}=f_j(V_0\times \widetilde{U})\subset Z_j$ is the subset of all lattices 
\[\left\{\latM=\langle u^{(n+2j+2)/2}e_1,u^{-(n+2j+2)/2}(qe_1+e_2)\rangle\mid q=a_0+\sum_{i=l+(p+1)j}^{n+2j+1}a_i u^i\right\}.\]
This is again an affine space and $f_j$ maps $V_0\times \widetilde{U}$ isomorphically onto $\widetilde{V}$.
Thus it is birational.\\
The case $m(\bfv)\not\equiv x_j\mod 2$ is similar. We have to consider the lattices 
\begin{align*}
&[l+(p+1)j,m(\bfv)-1]_z \hspace{5mm}\text{for}\ z\in\bar\F\\
&[-l-(p+1)j,m(\bfv)-1]_0 
\end{align*}
instead of $Q_j^z$. Now $Y_j\hookrightarrow {\rm{Grass}}\,M_{\bar\F}$ is the inclusion of the Schubert variety of lattices $\latM$ with
$d_2(\latM,\langle b_1,b_2\rangle)=1$ and the same condition on $d_1$ as above. The conclusion is now similar.\\
\noindent \emph{(iii)} For $i\geq 0$ we always have
\[Z_i\not\subset (\bigcup_{j=0}^{i-1} Z_j)\cup Z,\]
because for example $[n+2i+2,m(\bfv)]_0\in Z_i(\bar\F)$ but not in the latter union,
as we can see from the definitions. If $l\neq 2$, then
\[Z\not\subset \bigcup_{j\geq 0} Z_j\]
because for example $[n,m(\bfv)]_u\in Z(\bar\F)$ but not in the latter union.
The claim follows from that and the computation of the dimensions:
At first consider $Z\subset X_0^\bfv$. This is irreducible and its complement has dimension less or equal to $\dim\,Z$.
Now consider $Z\cup Z_0\not\supset Z$. The complement of this subscheme has dimension (strictly) less than $\dim\, Z_0$.
Proceeding by induction on $j$ yields the claim.\\
\noindent \emph{(iv)} In the case $l=2$ we have $Z\subset Z_0$: If $\latM=[x,m(\bfv)]_q\in Z(\bar\F)$ and if we assume again 
$v_u(q)>0$, then $d_1(\latM,\T)\leq n-1$ and hence
\[d_1(\latM,Q_0^0)=d_1(\latM,Q')+d_1(Q',Q_0^0)\leq n-1+l-1=n+2-l.\]
Thus each point $\latM\in Z(\bar\F)$ is also contained in $Z_0$.
The statement now follows by the same argument as in \emph{(iii)}.\\
\noindent \emph{(v)} This is a consequence of \emph{(i)-(iv)}.
\end{proof}
\begin{rem}
On each of the half lines $\calL_z$ for $z\in\mathbb{P}^1(\bar\F)$ we find Schubert varieties with decreasing dimensions.
This behavior is called \emph{"thinning tubes"} in \cite{phimod} 6.d.
compare loc. cit. B 2.
\end{rem}
\subsection{The case $(M_{\bar\F},\Phi)\cong (M_1,\Phi_1)\oplus (M_2,\Phi_2)$}
In this section we assume that 
\[M_{\bar\F}\sim\begin{pmatrix}au^s & 0\\ 0 & bu^t\end{pmatrix}\]
with $a,b\in\bar\F^\times$ and $0\leq s,t<p-1$. Further we assume $s\neq t$ or $a\neq b$.\\
Assume $s=t$ and let $n$ be the largest integer that is congruent to $m(\bfv)\mod 2$ and that satisfies
\begin{equation}\label{n2}
n\leq \tfrac{r_1-r_2}{p+1}.
\end{equation}
Denote by $l$ the smallest integer satisfying 
\begin{equation}\label{l2}
n+2\leq \tfrac{r_1-r_2+2l}{p+1}.
\end{equation}
Define the points 
\[Q_j^{\pm}=[\pm(l+(p+1)j),m(\bfv)]_0,\] 
and the subschemes $Z,Z_j^{\pm}\subset X_0^\bfv$ by:
\begin{align*}
& Z(\bar\F)\ \ =\{\latM\in\bar\B(m(\bfv))\mid d_1(\latM,P_{\rm{red}})\leq n\}\\
& Z_j^\pm(\bar\F)=\{\latM\in\bar\B(m(\bfv))\mid d_1(\latM,Q_j^\pm)\leq n+2-l-(p-1)j\}.
\end{align*}
\begin{prop}\label{propX02} 
Assume $s=t$ and define $n$ and $l$ as in $(\ref{n2})$ and $(\ref{l2})$.\\
\noindent (i) The scheme $Z$ is isomorphic to an $n$-dimensional Schubert variety.\\
\noindent (ii) The schemes $Z_j^\pm$ are isomorphic to $n+2-l-(p-1)j$ dimensional Schubert varieties.\\
\noindent (iii) If $l\neq 1$, then 
\[X_0^\bfv=Z\cup(\bigcup_{j\geq 0}Z_j^+)\cup(\bigcup_{j\geq 0}Z_j^-)\]
is the decomposition of $X^\bfv_0$ into its irreducible components.\\
\noindent (iv) If $l=1$, then
\[X_0^\bfv=(\bigcup_{j\geq 0}Z_j^+)\cup(\bigcup_{j\geq 0}Z_j^-)\]
is the decomposition of $X^\bfv_0$ into its irreducible components.\\
\noindent (v) The dimension of $X^\bfv_0$ is given by 
\[\dim\, X^\bfv_0=\begin{cases} n+1&\text{if}\ l=1\\ n& \text{if}\ l\neq 1.\end{cases}\]
\end{prop}
\begin{proof}
\emph{(i)} and \emph{(ii)} follow immediately from the definitions.\\
As in Lemma $\ref{lemX01}$ we easily find
\begin{equation}
X_0^\bfv=Z\cup(\bigcup_{j\geq 0}Z_j^+)\cup(\bigcup_{j\geq 0}Z_j^-)
\end{equation}
and as in Proposition $\ref{propX01}$ we find
\[Z_i^\pm\not\subset Z\cup(\bigcup_{j\geq 0}^{i-1}Z_j^+)\cup(\bigcup_{j\geq 0}^{i-1}Z_j^-).\]
Further 
\begin{align*}
& Z\not\subset Z_0^\pm \hspace{5mm} \text{if}\ l\neq 1\\
& Z\subset Z_0^\pm \hspace{5mm} \text{if}\ l=1.
\end{align*}
The computations are the same as in the proof of Proposition $\ref{propX01}$
with the only difference that we have to replace $\T$ by $\calL_0\cup\calL_\infty=\A_0\cap \bar\B(m(\bfv))$.
Part \emph{(iii)} and \emph{(iv)} now follow exactly as in the proof of Proposition $\ref{propX01}$.\\
Finally \emph{(v)} follows from \emph{(i)-(iv)}.
\end{proof}
In the case $s\neq t$ we have to distinguish more different cases. We only sketch the structure of the irreducible 
components.\\
Denote by $x_0=\lfloor\tfrac{t-s}{p-1}\rfloor$ the integral part of $\tfrac{t-s}{p-1}$.
Let $n_+$ be the largest integer congruent to $m(\bfv)\mod 2$ such that
\begin{equation}\label{n+red2}
n_+\leq \tfrac{t-s}{p-1}+\tfrac{1}{p+1}(r_1-r_2+2(x_0+1-\tfrac{t-s}{p-1})).
\end{equation}
Let $n_-$ be the smallest integer congruent to $m(\bfv)\mod 2$ such that
\begin{equation}
n_-\geq \tfrac{t-s}{p-1}-\tfrac{1}{p+1}(r_1-r_2+2(\tfrac{t-s}{p-1}-x_0)).
\end{equation}
By Lemma $\ref{lemred2}$, these numbers have the following meaning:
The maximal distance $d_1$ of a $\bfv$-admissible lattice in $\A_q\backslash\A_0$ with $v_u(q)=x_0+1$ from $P_{\rm{red}}$ is $n_+-\tfrac{t-s}{p-1}$;
the maximal distance $d_1$ of a $\bfv$-admissible lattice in $\A_q\backslash\A_0$ with $v_u(q)=x_0$ from $P_{\rm{red}}$ is $\tfrac{t-s}{p-1}-n_-$.
We define 
\begin{equation}
\begin{aligned}
& x_1=\tfrac{1}{2}(n_++n_-)\\
& n\ =\tfrac{1}{2}(n_+-n_-).
\end{aligned}
\end{equation}
Let $Z$ be the subscheme whose closed points are given by
\[Z(\bar\F)=\{\latM\in\bar\B(m(\bfv))\mid d_1(\latM,[x_1,m(\bfv)]_0)\leq n\}.\]
This is a $n$-dimensional Schubert variety.
Let $l_+$ be the smallest integer such that
\[n_++2\leq \tfrac{t-s}{p-1}+\tfrac{1}{p+1}(r_1-r_2+2(l_+-\tfrac{t-s}{p-1})),\]
i.e. the smallest integer such that there are $\bfv$-admissible lattices with $x$-coordinate $n_++2$ in the 
apartments branching of from $\A_0$ at the line $x=l_+$.\\
Similarly, let $l_-$ the largest integer such that 
\[n_--2\geq \tfrac{t-s}{p-1}-\tfrac{1}{p+1}(r_1-r_2+2(\tfrac{t-s}{p-1}-l_-)).\]
For $j\geq 0$ we define the following points 
\begin{align*}
& Q_j^+=[l_++(p+1)j,m(\bfv)]_0\\
& Q_j^-=[l_--(p+1)j,m(\bfv)]_0.
\end{align*}
Again, we define the following subschemes of $X_0^\bfv$:
\begin{equation}\label{Zj-red2}
\begin{aligned}
&Z_j^+(\bar\F)=\{\latM\in\bar\B(m(\bfv))\mid d_1(\latM,Q_j^+)\leq n_++2-l_+-(p-1)j\}\\
&Z_j^-(\bar\F)=\{\latM\in\bar\B(m(\bfv))\mid d_1(\latM,Q_j^-)\leq l_-+2-n_--(p-1)j\}.
\end{aligned}
\end{equation}
These subschemes are isomorphic to Schubert varieties. 
\begin{lem}\label{lemX03} 
With the above notation we have
\begin{equation}
X_0^\bfv(\bar\F)=Z(\bar\F)\cup(\bigcup_{j\geq 0}Z_j^+(\bar\F))\cup(\bigcup_{j\geq 0}Z_j^-(\bar\F)).
\end{equation}
\end{lem}
\begin{proof}
The proof of this fact is similar to the proof of Lemma $\ref{lemX01}$.\\
If $\latM=[x,m(\bfv)]_q$ is a lattice, we denote by $Q'$
the unique point in $\A_0\cap \bar\B(m(\bfv))$ with $d_1(\latM,Q')=d_1(\latM,\A_0)$.
We assume that $Q'=[x',m(\bfv)]_0$ with $x'>\tfrac{t-s}{p-1}$.\\
If $x'<l_+$, then $\latM$ is $\bfv$-admissible (and non-$\bfv$-ordinary) if and only if 
$\latM\in Z(\bar\F)$. This is a direct consequence of the definitions.\\
If $x'\geq l_+$, there is a unique integer $j\geq 0$ such that 
\[l_++(p+1)j\leq x'<l_++(p+1)(j+1).\]
In this case the definitions imply that $\latM$ is $\bfv$-admissible (and non-$\bfv$-ordinary)
if and only if $\latM\in Z_j^+(\bar\F)$
(compare the proof of Lemma $\ref{lemX01}$). \\
The conclusions for $x'<\tfrac{t-s}{p-1}$ are similar. In the set of coordinates considered above,
the computations become more complicated, but we can also deduce this result by interchanging $e_1$ and $e_2$.
\end{proof}
\begin{prop}
With the notations of $(\ref{n+red2})$-$(\ref{Zj-red2})$ :
\begin{equation}\label{decomp2}
X_0^\bfv=Z\cup(\bigcup_{j\geq 0}Z_j^+)\cup(\bigcup_{j\geq 0}Z_j^-).
\end{equation}
\noindent (i) The scheme $Z$ is isomorphic to an $n$-dimensional Schubert variety.\\
\noindent (ii) For $j\geq 0$ the schemes $Z_j^\pm$ are isomorphic to Schubert varieties of dimension 
\begin{align*}
& \dim\, Z_j^+=n_++2-l_+-(p-1)j\\
& \dim\, Z_j^-=l_-+2-n_--(p-1)j.
\end{align*}
\noindent (iii) If $Z\not\subset Z_0^+\cup Z_0^-$, then $(\ref{decomp2})$ is the decomposition of $X_0^\bfv$ into 
its irreducible components.\\
\noindent (iv) If $Z\subset Z_0^+\cup Z_0^-$, then 
\[X_0^\bfv=(\bigcup_{j\geq 0}Z_j^+)\cup(\bigcup_{j\geq 0}Z_j^-)\]
is the decomposition of $X_0^\bfv$ into its irreducible components. 
\end{prop}
In this case we would have to consider many cases in order to determine whether $Z\subset Z_0^+\cup Z_0^-$ or not from the
given integers $r_1,r_2,s,t$.  
\begin{proof}
By use of Lemma $\ref{lemX03}$, this is nearly the same as in the proof of Proposition $\ref{propX02}$.
\end{proof}
\begin{rem}
Again, we find that the Schubert varieties of decreasing dimension defined above correspond to the \emph{"thinning tubes"} in \cite{phimod} 6.d
along the $\Phi$-stable half lines $\{[x,m(\bfv)]_0\mid x\leq \tfrac{t-s}{p-1}\}$ and 
$\{[x,m(\bfv)]_0\mid x\geq \tfrac{t-s}{p-1}\}$ in the building for $PGL_2(\bar\F((u)))$.
\end{rem}
\subsection{The case of a non split extension.}
As in section 4.3, we assume that there is a basis $e_1,e_2$ of $M_{\bar\F}$ such that
\[M_{\bar\F}\sim\begin{pmatrix}au^s & \gamma \\ 0& bu^t\end{pmatrix}\]
for some $a,b\in\bar\F^\times,\,\gamma\in\bar\F((u))$ and $0\leq s,t<p-1$. We assume that $(M_{\bar\F},\Phi)$
is a non-split extension and use the notations of section 4.3.\\
Denote by $x_0$ the largest integer $x_0<\tfrac{k-s}{p}$. 
Let $n_+$ be the largest integer congruent $m(\bfv)\mod 2$ such that
\begin{equation}\label{n+3}
n_+\leq \tfrac{1}{p+1}(r_1-r_2-s-t+2k).
\end{equation} 
Let $n_-$ be the smallest integer congruent $m(\bfv)\mod 2$ such that
\begin{equation}\label{n-3}
n_-\geq \tfrac{t-s}{p-1}-\tfrac{1}{p+1}(r_1-r_2+2(\tfrac{t-s}{p-1}-x_0)).
\end{equation}
These numbers have the following meaning:
The integer $n_+$ is the maximal $x$-coordinate of a $\bfv$-admissible lattice in $\A_q$ with $v_u(q)\geq \tfrac{k-s}{p}$;
further $\tfrac{t-s}{p-1}-n_-$ is the maximal distance $d_1$ of a $\bfv$-admissible lattice in $\A_q\backslash \A_0$ with $v_u(q)=x_0$
from $P_{\rm{red}}$.
As above, we define the following integers 
\begin{equation}
\begin{aligned}
& x_1=\tfrac{1}{2}(n_++n_-)\\
& n\ =\tfrac{1}{2}(n_+-n_-).
\end{aligned}
\end{equation}
We have $x_1\in\{x_0,x_0+1\}$ which can be deduced from the equation
\[\tfrac{1}{p+1}(r_1-r_2-s-t+2k)-\tfrac{k-s}{p}
=\tfrac{k-s}{p}-(\tfrac{t-s}{p-1}-\tfrac{1}{p+1}(r_1-r_2+2(\tfrac{t-s}{p-1}-\tfrac{k-s}{p}))).\]
(Here, we compute the distance from the point $\tfrac{k-s}{p}$ and write $\tfrac{k-s}{p}$ instead of $x_0$
as in $(\ref{n-3})$). We define the following subset 
\begin{equation}
Z(\bar\F)=\{\latM\in\bar\B(m(\bfv))\mid d_1(\latM,[x_1,m(\bfv)]_0)\leq n\}.
\end{equation}
Let $l_-$ be the largest integer such that
\begin{equation}
n_--2\geq \tfrac{t-s}{p-1}-\tfrac{1}{p+1}(r_1-r_2+2(\tfrac{t-s}{p-1}-l_-)).
\end{equation}
For $j\geq 0$ we define the points $Q_j^-=[l_--(p+1)j,m(\bfv)]_0$ and the subsets
\begin{equation}\label{Zj-nonsplit}
Z_j^-(\bar\F)=\{\latM\in\bar\B(m(\bfv))\mid d_1(\latM,Q_j^-)\leq l_-+2-n_--(p-1)j\}.
\end{equation}
\begin{lem}
With the above notations
\[X_0^\bfv(\bar\F)=Z(\bar\F)\cup (\bigcup_{j\geq 0}Z_j^-(\bar\F)).\]
\end{lem}
\begin{proof}
Let $\latM=[x,m(\bfv)]_q$ be a lattice.\\
If $v_u(q)\geq \tfrac{k-s}{p}$ (or equivalently if $v_u(q)>x_0$), then (by Lemma $\ref{lemred3}$) 
$\latM$ is $\bfv$-admissible (and non-$\bfv$-ordinary) iff 
\[x\leq \tfrac{1}{p+1}(r_1-r_2-s-t+2k),\]
or equivalently iff $x\leq n_+$. \\
If $v_u(q)=x_0$ and $x>v_u(q)$, then (by Lemma $\ref{lemred3}$) $\latM$ is $\bfv$-admissible and non-$\bfv$-ordinary iff
\[d_1(\latM,[\tfrac{t-s}{p-1},m(\bfv)]_0)\leq \tfrac{1}{p+1}(r_1-r_2+2(\tfrac{t-s}{p-1}-x_0)),\]
or equivalently iff $d_1(\latM,[\tfrac{t-s}{p-1},m(\bfv)]_0)\leq \tfrac{t-s}{p-1}-n_-$.
By the definitions of $x_1$ and $n$ and the fact $x_1\in\{x_0,x_0+1\}$, we find that in both cases $\latM$ is
$\bfv$-admissible (and non $\bfv$-ordinary) iff $\latM\in Z(\bar\F)$.\\
For the $\bfv$-admissible lattices $\latM=[x,m(\bfv)]_q$ with $v_u(q)<x_0$ we proceed as in the proof of Lemma $\ref{lemX03}$.
\end{proof}
\begin{prop}
With the notations of $(\ref{n+3})$-$(\ref{Zj-nonsplit})$:
\begin{equation}\label{decomp100}
X_0^\bfv=Z\cup(\bigcup_{j\geq 0} Z_j^-).
\end{equation}
\noindent (i) The scheme $Z$ is isomorphic to an $n$-dimensional Schubert variety.\\
\noindent (ii) For $j\geq 0$ the schemes $Z_j^-$ are isomorphic to Schubert varieties of dimension
\[\dim\, Z_j^-=l_-+2-n_--(p-1)j.\]
\noindent (iii) If $Z\not\subset Z_0^-$, then $(\ref{decomp100})$ is the decomposition of $X_0^\bfv$ into
its irreducible components.\\
\noindent (iv) If $Z\subset Z_0^-$, then the decomposition of $X_0^\bfv$ into its irreducible components is given by
\[X_0^\bfv=\bigcup_{j\geq 0}Z_j^-.\]
\end{prop}
\begin{proof}
\noindent \emph{(i),(ii)} This follows from the definitions.\\
\noindent \emph{(iii),(iv)} As in the discussion of the other cases, our Schubert varieties are constructed in a way such that
\[Z_i^-\not\subset Z\cup(\bigcup_{j=0}^{i-1}Z_j^-)\]
for all $i\geq 0$. The claim follows from this and the computation of the dimension 
(compare the proof of Proposition $\ref{propX01}$).
\end{proof}
\begin{rem}
We find a sequence of Schubert varieties of decreasing dimension along the unique $\Phi$-stable
half line $\{[x,m(\bfv)]_0\mid x<\tfrac{k-s}{p}\}$ corresponding to the \emph{"thinning tubes"} in \cite{phimod} 6.d.\\
Further, we find a Schubert variety $Z$ corresponding to a ball with given radius around a given point 
as in loc. cit. A 3 resp. B 2.
\end{rem}
The discussion of this section implies the following result.
\begin{theo}
If $(M_{\bar\F},\Phi)$ is not isomorphic to the direct sum of two isomorphic one-dimensional $\phi$-modules, then the irreducible components of $X_0^\bfv$ are Schubert varieties.
Especially they are normal.
\end{theo}
\section{Relation to Raynaud's theorem}
In this section, we assume $p\neq 2$.
In \cite{Raynaud}, Raynaud introduces a partial order on the set of finite flat models for $V_\F$
(i.e. the set of $\F$-valued points of $\GR$) by defining $\mathcal{G}_1\preceq\mathcal{G}_2$ if there exists 
a morphism $\mathcal{G}_2\rightarrow \mathcal{G}_1$ inducing the identity on the generic fiber of ${\rm{Spec}}\,\mathcal{O}_K$.
By (loc. cit. 2.2.3 and 3.3.2), this order admits a minimal and maximal object (if the set is non-empty) which agree 
if $e<p-1$.\\
In our case, Raynaud's partial order is given by the inclusion of lattices in $M_\F$:
Inclusion of two admissible lattices is a morphism that commutes with the semi-linear map $\Phi$ and induces the identity
of $M_\F$ after inverting $u$. Here, a lattice $\latM$ is called {\emph {admissible}} if it defines a finite flat group scheme, i.e. if
$u^e\latM\subset \langle\Phi(\latM)\rangle\subset \latM$. 
\begin{prop}
Let $\rho:G_K\rightarrow V_\F$ be a continuous representation of $G_K$. Assume that there exists a finite extension
$\F'$ of $\F$ such that there is a finite flat group scheme model over ${\rm{Spec}}\,\OK$ for the induced $G_K$ representation
on $V_{\F'}=V_\F\otimes_\F\F'$. Then there exists a finite flat model for $V_\F$, i.e.
\[\GR\neq \emptyset\Rightarrow \GR(\F)\neq \emptyset.\] 
\end{prop}
\begin{proof}
Our assumptions imply $\GR(\F')\neq \emptyset$ for some finite extension $\F'$ of $\F$. 
Hence, by Raynaud's theorem, the set $\GR(\F')$ has a unique 
minimal element. The natural action of the Galois group $\rm{Gal}(\F'/\F)$ on $\GRloc(\F')$ preserves the partial order
(it preserves inclusion of lattices) and hence the minimal element is stable under this action. 
Consequently, the minimal object is already defined over $\F$. 
\end{proof} 
We now want to reprove Raynaud's theorem in our context: we will show that there is a minimal and a maximal lattice
for the order induced by inclusion on the set $\GR(\bar\F)$.
\begin{prop}
There exists a minimal and a maximal admissible lattice $\latM_{\min}$ resp. $\latM_{\max}$ for the order defined by the inclusion.
\end{prop}
\begin{proof}
We only prove the statement about the maximal lattice. The other one is analogue. We choose a basis $e_1,e_2$.
The proposition follows from the following two observations:\\
\emph{\noindent (a) There exists a unique admissible lattice with minimal $y$-coordinate.\\
\noindent (b) If $\latM$ is a admissible lattice with non-minimal $y$-coordinate,
then it is contained in an admissible lattice with strictly smaller $y$-coordinate.}\\
{\emph{Proof of (a)}:} First it is clear that the $y$-coordinates of admissible lattices are bounded below: 
If 
\[\langle e_1,e_2\rangle\sim A'=\begin{pmatrix}\alpha' &\beta' \\ \gamma' & \delta' \end{pmatrix}\]
and if $\latM=[x,y]_q$ is admissible, then $2e-d'=(p-1)y+v_u(\det\,A')$ with 
\[0\leq d'=\dim \Phiquot\leq 2e.\]
Suppose now that $\latM_1$ and $\latM_2$ are admissible lattices with the same $y$-coordinate. There is a basis $b_1,b_2$ such that
\begin{align*}
& \latM_1=\langle b_1,b_2\rangle\\
& \latM_2=\langle u^nb_1,u^{-n}b_2\rangle.
\end{align*}
for some $n\geq 0$. We have 
\begin{align*}
& \latM_1\sim A=\begin{pmatrix}\alpha & \beta \\ \gamma & \delta \end{pmatrix}\\
& \latM_2\sim B=\begin{pmatrix}u^{n(p-1)}\alpha & u^{-n(p+1)}\beta \\ u^{n(p+1)}\gamma & u^{-n(p-1)}\delta \end{pmatrix}
\end{align*}
for some $\alpha,\beta,\gamma,\delta\in\bar\F((u))$. Define 
\[\latM_3=\langle b_1,u^{-n}b_2\rangle=\begin{pmatrix}1&0\\0&u^{-n}\end{pmatrix}\latM_1.\]
Then $\latM_3$ has strictly smaller $y$-coordinate and is admissible. Indeed:
\[\latM_3\sim C=\begin{pmatrix}\alpha & u^{-np}\beta\\ u^{n}\gamma & u^{-n(p-1)}\delta \end{pmatrix}\]
and we have to show: $\min_{i,j} c_{ij}\geq 0$ and $v_u(\det\,C)-\min_{i,j} c_{ij}\leq e$. \\
But because $\latM_1$ and $\latM_2$
are admissible we know:
\begin{align*}
& v_u(\alpha)\geq 0 && v_u(u^{-np}\beta)\geq n\\
& v_u(u^{n}\gamma)\geq n && v_u(u^{-n(p-1)}\delta)\geq 0.
\end{align*} 
Similarly $v_u(\det\,C)=v_u(\det\,A)-(p-1)n=v_u(\det\,B)-(p-1)n$ and hence:
\begin{align*}
& v_u(\det\,C)-v_u(\alpha)\leq e-(p-1)n && v_u(\det\,C)-v_u(u^{-np}\beta)\leq e-pn\\
& v_u(\det\,C)-v_u(u^n\gamma)\leq e-pn && v_u(\det\,C)-v_u(u^{-n(p-1)}\delta)\leq e-(p-1)n.
\end{align*}
\emph{Proof of (b)}: Let $\latM$ be an admissible lattice with non-minimal $y$-coordinate.
Then there exists an admissible lattice $\latM'$ with strictly smaller $y$-coordinate.
These lattices are contained in a common apartment and hence there is a basis $b_1,b_2$ such that
\begin{align*}
& \latM=\langle b_1,b_2\rangle\\
& \latM'=\langle u^mb_1,u^nb_2\rangle
\end{align*}
for some integers $m,n$ with $m+n<0$, because the $y$-coordinate of $\latM'$ is strictly smaller than the $y$-coordinate
of $\latM$. \\
Without loss of generality, we assume $m-n\geq 0$. If $m\leq 0$, then $n\leq 0$ and we are done, since then $\latM\subset\latM'$.\\
If $m>0$ we claim that the lattice $\latM_1=\langle b_1, u^{m+n}b_2\rangle$ is admissible. Indeed
\begin{align*}
&\latM\sim\begin{pmatrix}\alpha &\beta\\ \gamma & \delta\end{pmatrix}\\
&\latM'\sim\begin{pmatrix}u^{(p-1)m}\alpha &u^{pn+m}\beta\\ u^{pm-n}\gamma & u^{(p-1)n}\delta\end{pmatrix}\\
&\latM_1\sim\begin{pmatrix}\alpha & u^{p(m+n)}\beta\\ u^{-m-n}\gamma & u^{(p-1)(m+n)}\delta\end{pmatrix} 
\end{align*}
and the claim follows by a similar argument as in the proof of \emph{(a)}.
\end{proof}
\begin{prop}
If $e<p-1$, then the minimal and the maximal lattice coincide.
\end{prop}
\begin{proof}
Denote the minimal lattice by $\latM_{\min}$, the maximal by $\latM_{\max}$. There is a apartment containing both lattices
and we may assume $\latM_{\max}=\langle e_1,e_2 \rangle=[0,0]_0$ and $\latM_{\min}=[x,y]_0$ for some $y\geq 0$.\\
Let $A\in GL_2(\bar\F((u)))$ be a matrix such that 
\[\latM_{\max}\sim A=\begin{pmatrix}\alpha & \beta\\ \gamma &\delta \end{pmatrix}.\] 
Define $d_{\min}=\dim\ \latM_{\min}/\langle\Phi(\latM_{\min})\rangle$ and similarly $d_{\max}$. Then 
\begin{align*}
& 2e-d_{\max}=2e-\dim\, \latM_{\max}/\langle\Phi(\latM_{\max})\rangle=v_u(\det\, A)\\
& 2e-d_{\min}=2e-\dim\, \latM_{\min}/\langle\Phi(\latM_{\min})\rangle=v_u(\det\, A)+(p-1)y.
\end{align*}
Thus we have $(p-1)y=d_{\min}-d_{\max}\leq 2e<2(p-1)$ and hence $y=0$ or $y=1$.
If $y=0$ we are done, as $\latM_{\max}$ is the unique lattice with minimal $y$-coordinate.\\
Assume $y=1$. In this case $\latM_{\min}\subset \latM_{\max}$ implies $\latM_{\min}=[\pm 1,1]_0$.
Without loss of generality, we assume $\latM_{\min}=[-1,1]_0=\langle u e_1,e_2\rangle$. Then
\[\latM_{\min}\sim B=\begin{pmatrix} u^{p-1}\alpha & u^{-1}\beta \\ u^p\gamma & \delta\end{pmatrix}.\]
As both lattices are admissible, we have 
\[\max\{v_u(\det\, B)-v_u(u^{-1}\beta),v_u(\det\, B)-v_u(\delta)\}\leq e\]
and hence:
\begin{align*}
& v_u(\alpha)\geq 0 && v_u(\beta)-1\geq v_u(\det\, B)-e=v_u(\det\, A)+(p-1)-e>v_u(\det\, A)\\
& v_u(\gamma)\geq 0 && v_u(\delta)\geq v_u(\det\, B)-e=v_u(\det\, A)+(p-1)-e>v_u(\det\, A)
\end{align*}
It follows that
\[v_u(\det\, A)=v_u(\alpha\delta-\beta\gamma)\geq \min\{v_u(\alpha)+v_u(\delta),v_u(\beta)+v_u(\gamma)\}>v_u(\det\, A).\]
Contradiction.
\end{proof}
Finally we want to determine the elementary divisors of $\langle\Phi(\latM)\rangle$ with respect to $\latM$
for the minimal and the maximal lattice in the cases where $(M_{\bar\F},\Phi)$ is simple resp. split reducible.
If $(M_{\bar\F},\Phi)$ is non-split reducible, the computations turn out to be very difficult
and are omitted.
\subsection{The absolutely simple case} As in section $3$, we fix a basis $e_1,e_2$ such that
\[M_{\bar\F}\sim \begin{pmatrix} 0& au^s\\ 1 & 0\end{pmatrix}\]
with $a\in\bar\F^\times$ and $0\leq s<p^2-1$. Let $\latM_{\min}$ be the minimal and $\latM_{\max}$ be the maximal lattice.
Denote by $s_1,s_2$ the unique integers $0\leq s_1<p+1$ resp. $0\leq s_2<p-1$ with 
\begin{align*}
& s_1\equiv s\mod (p+1) \\
& s_2\equiv s\mod (p-1).
\end{align*}
Because $p-1$ and $p+1$ are both even, we find $s_1\equiv s_2\mod 2$.\\
Similarly let $s_2'$ be the unique integer 
$0\leq s_2'<p-1$ with $2e-s\equiv s_2'\mod (p-1)$.
\begin{prop}
Denote by $m=\tfrac{s-s_1}{p+1}$ the integral part of $\tfrac{s}{p+1}$ and by $l=\tfrac{s-s_2}{p-1}$
resp. $l'=\tfrac{2e-s-s_2'}{p-1}$ the integral part of $\tfrac{s}{p-1}$ resp. $\tfrac{2e-s}{p-1}$.\\
\noindent (i)  The elementary divisors $(a_{\max},b_{\max})$ of $\langle\Phi(\latM_{\max})\rangle$ with respect to 
$\latM_{\max}$ are given by
\[
\begin{cases}
(\tfrac{s_1+s_2}{2},\tfrac{s_2-s_1}{2}) &\text{if}\ \ l+m\in 2\mathbb{Z}\ \text{and}\ s_2\geq s_1\\
(\tfrac{s_2-s_1}{2}+p,\tfrac{s_1+s_2}{2}-1) &\text{if}\ \ l+m\in 2\mathbb{Z}\ \text{and}\ s_2< s_1\\
(\tfrac{s_2-s_1+(p+1)}{2},\tfrac{s_1+s_2-(p-1)}{2}) &\text{if}\ \ l+m\notin 2\mathbb{Z}\ \text{and}\ s_1+s_2\geq p+1\\
(\tfrac{s_1+s_2+(p-1)}{2},\tfrac{s_2-s_1+(p-1)}{2}) &\text{if}\ \ l+m\notin 2\mathbb{Z}\ \text{and}\ s_1+s_2< p+1.
\end{cases}
\]
\noindent (ii) The elementary divisors $(a_{\min},b_{\min})$ of $\langle\Phi(\latM_{\min})\rangle$ with respect to 
$\latM_{\min}$ are given by
\[
\begin{cases}
(e+\tfrac{s_1-s_2'}{2},e-\tfrac{s_1+s_2}{2}) &\text{if}\ \ l'+m\in 2\mathbb{Z}\ \text{and}\ s_1\leq s_2'\\
(e+1-\tfrac{s_1+s_2}{2},e-p+\tfrac{s_1-s_2}{2}) &\text{if}\ \ l'+m\in 2\mathbb{Z}\ \text{and}\ s_1>s_2'\\
(e+\tfrac{(p+1)-s_1-s_2'}{2},e+\tfrac{s_1-s_2'-(p+1)}{2}) &\text{if}\ \ l'+m\notin 2\mathbb{Z}\ \text{and}\ s_1+s_2'\geq p+1\\
(e+\tfrac{s_1-s_2'-(p-1)}{2},e+\tfrac{-s_1-s_2'-(p-1)}{2}) &\text{if}\ \ l'+m\notin 2\mathbb{Z}\ \text{and}\ s_1+s_2<p+1.
\end{cases}
\]
\end{prop}
\begin{proof}
We only prove the statement on $\latM_{\max}$. From Corollary $\ref{singptirred}$ we know that
the lattice $\latM_{\max}$ is contained in the apartment defined by $e_1,e_2$. 
For a lattice $\latM$, denote the elementary divisors of $\langle\Phi(\latM)\rangle$
with respect to $\latM$ by $(a,b)$.
Then $\latM$ is admissible if $0\leq b,a\leq e$. As we are assuming that there exist admissible lattices
we only have to check the condition $a,b\geq 0$ for $\latM_{\max}$ (and the condition $a,b\leq e$ for $\latM_{\min}$).\\
If $l+m\in 2\mathbb{Z}$, the candidate for the maximal lattice is $[m,-l]_0$.
If this is not admissible, then we take $[m+1,-l+1]_0$.
Computing the elementary divisors $a,b$ by use of Lemma $\ref{lemirred}$ and Definition $\ref{defd1d2}$
we find the above expressions.
In the case $l+m\notin 2\mathbb{Z}$ we deal with the lattices $[m+1,-l]_0$ and $[m,-l+1]_0$.
\end{proof}
\subsection{The split reducible case}
As in section $4$, we fix a basis $e_1,e_2$ such that
\[M_{\bar\F}\sim\begin{pmatrix}au^s & 0 \\ 0& bu^t\end{pmatrix}\]
with $a,b\in\bar\F^\times$ and $0\leq s,t<p-1$. 
\begin{prop}\noindent (i) The elementary divisors $(a_{\max},b_{\max})$ of $\langle\Phi(\latM_{\max})\rangle$ with respect to 
$\latM_{\max}$ are given by 
\[(a_{\max},b_{\max})=\begin{cases}(t,s)& \text{if}\ t\geq s\\ (s,t) & \text{if}\ s\geq t .\end{cases}\]
\noindent (ii) The elementary divisors $(a_{\min},b_{\min})$ of $\langle\Phi(\latM_{\min})\rangle$ with respect to 
$\latM_{\min}$ are given by 
\[(a_{\min},b_{\min})=
\begin{cases}((p-1)\lfloor\tfrac{e-t}{p-1}\rfloor+t,(p-1)\lfloor\tfrac{e-s}{p-1}\rfloor+s) 
& \text{if}\ \tfrac{t-s}{p-1}\geq \lfloor\tfrac{e-s}{p-1}\rfloor-\lfloor\tfrac{e-t}{p-1}\rfloor \\
((p-1)\lfloor\tfrac{e-s}{p-1}\rfloor+s,(p-1)\lfloor\tfrac{e-t}{p-1}\rfloor+t)
& \text{if}\ \tfrac{t-s}{p-1}\leq \lfloor\tfrac{e-s}{p-1}\rfloor-\lfloor\tfrac{e-t}{p-1}\rfloor.
\end{cases}\]
\end{prop}
\begin{proof}
From Corollary $\ref{oneptred1}$ and Corollary $\ref{oneptred2}$ 
we know that the minimal and the maximal lattice are contained in the apartment 
defined by $e_1,e_2$. Now
\begin{align*}
&\Phi(u^me_1)=au^{(p-1)m+s}(u^me_1)\\
&\Phi(u^ne_2)=bu^{(p-1)n+t}(u^ne_2).
\end{align*}
The first part of the Proposition follows from the fact that $s$ (resp. $t$) are the smallest positive integers
that are congruent to $s$ (resp. $t$) modulo $p-1$.\\
The second part follows from the fact that $(p-1)\lfloor\tfrac{e-s}{p-1}\rfloor+s$ is the largest integer smaller than $e$
that is congruent to $s$ modulo $p-1$ (and similar for $t$).
\end{proof}

\vspace{0.5cm}

Eugen Hellmann\\
Mathematisches Institut der Universit\"at Bonn\\
Beringstr. 1\\
D-53115 Bonn\\
E-mail: hellmann@math.uni-bonn.de

\end{document}